\documentclass[preprint,11pt]{elsarticle}

\usepackage{amsmath,amssymb,bm,subfigure}
\usepackage{accents}

\biboptions{sort&compress}
\usepackage[usenames,dvipsnames,svgnames,table]{xcolor}

\usepackage{color}

\usepackage{tikz}
\usetikzlibrary{positioning}
\usetikzlibrary{decorations.pathreplacing}
\usetikzlibrary{matrix}

\usepackage{graphicx}
\usepackage{calc,fancyvrb}
\usepackage[margin=1.in]{geometry}

\usepackage{algorithm,algpseudocode}
\usepackage{algorithmicx}
\algrenewcommand\alglinenumber[1]{\footnotesize #1:} 
\newcommand{\algFontSize}{\small}

\newcommand{\bibtex}{ ${\mathrm {B{\scriptstyle {IB}}\!T\!_{\displaystyle E}\!X} }$}

\usepackage[bookmarks=true,colorlinks=true,linkcolor=blue]{hyperref}

\newcommand{\blue}{\color{blue}}

\newlength{\tfwidth}
\newlength{\tfheight}
\newlength{\tfxa}
\newlength{\tfxb}
\newlength{\tfya}
\newlength{\tfyb}
\newcommand{\trimFigWithBox}[6]{%
\setlength\fboxsep{0pt}%
\setlength\fboxrule{1.0pt}
\fbox{\includegraphics[width=#2, clip, trim=#3 #4 #5 #6]{#1}}%
}
\newcommand{\trimFigNoBox}[6]{%
\setlength\fboxsep{1pt}
\setlength\fboxrule{0.0pt}
\fbox{\includegraphics[width=#2, clip, trim=#3 #4 #5 #6]{#1}}%
}
\newcommand{\trimFigHeightWithBox}[6]{%
\setlength\fboxsep{0pt}%
\setlength\fboxrule{1.0pt}
\fbox{\includegraphics[height=#2, clip, trim=#3 #4 #5 #6]{#1}}%
}
\newcommand{\trimFigHeightNoBox}[6]{%
\setlength\fboxsep{1pt}
\setlength\fboxrule{0.0pt}
\fbox{\includegraphics[height=#2, clip, trim=#3 #4 #5 #6]{#1}}%
}

\newsavebox\figBox

\newcommand{\trimw}[6]{%
\sbox\figBox{\includegraphics{#1}}
\setlength{\tfwidth}{\the\wd\figBox}
\setlength{\tfheight}{\the\ht\figBox}
\setlength{\tfxa}{\tfwidth*\real{#3}}%
\setlength{\tfxb}{\tfwidth*\real{#4}}%
\setlength{\tfya}{\tfheight*\real{#5}}%
\setlength{\tfyb}{\tfheight*\real{#6}}%
\trimFigNoBox{#1}{#2}{\tfxa}{\tfya}{\tfxb}{\tfyb}%
}

\newcommand{\trimwb}[6]{%

\sbox\figBox{\includegraphics{#1}}
\setlength{\tfwidth}{\the\wd\figBox}
\setlength{\tfheight}{\the\ht\figBox}
\setlength{\tfxa}{\tfwidth*\real{#3}}%
\setlength{\tfxb}{\tfwidth*\real{#4}}%
\setlength{\tfya}{\tfheight*\real{#5}}%
\setlength{\tfyb}{\tfheight*\real{#6}}%
\trimFigWithBox{#1}{#2}{\tfxa}{\tfya}{\tfxb}{\tfyb}%
}

\newcommand{\trimh}[6]{%
\sbox\figBox{\includegraphics{#1}}
\setlength{\tfwidth}{\the\wd\figBox}
\setlength{\tfheight}{\the\ht\figBox}
\setlength{\tfxa}{\tfwidth*\real{#3}}%
\setlength{\tfxb}{\tfwidth*\real{#4}}%
\setlength{\tfya}{\tfheight*\real{#5}}%
\setlength{\tfyb}{\tfheight*\real{#6}}%
\trimFigHeightNoBox{#1}{#2}{\tfxa}{\tfya}{\tfxb}{\tfyb}%
}

\newcommand{\trimhb}[6]{%

\sbox\figBox{\includegraphics{#1}}
\setlength{\tfwidth}{\the\wd\figBox}
\setlength{\tfheight}{\the\ht\figBox}
\setlength{\tfxa}{\tfwidth*\real{#3}}%
\setlength{\tfxb}{\tfwidth*\real{#4}}%
\setlength{\tfya}{\tfheight*\real{#5}}%
\setlength{\tfyb}{\tfheight*\real{#6}}%
\trimFigHeightWithBox{#1}{#2}{\tfxa}{\tfya}{\tfxb}{\tfyb}%
}
\usepackage{xargs}

\newcommandx{\figByHeight}[9][5=0, 6=0, 7=0, 8=0,9=]{
\draw (#1,#2) node[anchor=south west,xshift=-16pt,yshift=-4pt] {\trimh{#3}{#4}{#5}{#6}{#7}{#8}};}
\newcommandx{\figByHeightb}[9][5=0, 6=0, 7=0, 8=0,9=]{
\draw (#1,#2) node[anchor=south west,xshift=-16pt,yshift=-4pt] {\trimhb{#3}{#4}{#5}{#6}{#7}{#8}};}

\newcommandx{\figByHeightWithLabel}[9][5=0, 6=0, 7=0, 8=0,9=]{
\draw (#1,#2) node[anchor=south west,xshift=-16pt,yshift=-4pt] {\trimh{#3}{#4}{#5}{#6}{#7}{#8}} node[draw=white,fill=white,inner sep=1pt,anchor=south west] {#9};}
\newcommandx{\figByHeightWithLabelb}[9][5=0, 6=0, 7=0, 8=0,9=]{
\draw (#1,#2) node[anchor=south west,xshift=-16pt,yshift=-4pt] {\trimhb{#3}{#4}{#5}{#6}{#7}{#8}} node[draw=white,fill=white,inner sep=1pt,anchor=south west] {#9};}

\newcommandx{\figByWidth}[9][5=0, 6=0, 7=0, 8=0,9=]{
\draw (#1,#2) node[anchor=south west,xshift=-16pt,yshift=-4pt] {\trimw{#3}{#4}{#5}{#6}{#7}{#8}};}
\newcommandx{\figByWidthb}[9][5=0, 6=0, 7=0, 8=0,9=]{
\draw (#1,#2) node[anchor=south west,xshift=-16pt,yshift=-4pt] {\trimwb{#3}{#4}{#5}{#6}{#7}{#8}};}

\newcommandx{\figByWidthWithLabel}[9][5=0, 6=0, 7=0, 8=0,9=]{
\draw (#1,#2) node[anchor=south west,xshift=-16pt,yshift=-4pt] {\trimw{#3}{#4}{#5}{#6}{#7}{#8}} node[draw=white,fill=white,inner sep=1pt,anchor=south west] {#9};}
\newcommandx{\figByWidthWithLabelb}[9][5=0, 6=0, 7=0, 8=0,9=]{
\draw (#1,#2) node[anchor=south west,xshift=-16pt,yshift=-4pt] {\trimwb{#3}{#4}{#5}{#6}{#7}{#8}} node[draw=white,fill=white,inner sep=1pt,anchor=south west] {#9};}

\newcommand{\plotTwoFigsByWidth}[5]{%
\begin{figure}[htb]
\begin{center}
\begin{tikzpicture}[scale=1]
  \useasboundingbox (0,.25) rectangle (16,#5);  

   \draw(0.0,0.0) node[anchor=south west,xshift=-16pt,yshift=-5pt] {\trimw{#1}{#5}{.0}{.0}{.0}{.0}};
   \draw(8.0,0.0) node[anchor=south west,xshift=-16pt,yshift=-5pt] {\trimw{#2}{#5}{.0}{.0}{.0}{.0}};
\end{tikzpicture}
\end{center}
\caption{#3}
\label{#4}
\end{figure}
}

\renewcommand{\url}[1]{}
\newcommand{\citeCount}[1]{}

\newcommand{\bogus}[1]{{}}

\newcommand{\bni}{\bigskip\noindent}
\newcommand{\mni}{\medskip\noindent}

\newlength{\ycbTop}
\newlength{\ycbMid}%

\newcommand{\p}{\partial}

\newcommand{\f}[2]{\frac{#1}{#2}}

\def\ba#1\ea{\begin{align}#1\end{align}}

\def\bas#1\eas{\begin{align*}#1\end{align*}}

\def\bat#1\eat{\begin{alignat}{3}#1\end{alignat}}

\def\bats#1\eats{\begin{alignat*}{3}#1\end{alignat*}}

\newcommand{\bse}{\begin{subequations}}
\newcommand{\ese}{\end{subequations}}

\newcommand{\Dzt}{D_{0t}}
\newcommand{\Dpt}{D_{+t}}
\newcommand{\Dmt}{D_{-t}}
\newcommand{\dt}{\Delta t}

\newcommand{\defeq}{\overset{{\rm def}}{=}}
\newcommand{\eqdef}{\overset{{\rm def}}{=}}

\newcommand{\av}{\mathbf{ a}}
\newcommand{\bv}{\mathbf{ b}}

\newcommand{\fv}{\mathbf{ f}}
\newcommand{\gv}{\mathbf{ g}}

\newcommand{\iv}{\mathbf{ i}}
\newcommand{\jv}{\mathbf{ j}}

\newcommand{\nv}{\mathbf{ n}}

\newcommand{\rv}{\mathbf{ r}}

\newcommand{\xv}{\mathbf{ x}}

\newcommand{\Ev}{\mathbf{ E}}

\newcommand{\Gv}{\mathbf{ G}}
\newcommand{\Hv}{\mathbf{ H}}
\newcommand{\Iv}{\mathbf{ I}}

\newcommand{\Nv}{\mathbf{ N}}

\newcommand{\Pv}{\mathbf{ P}}

\newcommand{\Wv}{\mathbf{ W}}

\newcommand{\half}{{1\over2}}

\newcommand{\Real}{{\mathbb R}}

\newcommand{\Bc}{{\mathcal B}}

\newcommand{\Ec}{{\mathcal E}}

\newcommand{\Gc}{{\mathcal G}}

\newcommand{\Kc}{{\mathcal K}}

\newcommand{\Nc}{{\mathcal N}}
\newcommand{\Oc}{{\mathcal O}}

\newcommand{\Tc}{{\mathcal T}}

\newcommand{\alphav}{\boldsymbol{\alpha}}
\newcommand{\betav}{\boldsymbol{\beta}}

\newcommand{\grad}{\nabla}

\newcommand{\nd}{n_d}

\newcommand{\eps}{\epsilon}

\newcommand{\predE}{\Ev^{n+1,*}_{\jv}}
\newcommand{\predP}{\Pv^{n+1,*}_{m,\jv}}
\newcommand{\predN}{N^{n+1,*}_{\ell,\jv}}

\DeclareMathOperator{\sech}{sech}

\newcommand{\matr}[1]{\underline{\underline{#1}}}

\newcommand{\am}{\matr{\av}}
\newcommand{\bmm}{\matr{\bv}}
\newcommand{\alpham}{\matr{\alphav}}
\newcommand{\betam}{\matr{\betav}}

\newcommand{\alphaMax}{\alpha_{\rm max}}
\newcommand{\Ca}{C_{\alpha}}

\newcommand{\SumiNn}{\sum_{i=0}^{\Nc_n-1}}
\newcommand{\SumkNn}{\sum_{k=0}^{\Nc_n-1}}

\usepackage{mathtools}

\usepackage{listings}
\usepackage{listings}
\newtheorem{theorem}{Theorem}

\newtheorem{remark}{Remark}[section]

\numberwithin{equation}{section}

\begin{document}

\begin{frontmatter}
 \title{High-order accurate schemes for Maxwell's equations with nonlinear active media and material interfaces}


\author[rpi]{Qing~Xia\fnref{DARPA}}
\ead{xiaq2@rpi.edu}

\author[rpi]{Jeffrey~W.~Banks\fnref{DARPA}}
\ead{banksj3@rpi.edu}

\author[rpi]{William~D.~Henshaw\corref{cor}\fnref{DARPA}}
\ead{henshw@rpi.edu}

\author[purdue]{Alexander~V.~Kildishev\fnref{DARPA}}
\ead{kildishev@purdue.edu}

\author[rpi]{Gregor~Kova\v ci\v c\fnref{DARPA}}
\ead{kovacg@rpi.edu}

\author[purdue]{Ludmila~J.~Prokopeva\fnref{DARPA}}
\ead{lprokop@purdue.edu}

\author[rpi]{Donald~W.~Schwendeman\fnref{DARPA}}
\ead{schwed@rpi.edu}

\address[rpi]{Department of Mathematical Sciences, Rensselaer Polytechnic Institute, Troy, NY 12180, USA}
\address[purdue]{School of Electrical and Computer Engineering, Purdue University, West Lafayette, IN 47907, USA}

\cortext[cor]{Corresponding author}

\fntext[DARPA]{This work was partially funded by the DARPA Defense Sciences Office, Award HR00111720032.}




\begin{abstract}
  We describe a fourth-order accurate finite-difference time-domain scheme for solving dispersive Maxwell's equations 
  with nonlinear multi-level carrier kinetics models. The scheme is based on an efficient single-step three time-level modified equation approach
  for Maxwell's equations in second-order form for the electric field coupled to ODEs for the polarization vectors and
  population densities of the atomic levels. The resulting scheme has a large CFL-one time-step.
 Curved interfaces between different materials are accurately treated with 
  curvilinear grids and compatibility conditions. A novel hierarchical modified equation approach leads to an explicit scheme
  that does not require any nonlinear iterations. The hierarchical approach at interfaces leads to local updates at the interface
  with no coupling in the tangential directions. Complex geometry is treated with overset grids. Numerical stability is maintained
  using high-order upwind dissipation designed for Maxwell's equations in second-order form.
  The scheme is carefully verified for a number of two and three-dimensional problems.
  The resulting numerical model with generalized dispersion and arbitrary nonlinear multi-level system can be used for many plasmonic applications such as
  for ab initio time domain modeling of nonlinear engineered materials for nanolasing applications, 
  where nano-patterned plasmonic dispersive arrays are used to enhance otherwise weak nonlinearity in the active media.
\end{abstract}


\end{frontmatter}

\tableofcontents

\clearpage 


\section{Introduction}\label{sec:intro}
The overall objective of this work is the development of efficient and high-order accurate numerical schemes for modeling light-matter interaction with nonlinear materials. 
We develop novel high-order finite-difference-time-domain (FDTD) numerical schemes for nonlinear active materials with carrier kinetics modeled by real-valued rate equations and the auxiliary differential equation (ADE) approach. The target nonlinear multilevel dispersive models and the geometry under consideration are universal in the sense that the number of atomic levels and the number of transitions (polarization vectors) are arbitrary, and the geometry can be of complex shapes with material interfaces in 2D or 3D, which essentially overcome the limitations of existing methods for full wave simulations in nonlinear active materials in the literature.
 There are a number of novel and attractive features of our schemes. 
 (1) Modified equation time-stepping leads to an extremely efficient three-level single-step scheme that is fourth-order accurate in space and time
 and has a large CFL-one time-step. 
 (2) High-order accuracy at curved boundaries and interfaces is achieved
 using conforming grids and compatibility conditions. (3) Interfaces are treated in an efficient way using a hierarchical modified equation (HIME) approach
 that provides local updates to the interface ghost points and requires no nonlinear solves.
 (4) Numerical stability on overset grids is achieved using a high-order upwind dissipation for Maxwell's equation in second-order form. 
 Upwind dissipation for wave equations in second-order form was first proposed in~\cite{sosup2012} and extended
 the ideas of Godunov's upwind scheme for first-order systems of equations. 
 The original scheme in~\cite{sosup2012}  was extended to Maxwell's equations in~\cite{mxsosup2018}.
 An optimized version of this latter approach is used in the numerical scheme given in this article.

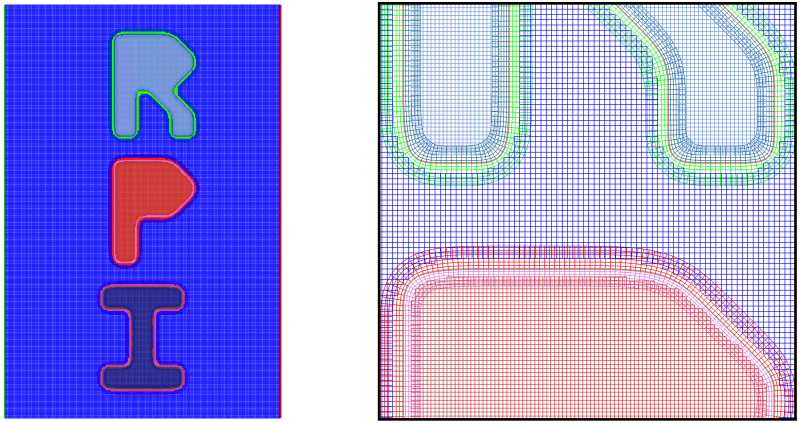
\begin{figure}[htb]
\begin{center}
\begin{tikzpicture}
   \useasboundingbox (0,.65) rectangle (10,5.5);  
  \figByHeight{   0}{0}{fig/rpiGridG4}{5.5cm}[0][0][0.][0]
  \figByHeightb{5}{0}{fig/rpiGridG4Zoom}{5.5cm}[0][0][0.][0]
\end{tikzpicture}
\end{center}
\caption{Left: overset grids for three meta-atoms. Right: zoomed-in view of the overset grids showing the curvilinear interface-fitted grids.}
\label{fig:rpiGrid}
\end{figure}

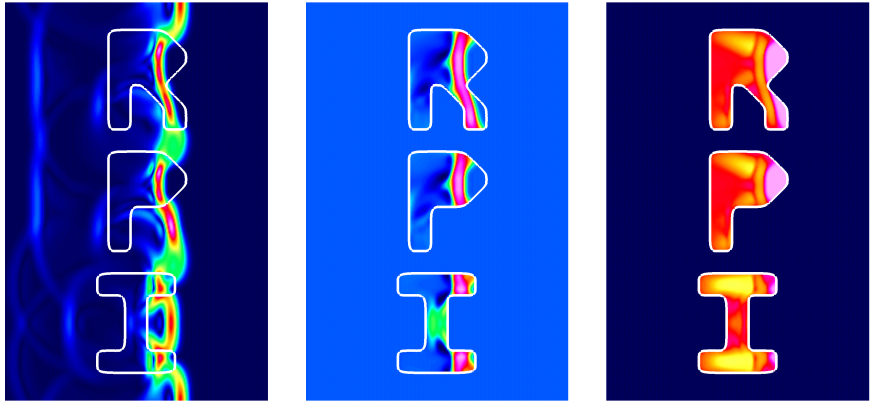
\begin{figure}[htb]
\begin{center}
\begin{tikzpicture}
   \useasboundingbox (0,.5) rectangle (12,5.5);  
  \figByWidth{   0}{0}{fig/rpiMLAEfieldNormt3}{3.5cm}[0][0][0.][0]
  \figByWidth{4}{0}{fig/rpiMLAPyt3}{3.5cm}[0][0][0.][0]
  \figByWidth{8}{0}{fig/rpiMLAN0t3}{3.5cm}[0][0][0.][0]
\end{tikzpicture}
\end{center}
\caption{Gaussian plane wave hitting three four-level active-material meta-atoms; electric field norm $\| \Ev\|$ (left),
$y$-component of the total polarization $P_y$ (middle) and ground-state population density $N_0$ (right)}
\label{fig:rpiFields}
\end{figure}

The current work substantially expands the previous work \cite{angel2019high,banks2020high} on general linear dispersive materials by employing multi-level carrier kinetics to account for diverse nonlinear effects. Similar to \cite{angel2019high,banks2020high}, the proposed schemes are implemented in the Overture framework, and then tested using overlapping grids over complex two- and three-dimensional geometries. 
Figures~\ref{fig:rpiGrid} and~\ref{fig:rpiFields} 
depict an example of a two-dimensional geometry and a zoom-in view of the associated overset grids.

Following the early theoretical works (see, e.g. \cite{siegman1986lasers}) and their original numerical FDTD approximations \cite{ziolkowski1995ultrafast, nagra1998FDTD, chang2004finite, azzam_2020lsa}, propagation and scattering of light in all our schemes are modeled classically by the Maxwell's equations, while the nonlinear medium is described by real-valued multi-level rate equations~\cite{ trieschmann2011experimental, azzam_2018lpr, azzam2021lpr, azzam2018sa, azzam2018rsa, azzam2020chapter6}. The interaction between light and matter is then coupled by real-valued atomic dipole moment (polarization) equations using the ADE technique, which is also adopted in the models of general linear dispersive materials in \cite{angel2019high,banks2020high}. Such nonlinear multilevel models are generalization of the two-level systems \cite{pantell1969fundamentals,ziolkowski1995ultrafast}. 

Perturbation theory has been and still remains one of the most popular approaches to approximate nonlinear material responses in optics. With this approach, the optical material response in the time and frequency domains is modeled employing the power series expansions of a weak nonlinear part of susceptibility. The method has become imperative for the initial numerical analysis of Stockman’s spasers \cite{bergman2003surface, stockman2008nat-photon}. Classical electrodynamics with the quantum-mechanical effects of the gain medium introduced through the perturbation nonlinear susceptibility terms have been initially adopted to describe spasing (see, e.g. \cite{Li2010design, baranov2013exactly, arnold2015spasers, kristanz2018power}). Thus, Li and Yu \cite{Li2010design} derived and computed the gain threshold requirements for core-shell single-particle spasers, accounting for the interband transitions of the plasmonic metal core. Kristanz et al. \cite{kristanz2018power} analyzed the power balance and heating to guide the spaser design in terms of the allowed pumping intensities, duration, and expected output radiation and thermal load. These studies have been of ultimate importance for analyzing the parameters affecting the threshold, including the resonant wavelength, the refractive index of the background host material, and the dimensions of the core and shell of regular-shape (mainly spherical or spheroidal) spasers.  While such models are capable of adequately predicting the conditions for loss compensation and the transition to the spasing regime for simplified geometries and operation regimes, as the designs of spaser systems are becoming more involved, full-wave numerical analysis that can unlock the temporal and spatial details of a given spaser are required. In general, 
the perturbation theory has many restrictions. For example, the modeling techniques employing this classical approach are capable of neither capturing complete transient and irreversible effects nor accounting for many critical quantum phenomena. They may also fail to converge in some crucial real-life cases \cite{Bravo-Abad2007modeling} and are inadequate for modeling several distinct classes of nonlinearities, e.g., epsilon-near-zero materials \cite{reshef2017nonlinear}.
In contrast, the time-domain multiphysics techniques are considered amongst the most accurate numerical frameworks that can account for the quantum-mechanical nature of the gain and plasmonic materials, naturally combining nonlinear and thermal effects in a single computational domain with complex structural and material composition \cite{nagra1998FDTD,chang2004finite,chua2014modeling,pusch2012coherent,trivedi2017model, azzam_2018lpr}.

Early carrier kinetics approaches to multiphysics modeling of nonlinear light-matter interaction were introduced for describing gain media in response to external pulsed excitation \cite{siegman1986lasers,ziolkowski1995ultrafast}, aiming at the numerical analysis of the dynamics of pumping, population inversion, and saturation. The multi-level rate equation technique has been widely used for simulating various atomic systems, for example, in modeling 1-electron system with 4 levels and 6 levels~\cite{trieschmann2011experimental}, 2-electron system with 4 levels~\cite{chang2004finite}, saturable absorption \cite{azzam2018sa}, reverse saturable absorption~\cite{azzam2018rsa}, and 2-photon absorption \cite{azzam2020chapter6}. Modeling multi-level active medium using the rate equations together with the ADE-type polarizations is equivalent to the first-order optical Bloch equations formulated using density matrix \cite{boyd2020nonlinear,allen1987optical} for two-level systems, or multi-level systems consisting of pairwise atomic level transitions \cite{chang2004finite, taflove2006erratum, taflove2005book}. For general multi-level atomic systems, such as those with $V$, $\Lambda$ or cascade configurations \cite{maimistov1990present}, the equivalence does not hold. However, one could fit the multilevel models using, for example, experimental data by leaving out non-essential transitions \cite{trieschmann2011experimental}.


There have been many numerical methods that were developed for the complex-valued optical Bloch equations based on the density matrix, see the review paper \cite{jirauschek2019optoelectronic} for example, among which the finite difference methods have been prevailing in the time-domain multiphysics techniques. In the FDTD regime, Yee's scheme \cite{yee1966numerical} was widely used. For instance, in \cite{bidegaray2003time}, a weakly decoupled and Strang splitting time discretizations of Maxwell-Bloch system that preserves the carrier populations was discussed, with a feature of using different time marching for diagonal and off-diagonal entries of the density matrix. Yee's scheme was extended to the light-matter interactions with ultrashort pulse in anisotropic media for the unidimensional case in \cite{saut2004computational}, and bidimensional case in \cite{bourgeade2006numerical}, where a pesudo-spectral time-domain method was also discussed, along with similar splitting schemes for Bloch equations. In \cite{2020arXiv200505412R}, a Maxwell-Bloch solver for two-level atomic systems was developed. A scrutiny of FDTD numerical methods based on Yee's scheme for nonlinear active materials revealed that they have limitations and restrictions to low-order accuracy, low dimensions, or only two-level atomic systems. Nor are they capable of handling geometry with complex shaped boundaries and interfaces.

The rest of the paper is outlined as follows. In Section~\ref{sec:eqns}, we present the mathematical models of nonlinear dispersive materials that use multi-level rate equations, employing the ADE technique. Here, we prescribe the initial/boundary conditions, and interface jump conditions for the second-order Maxwell's equations. In Section~\ref{sec:scheme}, the second-order finite difference time-stepping schemes for nonlinear models in homogeneous materials are discussed first, whereas the fourth-order schemes that employ the second-order results and the modified equation approach are presented in sequel. The second- and fourth-order numerical interface treatments are discussed in Section~\ref{sec:interface}. 
Lastly, several numerical examples in both 2D and 3D are given in Section~\ref{sec:numerics} and concluding remarks are given in Section~\ref{sec:conclusion}.

\section{Governing equations}\label{sec:eqns}

We consider the solution to the initial-boundary-value (IBVP) problem for Maxwell's equations in a domain $\Omega \subset \Real^{n_d}$ in $n_d$ space dimensions.
The domain consists of $N_k$ different material regions $\Omega_k$ with $\Omega=\cup_{k=1}^{N_k} \Omega_k$. Let $\Gamma_{k,k'}$ denote the
interface between material $k$ and $k'$. 
A given region $\Omega_k$ may be governed by the isotropic Maxwell equations, the linear dispersive Maxwell's equations as discussed
in~\cite{banks2020high} or the following nonlinear equations, 
\begin{subequations}\label{eqn:mbe}
\bat
  &  \partial_t^2\Ev = c^2\Delta \Ev-\epsilon^{-1}_0\sum_{m=1}^{\Nc_p}\partial_t^2\Pv_{m}   ,    \label{eqn:nonlin_e}\\
  &  \partial_t^2\Pv_m + b_{1,m}\partial_t\Pv_m + b_{0,m} \Pv_m = \sum_{\ell=1}^{\Nc_n} a_{m,\ell} N_\ell \Ev ,  
              \qquad && m=1,2,\ldots,\Nc_p,  \label{eqn:nonlin_p}\\
  &  \partial_t N_\ell = \sum_{\hat{\ell}=0}^{\Nc_n-1} \alpha_{\ell,\hat{\ell}}N_{\hat{\ell}} 
       + \sum_{m=1}^{\Nc_p} \beta_{\ell,m}\Ev\cdot\partial_t \Pv_m ,    \quad&& \ell=0,1,2,\ldots,\Nc_\ell-1 .\label{eqn:nonlin_n}
\eat
\end{subequations}
Equations~\eqref{eqn:mbe}, called the {\em Maxwell-MLA} system, define Maxwell's equations 
in second-order form for the electric field coupled to a multi-level carrier kinetic model.
The kinetic model consists of $\Nc_n$ atomic levels for population densities $N_{\ell}$ and $\Nc_p$ polarization vectors $\Pv_m$.
The parameter $c=1/\sqrt{\eps_0 \mu_0}$ is the speed of light in a vacuum, with $\eps_0$ and $\mu_0$ the vacuum permittivity and permeability,
respectively. The parameters 
$b_{1,m}$, $b_{0,m}$, $a_{m,\ell}$, $\alpha_{\ell,k}$ , and $\beta_{\ell,m}$ are all real and chosen to model the transitions in
a particular active material, either based on theory or experimental data. 
In a typical case the sum of the population densities $N_\ell$ will be constant, often normalized to be one.
Note that in subsequent discussions the bounds on the sums in equations in~\eqref{eqn:mbe} will often be suppressed for notational brevity.
Also note that, following our previous work~\cite{angel2019high,banks2020high,max2006b},
we solve for $\Ev$ using Maxwell's equations in second-order form.
The advantages of using the second-order form Maxwell's equations are described, for example, in~\cite{max2006b}.

\begin{figure}[H]
\centering
\begin{tikzpicture}[scale=1]
\draw[dashed,thick] (0,0) -- (4.5,0);
\draw[dashed,thick] (2,1) -- (4.5,1);
\draw[dashed,thick] (1,2) -- (4.5,2);
\draw[dashed,thick] (0,3) -- (4.5,3);
\draw[->,thick] (0.25,0) -- (0.25,3);
\draw[->,dotted,thick] (0.5,3) -- (0.75,0);
\draw[->,dotted,thick] (1.5,3) -- (1.75,2);
\draw[->,dotted,thick] (2.5,1) -- (2.75,0);
\draw[->,dotted,thick] (3.75,2) -- (4,1);
\draw[->,thick] (3.6,1) -- (3.6,2);
\draw (0.25+0.1,2.5) node[left] {$P_{30}$};
\draw (0.5,2.5) node[right] {$\tau_{30}$};
\draw (1.5+0.1,2.5) node[right] {$\tau_{32}$};
\draw (2.9+0.,1.5) node[right] {$P_{21}$};
\draw (3.75+0.1,1.5) node[right] {$\tau_{21}$};
\draw (2.5+0.1,0.5) node[right] {$\tau_{10}$};
\draw (4.5+0.1,0) node[right] {$N_0$};
\draw (4.5+0.1,1) node[right] {$N_1$};
\draw (4.5+0.1,2) node[right] {$N_2$};
\draw (4.5+0.1,3) node[right] {$N_3$};
\end{tikzpicture}
\caption{Jablonski diagram for a 4-level atomic system}
\label{fig:4-level}
\end{figure}
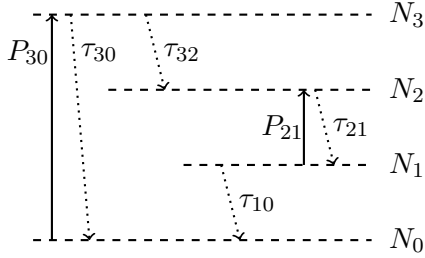

The kinetic model in~\eqref{eqn:nonlin_n} is quite general and can represent transitions in a variety of 
multi-level atomic systems. Consider, for example, the four-level system depicted in Figure~\ref{fig:4-level}.
This system consists of four energy levels with population densities $N_\ell$, $\ell=0,1,2,3$. 
Here we change notation slightly to be consistent with the literature.
This is a model for lasing in a gain medium.
Energy is pumped into the ground level, where the electrons are excited to the highest energy level 3, then relaxed to lower energy levels. With proper constraints on the relaxation time $\tau$'s, a population inversion (more populations at a higher energy level) between levels 1 and 2 will take place, which leads to lasing. 
The corresponding rate equations for the population densities $N_\ell$ are given by
\newcommand{\FF}{\vphantom{\begin{bmatrix}
    N_0\\
    N_1
\end{bmatrix}}}
\bas
\begin{bmatrix}
    \partial_t N_0\\
    \partial_t N_1\\
    \partial_t N_2\\
    \partial_t N_3
\end{bmatrix} =&
\begin{bmatrix}
    0 & \tau^{-1}_{10} & 0 & \tau^{-1}_{30}\\
    0 & -\tau^{-1}_{10} & \tau^{-1}_{21} & 0\\
    0 & 0 & -\tau^{-1}_{21} & \tau^{-1}_{32}\\
    0 & 0 & 0 & -\tau^{-1}_{30}-\tau^{-1}_{32}\\
\end{bmatrix}
\begin{bmatrix}
    N_0\\
    N_1\\
    N_2\\
    N_3
\end{bmatrix}   
   + 
\begin{bmatrix}
   -(\hbar\omega_{30})^{-1}& 0\\
   0&-(\hbar\omega_{21})^{-1}\\
   0&(\hbar\omega_{21})^{-1}\\
   (\hbar\omega_{30})^{-1}&0\\
\end{bmatrix}
\begin{bmatrix}
    \Ev\cdot\partial_t\Pv_{30} \FF  \\
    \Ev\cdot\partial_t\Pv_{21} \FF
\end{bmatrix} ,
\eas
where $\hbar$ is the reduced Planck constant, $\omega_{21},\omega_{30}$ are the transition frequencies between the paired levels, and $\tau_{10},\tau_{30},\tau_{21},\tau_{32}$ are the relaxation times from the higher energy level to corresponding low level respectively.
The associated polarizations can be expressed as
\ba
    \partial_t^2\Pv_{ji} + \gamma_{ji}\partial_t\Pv_{ji} + \omega_{ji}^2 \Pv_{ji} 
      = \kappa_{ji} (N_i-N_j) \Ev,\qquad ji = \{30,21\}, 
\ea
and thus the tensor $\am$ with entries $a_{m,n}$  in the polarization equation~\eqref{eqn:nonlin_p} is given by
\ba
  \am = 
  \begin{bmatrix}
   \kappa_{30} & 0 & 0 & -\kappa_{30} \\
   0 & \kappa_{21} & -\kappa_{21} & 0\\
  \end{bmatrix} ,
\ea
while $\Pv_1=\Pv_{30}$, and $\Pv_{2} = \Pv_{21}$.

To define a well-posed IBVP, the Maxwell-MLA equations~\eqref{eqn:mbe} are augmented with appropriate initial conditions, 
boundary conditions and interface conditions. Initial conditions are required for $\Ev$, $\p_t\Ev$, $\Pv_m$, $\p_t\Pv_m$ and
$N_\ell$. 
For the purposes of this article, the nonlinear materials will be bounded by linear materials and thus will not require
boundary conditions. Boundary conditions at physical or far-field boundaries for linear materials will be specified in 
the usual way as discussed in~\cite{banks2020high}. Note that the second-order form of the equations for $\Ev$ uses
the additional boundary condition $\grad\cdot\Ev=0$.
At an interface $\Gamma_{k,k'}$
between sub-domains $\Omega_k$ and $\Omega_{k'}$ the following primary interface conditions hold
\bse
\label{eqn:primary-ic}
\ba
   \big[  \nv\times\Ev                        \big]_{\Gamma_{k,k'}}  =0 , \label{eqn:ic-1}\\
   \big[  \nv\cdot (\epsilon_0\Ev+\Pv)        \big]_{\Gamma_{k,k'}}  =0 , \label{eqn:ic-2}\\
   \big[  \mu_0^{-1}\nv\times\nabla\times\Ev  \big]_{\Gamma_{k,k'}}  =0 , \label{eqn:ic-3}\\
   \big[  \nabla\cdot\Ev                      \big]_{\Gamma_{k,k'}}  =0 , \label{eqn:ic-4}
\ea
\ese
where $\nv$ is defined to be the normal that points from domain $\Omega_k$ into $\Omega_{k'}$.

  The well-possedness of the IBVP for the Maxwell-MLA system~\eqref{eqn:mbe} is discussed in~\ref{sec:stability}.
The problem is well-posed and the solutions to these nonlinear equations will exist for at least short times.
Long-time existence can be shown for a restricted class of commonly used systems, such as the four-level system
described above. For such systems an $L_2$-energy can be found that shows the solutions have at most
bounded exponential growth in time.

\section{Numerical Scheme}\label{sec:scheme}

\bogus{
{\color{blue}First Questions/Comments
\begin{enumerate}
  \item In the model, are there exterior BCs considered, or only interfaces?
  \item Should the interface discretization section go as a subsection in Section 3?
  \item Can we change the \bibtex names of MR3881594,MR4081505?
  \item Are we planning to pursue any stability analysis? If we do the energy estimate for the continuous setting, perhaps stability for the discrete system follows since the ML equations don't add any difficulties at the boundaries?
\end{enumerate}}
}

The basic approach to discretization of the MLA equations~\eqref{eqn:mbe} uses finite-difference approximations and modified-equation time-stepping. This approach 
follows the path previously advocated for the 
nondispersive isotropic Maxwell's equations in~\cite{max2006b}, later extended to linear dispersive materials in~\cite{angel2019high}, and subsequently to linear dispersive materials with interfaces in~\cite{banks2020high}. 
The primary developments described in the present article are the formulation and application of methods for the equations of {\em nonlinear} electromagnetics for active media, and a novel approach to the treatment of interfaces that eliminates the need for the solution of coupled nonlinear systems of equations along material interfaces.
 The treatment of complex geometry will again make use of overlapping grids, which is discussed briefly in Section~\ref{sec:overlapping}. 
 Second-order and fourth-order accurate discretizations are then discussed in Sections~\ref{sec:2nd} and~\ref{sec:4th} respectively. 
 The discretization of interface equations is then presented in Section~\ref{sec:interface}. 
 See Algorithm~\ref{alg:overview} for the overview of the developed algorithms and the arrangements of this section.


\begin{algorithm}
\algFontSize 
\caption{Overview of the developed algorithms}
\begin{algorithmic}[1]
    \State Generate overset grids for the geometry; \Comment Sect.~\ref{sec:overlapping}
    \State Initialization;
    \While{$t<T_{final}$} \Comment Begin time-stepping loop
        \For{$\bm{j}$ in each grid $G$} \Comment Sect.~\ref{sec:2nd},~\ref{sec:4th}
          \For{$m=1,\dots, \mathcal{N}_p$}
              \State Update $\Pv^{n+1}_{m,\bm{j}}$;
          \EndFor
          \State Update $\Ev^{n+1}_{\bm{j}}$;
          
          \For{$\ell=0,\dots, \mathcal{N}_n-1$}
              \State Update $\Nv^{n+1}_{\ell,\bm{j}}$;
          \EndFor
        \EndFor 
        \State Apply boundary and interface conditions; \Comment Sect.~\ref{sec:ic2},~\ref{sec:ic4}

        \State $t^{n+1} = t^{n} + \dt$, $n = n + 1$;
    \EndWhile    \Comment End time-stepping loop
\end{algorithmic} 
\label{alg:overview}
\end{algorithm}


\subsection{Overlapping grids}\label{sec:overlapping}

As indicated in the introduction, geometric complexities in the simulation domain will be addressed using {\em overlapping} (sometimes referred to as {\em overset}, or {\em chimera}) grids. An example is depicted in Figs.~\ref{fig:rpiGrid} and~\ref{fig:rpiFields}, where a light pulse is propagated through active material in the shape of the letters ``R'', ``P'', and ``I''. Each computational subdomain is discretized using a composite overlapping grid consisting of a set of thin boundary fitted grids overlaying a Cartesian background grid, see e.g. Fig.~\ref{fig:rpiGrid}. The primary motivation for our use of composite overlapping grids is to enable the use of efficient finite difference schemes on structured grids, while simultaneously treating complex geometry with high-order accuracy up to and including boundaries and material interfaces.
\begin{figure}[hbt]
\begin{center}
\begin{tikzpicture}[scale=.7]
\useasboundingbox (.75,1) rectangle (15.5,6.25);  
%
\begin{scope}[xshift=1cm,yshift=1cm]
\fill[black!10!white,xshift=.5cm,yshift=.5cm] (0,0) -- (2.583333,0) arc (180:90:1.416667) -- (4.,4.) -- (0,4.) -- (0,0);
\draw[-,thick,blue,yshift=.0 cm] 
   \foreach \x/\y in {1.5/0,1.5/.5,2/1,2/1.5,2.5/2,3/2.5,4/3,5/3.5,5/4,5/4.5,5/5}{ (0,\y) -- (\x,\y) }
   \foreach \x/\y in {0/0,.5/0,1/0,1.5/0,2/1,2.5/2,3/2.5,3.5/3,4/3,4.5/3.5,5/3.5}{ (\x,\y) -- (\x,5) };
  \begin{scope}[xshift=4.5cm,yshift=0.5cm]
    \draw[thick,green] \foreach \r in {1.000000,1.416667,1.833333,2.250000,2.666667,3.083333,3.500000}{ (0,\r) arc (90:190:\r)  (0,\r) arc (90:80:\r) };
    \draw[thick,green]
     (0.173648,0.984808)  -- (0.607769,3.446827)
     (0.000000,1.000000)  -- (0.000000,3.500000)
     (-0.173648,0.984808) -- (-0.607769,3.446827)
     (-0.342020,0.939693) -- (-1.197071,3.288924)
     (-0.500000,0.866025) -- (-1.750000,3.031089)
     (-0.642788,0.766044) -- (-2.249757,2.681156)
     (-0.766044,0.642788) -- (-2.681156,2.249757)
     (-0.866025,0.500000) -- (-3.031089,1.750000)
     (-0.939693,0.342020) -- (-3.288924,1.197071)
     (-0.984808,0.173648) -- (-3.446827,0.607769)
     (-1.000000,0.000000) -- (-3.500000,0.000000)
     (-0.984808,-0.173648) -- (-3.446827,-0.607769);
  \end{scope}
  \draw[very thick,red,xshift=.5cm,yshift=.5cm] (0,0) -- (2.583333,0) arc (180:90:1.416667) -- (4.,4.) -- (0,4.) -- (0,0);
%
   \filldraw[green] (1.5,.5)  circle (3pt)
                 (1.5,1 )  circle (3pt)
                 (2  ,1 )  circle (3pt)
                 (2 ,1.5)  circle (3pt)
                 (2 , 2 )  circle (3pt)
                 (2.5,2 )  circle (3pt)
                 (2.5,2.5) circle (3pt)
                 (3 , 2.5) circle (3pt)
                 (3  ,3 )  circle (3pt)
                 (3.5,3 )  circle (3pt)
                 (4  ,3. ) circle (3pt)
                 (4  ,3.5) circle (3pt)
                 (4.5,3.5) circle (3pt);
%
  \begin{scope}[xshift=4.5cm,yshift=0.5cm]
      \filldraw[blue]
       (0.000000,3.500000)    circle (3pt)
       (-0.607769,3.446827)   circle (3pt)
       (-1.197071,3.288924)  circle (3pt) 
       (-1.750000,3.031089)  circle (3pt) 
       (-2.249757,2.681156)  circle (3pt) 
       (-2.681156,2.249757)  circle (3pt) 
       (-3.031089,1.750000)  circle (3pt) 
       (-3.288924,1.197071)  circle (3pt) 
       (-3.446827,0.607769)  circle (3pt) 
       (-3.500000,0.000000)  circle (3pt);
  \end{scope}
   \draw (1.25,3.5) node[thick,draw=blue,fill=white] {\large$\Gv_1$};
   \draw (3.25,1.95) node[thick,draw=green,fill=white] {\large$\Gv_2$};
\end{scope}
%
\definecolor{ghostColour}{named}{DodgerBlue}
\newcommand{\mytrix}{(\x,-.15) -- ++(.3,0) -- ++(-.15,.26) -- (\x,-.15)}
\newcommand{\mytriy}{(-.15,\y) -- ++(.3,0) -- ++(-.15,.26) -- (-.15,\y)}
\begin{scope}[xshift=7cm,yshift=2.25cm,scale=.75]
\draw[-,thick,blue,yshift=.0 cm] 
   \foreach \x in {0,.5,...,5}{ (\x,0) -- (\x,5) }
   \foreach \y in {0,.5,...,5}{ (0,\y) -- (5,\y) };
  \draw[very thick,red,xshift=.5cm,yshift=.5cm] (1.,0) -- (.0,0) -- (.0,4.) -- (4.,4.) -- (4.,3.);
   \filldraw[green] (1.5,.5)  circle (3pt)
                 (1.5,1 )  circle (3pt)
                 (2  ,1 )  circle (3pt)
                 (2 ,1.5)  circle (3pt)
                 (2 , 2 )  circle (3pt)
                 (2.5,2 )  circle (3pt)
                 (2.5,2.5) circle (3pt)
                 (3 , 2.5) circle (3pt)
                 (3  ,3 )  circle (3pt)
                 (3.5,3 )  circle (3pt)
                 (4  ,3. ) circle (3pt)
                 (4  ,3.5) circle (3pt)
                 (4.5,3.5) circle (3pt);
  \filldraw[fill=white,draw=black]  \foreach \x in {2,2.5,...,5}{ (\x,.0) circle (3.5pt) };
  \filldraw[fill=white,draw=black]  \foreach \x in {2,2.5,...,5}{ (\x,.5) circle (3.5pt) };
  \filldraw[fill=white,draw=black]  \foreach \x in {2.5,3,...,5}{ (\x,1.) circle (3.5pt) };
  \filldraw[fill=white,draw=black]  \foreach \x in {2.5,3,...,5}{ (\x,1.5) circle (3.5pt) };
  \filldraw[fill=white,draw=black]  \foreach \x in {3,3.5,...,5}{ (\x,2.0) circle (3.5pt) };
  \filldraw[fill=white,draw=black]  \foreach \x in {3.5,4,...,5}{ (\x,2.5) circle (3.5pt) };
  \filldraw[fill=white,draw=black]  \foreach \x in {4.5,5}      { (\x,3.0) circle (3.5pt) };
  \draw[fill=ghostColour,xshift=-.15cm,yshift=0cm]  \foreach \x in {.5,1.,1.5}{ \mytrix };  
  \draw[fill=ghostColour,xshift=-.15cm,yshift=5cm]  \foreach \x in {.5,1.,...,5}{ \mytrix };  
  \draw[fill=ghostColour,xshift=0cm,yshift=-.15cm]  \foreach \y in {0,.5,...,5}{ \mytriy };
  \draw[fill=ghostColour,xshift=5cm,yshift=-.15cm]  \foreach \y in {3.5,4,4.5}{ \mytriy };
   \draw (1.25,3.5) node[thick,draw=blue,fill=white] {\large$\Gv_1$};
\end{scope}
\begin{scope}[xshift=11.5cm,yshift=2.25cm,scale=.75]
\draw[-,thick,green,yshift=.0 cm] 
   \foreach \x in {0,.454545,...,5}{ (\x,0) -- (\x,5) }
   \foreach \y in {0,.833333,...,5}{ (0,\y) -- (5,\y) };
 \draw[very thick,red,xshift=.454545cm,yshift=.833333cm] (0.,4) -- (.0,0) -- (4.0909,0.) -- (4.0909,4);
 \filldraw[blue]  \foreach \x in {.454545,.909090,...,4.545454}{ (\x,5) circle (3.5pt) };
 \draw[fill=ghostColour,xshift=-.15cm]  \foreach \x in {.454545,.909090,...,4.545454}{ \mytrix };
 \draw[fill=ghostColour,yshift=-.15cm]  \foreach \y in {0,.833333,...,5}{ \mytriy };
 \draw[fill=ghostColour,xshift=5cm,yshift=-.15cm]  \foreach \y in {0,.833333,...,5}{ \mytriy };
\end{scope}
\begin{scope}[xshift=7cm,yshift=.6cm]
  \fill[Goldenrod!20!white,draw=black,xshift=-.1cm,yshift=-.25cm] (-.1,-.1) -- (4.1,-.1) -- (4.1,1.4) -- (-.1,1.4) -- (-.1,-.1);
  \filldraw[green,xshift=.0cm,yshift=.8cm] (.25,.0)  circle (3pt);
  \filldraw[blue,xshift=.3cm,yshift=.8cm] (.25,.0)  circle (3pt);
  \draw[xshift=.0cm,yshift=.8cm] (.5,0) node[anchor=west,xshift=6] {\small interpolation};
  \draw[fill=ghostColour,xshift=.0cm,yshift=.4cm] (.35,0) \foreach \x in {.1}{ \mytrix } node[anchor=west,xshift=12,yshift=3] {\small ghost};
  \draw[fill=white,draw=black,xshift=.0cm,yshift=.0cm] (.25,0) circle (3.5pt) node[anchor=west,xshift=6] {\small unused};
\end{scope}
\begin{scope}[xshift=11.5cm,yshift=2.25cm,scale=.75]
   \draw (1.6,3.27) node[thick,draw=green,fill=white] {\large$\Gv_2$};
\end{scope}
\end{tikzpicture}
\end{center}
\caption{Left: an overlapping grid consisting of two
structured curvilinear component grids, $\xv=G_1(\rv)$ and $\xv=G_2(\rv)$. Middle and right: 
component grids for the square and annular grids in the unit square parameter space $\rv$. Grid
 points are classified as discretization points, interpolation points or unused points. Ghost points
 are used to apply boundary conditions.}    \label{fig:overlappingGridCartoon}
\end{figure}
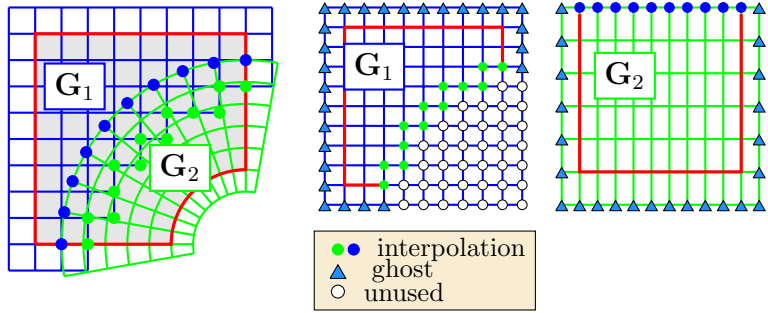
In the composite overlapping grid approach, the simulation domain $\Omega$, is divided into its geometric components $\Omega_k$ (e.g. the ``R'' domain in Fig.~\ref{fig:rpiGrid}). Each sub-domain $\Omega_k$ is then covered by a composite grid $\Gc_k$, consisting of a set of component grids $\Gc_{k,g}$, $g=1,\ldots,N_k$. A simple example composite grid in two space dimensions is illustrated in Fig.~\ref{fig:overlappingGridCartoon}. Each component grid $\Gc_{k,g}$ is a logically rectangular, curvilinear grid defined by a smooth mapping from a reference domain $\rv\in[0,1]^{\nd}$ (i.e. the unit square in 2D or unit cube in 3D) to physical space~$\xv$,
\begin{equation}
  \xv = \Gv_{k,g}(\rv),\qquad \rv\in[0,1]^{\nd},\qquad \xv\in\Real^{\nd}.
\label{eq:gridMapping}
\end{equation}
Grid points are classified as discretization points (where the PDE or boundary/interface conditions are applied), interpolation points (where solutions are interpolated from other component grids) or unused points. Throughout the present work, ghost points are use to implement boundary and interface conditions. The overlapping grid generator {\bf Ogen}~\cite{ogen} from the {\it Overture} framework is used to construct the overlapping grids. Overlapping grid interpolation is performed using a tensor-product Lagrange basis with quadratic polynomials for the  second-order accurate scheme, and quartic polynomials fourth-order scheme, as required to maintain accuracy~\cite{CGNS}.

\subsection{Discretizations for Nonlinear and Active Media}

Discretization of the governing PDE system~\eqref{eqn:mbe}, on each component grid $\Gc_{k,g}$, is performed in the reference coordinate system $\rv$. The overall approach taken here follows the general principles described in~\cite{max2006b,angel2019high}, and uses a single-step modified equation (ME) scheme (sometimes referred to as a {\em space-time} or {\em Lax-Wendroff} time stepper). To describe the schemes, denote $\xv_{\jv}\in\mathbb{R}^{\nd}$ as a point on a component grid, where $\jv=(j_1,\dots,j_d)\in\mathbb{Z}^{\nd}$ is a multi-index. 
Generically, ME time stepping schemes are based on a Taylor expansion of discrete approximations to temporal derivatives. For example, discretization of the leading second derivative terms in \eqref{eqn:nonlin_e} and \eqref{eqn:nonlin_p} can be based on the usual 3-level approximation of the second time derivative along with appropriate correction terms to obtain the required accuracy. Specifically, for schemes of order $p=2q$
\newcommand{\nn}{\nu}
\begin{align}
  \frac{\Wv_{\jv}(t+\dt)-2\Wv_{\jv}(t)+\Wv_{\jv}(t-\dt)}{\dt^2}   = \sum_{\nn=1}^{q} \frac{2\dt^{2(\nn-1)}}{(2\nn)!}\partial_t^{2\nn} \Wv_{\jv}(t) + \Oc(\dt^{p}), \label{eqn:expansion2} 
\end{align}
where $\dt$ is a time step size, and $\Wv_{\jv}(t)$ indicates a generic grid function and could be any of $\Ev_{\jv}(t),\Pv_{m,\jv}(t), N_{\ell,\jv}(t)$. 
On the other hand, \eqref{eqn:nonlin_n} is a first-order ODE 
and so the schemes are based on the forward difference approximation to the first derivative  
\begin{align}
  \frac{\Wv_{\jv}(t+\dt)-\Wv_{\jv}(t)}{\dt}   = \sum_{\nn=1}^{p} \frac{\dt^{\nn-1}}{\nn!}\partial_t^{\nn} \Wv_{\jv}(t) + \Oc(\dt^{p}).\label{eqn:expansion1} 
\end{align}
Repeated time differentiation of the PDE system~\eqref{eqn:mbe} is then used to define the various terms on the right-hand-side of \eqref{eqn:expansion2} and \eqref{eqn:expansion1}. To obtain a fully discrete scheme, spatial derivative operators are then replaced with difference approximations\footnote{Both conservative and non-conservative representation of the discrete Laplacian may be used, see~\cite{max2006b} for additional  details.}. See~\cite{max2006b,angel2019high} for additional details of this derivation for the non-dispersive and linearly dispersive Maxwell's equations respectively.


While the ME prescription above is correct and effective, straightforward implementation for the nonlinear dispersive equations \eqref{eqn:mbe} will necessitate the need to solve a globally coupled system of nonlinear equations at each time level. A similar effect was observed for linear dispersive materials in~\cite{angel2019high}, which led to the development of a redesigned ME scheme using a predictor-corrector methodology using a mixture (a hierarchy) of second- and fourth-order approximations. This methodology, subsequently referred to as a hierarchical modified equation (HIME) scheme, is adapted here for the nonlinear equations, and schemes of order 2 and 4 are described in Sections~\ref{sec:2nd} and~\ref{sec:4th} respectively.



\subsubsection{Second-order accurate scheme}
\label{sec:2nd}

  In this section we outline the second-order accurate Maxwell-MLA scheme.
A pseudo-code version of this algorithm is given in~\ref{sec:algorithmOrder2}.
In describing the discretization of~\eqref{eqn:mbe}, let $\Ev^n_{\jv},\Pv^n_{m,\jv}, N^n_{\ell,\jv}$ be approximations of $\Ev(\xv_{\jv},t^n),\Pv_m(\xv_{\jv},t^n)$, and $N_{\ell}(\xv_{\jv},t^n)$, respectively, at time $t^n=n\dt$. Further, let $\Delta_{ph}$ denote a $p$th-order accurate approximation to the Laplace operator $\Delta$, and
$\Dpt$, $\Dmt$, and $\Dzt$ denote the usual forward, backward, and central divided difference approximations
to the time derivative, as given by
\ba
   \Dpt \Wv_{\jv}^n \eqdef \f{\Wv_{\jv}^{n+1} - \Wv_{\jv}^n}{\dt}, \quad 
   \Dmt \Wv_{\jv}^n \eqdef \f{\Wv_{\jv}^{n} - \Wv_{\jv}^{n-1}}{\dt}, \quad 
   \Dzt \Wv_{\jv}^n \eqdef \f{\Wv_{\jv}^{n+1} - \Wv_{\jv}^{n-1}}{2\dt},
\ea
for a generic grid function $\Wv_{\jv}^{n}\approx\Wv(\xv_{\jv},t^n)$. 

Using this notation, and the expansions \eqref{eqn:expansion2} and \eqref{eqn:expansion1}, second-order accurate approximations to \eqref{eqn:nonlin_e}--\eqref{eqn:nonlin_n} are straight forward, with perhaps the simplest\footnote{Other discretizations involving alternate temporal weighting, e.g.  $b_{0,m}\left(\frac{1}{4}\Pv^{n+1}_{m,\jv} + \frac{1}{2}\Pv^n_{m,\jv} + \frac{1}{4}\Pv^{n-1}_{m,\jv}\right)$ in place of $b_{0,m}\Pv^n_{m,\jv}$, are also possible and could serve as the basis for higher-order schemes, but these are not pursued here.} being given by 
\begin{subequations}\label{eqn:mbe2}
  \begin{align}
    &D_{+t}D_{-t}\Ev^n_{\jv} = c^2\Delta_{2h}\Ev^n_{\jv}-\epsilon^{-1}_0D_{+t}D_{-t}\Pv^n_{\jv},\label{eqn:E-update-2}\\
    &D_{+t}D_{-t}\Pv^n_{m,\jv} +b_{1,m} D_{0t}\Pv^n_{m,\jv}+b_{0,m}\Pv^n_{m,\jv} = \sum_{\ell}a_{m,\ell}N_{\ell}^n\Ev^n_{\jv},\label{eqn:P-update-2}\\
    &D_{+t}N^{n}_{\ell,\jv} = \left.D_{2t}N_\ell\right|^{n}_{\jv} +\frac{\dt}{2} \left.D_{2tt}N_\ell\right|^{n}_{\jv},\label{eqn:n-update-2}
  \end{align}
\end{subequations}
for $m=1,\dots,\mathcal{N}_p$ and $\ell=0,\dots,\mathcal{N}_n-1$, with $\Pv_{\jv}^n \defeq \sum_{m}\Pv_{m,\jv}^n$. The notation $\left.D_{2t}N_\ell\right|^{n}_{\jv}$ and $\left.D_{2tt}N_\ell\right|^{n}_{\jv}$, used in \eqref{eqn:n-update-2}, is meant to indicate 2nd-order accurate approximation to the continuous time derivatives $\left.\partial_tN_\ell\right|^{n}_{\jv}$ and $\left.\partial_t^2N_\ell\right|^{n}_{\jv}$ respectively. In principle there are a number of choices for these approximations, e.g. backward differencing or implicit temporal averaging. However, in this work we use the Taylor series approach with
\begin{subequations}
\begin{align}
  \left.D_{2t}N_\ell\right|^n_{\jv} &\defeq \sum_{\hat{\ell}}\alpha_{\ell,\hat{\ell}}N^n_{\hat{\ell},\jv} + \sum_{m}\beta_{\ell,m}\Ev^n_{\jv}\cdot D_{0t} \Pv^n_{m,\jv},\label{eqn:2-accurate-Nt}\\
  \left.D_{2tt}N^{n}_\ell\right|_{\jv}^n &\defeq \sum_{\hat{\ell}}\alpha_{\ell,\hat{\ell}}\left.D_{2t}N_{\hat{\ell}}\right|^n_{\jv} + \sum_{m}\beta_{\ell,m}D_{0t}\Ev^n_{\jv}\cdot D_{0t} \Pv^n_{m,\jv}+ \sum_{m}\beta_{\ell,m}\Ev^n_{\jv}\cdot D_{+t}D_{-t} \Pv^n_{m,\jv},\label{eqn:2-accurate-Ntt}
\end{align}
\end{subequations}
since it avoids the need for additional storage (as in backward differencing), or the solution to a nonlinear system (as in temporal averaging).

The fully discrete system \eqref{eqn:mbe2} is a complete set of nonlinear equations defining the solution state at the new time, $t^{n+1}$. However, from the perspective of implementation, there is significant benefit in realizing that $\Pv^{n+1}_{m,\jv}$ are decoupled from other quantities at the next time level, i.e. $\Ev^{n+1}_{\jv}$ and $N^{n+1}_{\ell,\jv}$. As a result, they can be updated independently using \eqref{eqn:P-update-2} as 
\begin{align}\label{eqn:p-update-2}
  \Pv^{n+1}_{m,\jv}=\frac{1}{1+b_{1,m}\frac{\dt}{2}}\left(2\Pv^n_{m,\jv}-\Pv^{n-1}_{m,\jv} +b_{1,m} \frac{\dt}{2}\Pv^{n-1}_{m,\jv}-\dt^2b_{0,m}\Pv^n_{m,\jv} + \dt^2\sum_{\ell}a_{m,\ell}N_{\ell}^n\Ev^n_{\jv}\right).
\end{align}
Subsequently, $\Ev^{n+1}_{\jv}$ can be trivially determined from \eqref{eqn:E-update-2}. 
%
Finally, $N^{n+1}_{\ell,\jv}$ can be obtained using \eqref{eqn:n-update-2}, \eqref{eqn:2-accurate-Nt} and \eqref{eqn:2-accurate-Ntt}, 
where all terms on the right-hand-side of \eqref{eqn:n-update-2} are known because $\Pv^{n+1}_{m,\jv}$ and $\Ev^{n+1}_{\jv}$ have been previously computed.
%
This decoupling is a major difference from the schemes developed in~\cite{angel2019high,banks2020high}, where the solution update required the solution of a coupled system of linear equations locally at each grid cell $x_{\jv}$. In the present work this would translate to a nonlinear system, which may introduce numerical subtleties such as solver tolerances, choice of nonlinear root, etc.

\subsubsection{Fourth-order accurate scheme}
\label{sec:4th}
Following the ME approach, higher-order accurate approximations to \eqref{eqn:nonlin_e}--\eqref{eqn:nonlin_n} can be obtained by retaining additional correction terms in the Taylor expansions of the discrete temporal operators, e.g. \eqref{eqn:expansion2} and \eqref{eqn:expansion1}. Typically, the governing equations would be used to exchange temporal for spatial derivatives, or in the case of ODEs to successively reduce the order of temporal derivation. However as previously mentioned, this would lead to a globally coupled nonlinear system that would need to be solved at each time step. An alternative, first discussed in~\cite{angel2019high} for linear dispersive materials, uses predictions from lower-order schemes, e.g. \eqref{eqn:mbe2}, to approximate the correction terms to the requisite accuracy, and thereby enable a local explicit solution update at each grid point. This approach is dubbed HIME, for Hierarchical Modified Equation.
A pseudo-code algorithm for the fourth-order accurate Maxwell-MLA scheme is given in~\ref{sec:algorithmOrder4}.

To describe the fourth-order accurate HIME scheme, we first make some convenient notational definitions. Let the predicted approximation to the fields, polarization vectors, and carrier population densities at the new time $t^{n+1}$, as defined by the second-order accurate scheme \eqref{eqn:mbe2}, be denoted $\predE$, $\predP$, and $\predN$ respectively. Difference approximations based on these predictions will be colorized {\em \color{blue}blue} for clarity, and the difference operators will be applied to predicted ``star'' quantities, e.g. 
\begin{subequations}
\begin{align*}
  {\color{blue}D_{0t} \Ev^{n,*}_{\jv}}\eqdef\frac{\Ev^{n+1,*}_{\jv}-\Ev^{n-1}_{\jv}}{2\dt},\qquad
  {\color{blue}D_{+t}D_{-t}\Delta_{2h}\Pv^{n,*}_{\jv}} \eqdef \frac{\Delta_{2h}\Pv^{n+1,*}_{\jv}-2\Delta_{2h}\Pv^{n}_{\jv}+\Delta_{2h}\Pv^{n-1}_{\jv}}{\dt^2}.
\end{align*} 
\end{subequations}
With this notation, the fourth-order accurate HIME scheme, with any terms involving predicted values highlighted in {\em \color{blue}blue} for clarity, is 
\begin{subequations}\label{eqn:mbe4}
\begin{align}
  &D_{+t}D_{-t}\Ev^n_{\jv}
    -\frac{\dt^2}{12}\left(c^4\Delta_{2h}^2\Ev^n_{\jv}
      -\epsilon^{-1}_0c^2{\color{blue}D_{+t}D_{-t}\Delta_{2h}\Pv^n_{\jv}}
      -\epsilon^{-1}_0{\color{red}(D_{+t}D_{-t})^2\Pv^n_{\jv}}\right) \nonumber \\ 
  & \hspace{2in} = c^2\Delta_{4h}\Ev^n_{\jv}-\epsilon^{-1}_0D_{+t}D_{-t}\Pv^n_{\jv}
      +\frac{\dt^2}{12}\epsilon^{-1}_0{\color{red}(D_{+t}D_{-t})^2\Pv^n_{\jv}},\label{eqn:mbe4E}\\
  &D_{+t}D_{-t}\Pv^n_{m,\jv}-\frac{\dt^2}{12}{\color{blue}\left.D_{2tttt}\Pv_m\right|_{\jv}^n}+b_{1,m} \left(D_{0t}\Pv^n_{m,\jv}-\frac{\dt^2}{6}{\color{blue}\left.D_{2ttt}\Pv_m\right|_{\jv}^n}\right)  \nonumber\\
  & \hspace{2in} + b_{0,m}\Pv^n_{m,\jv} = \sum_{\ell}a_{m,\ell}N_{\ell,\jv}^n\Ev^n_{\jv},\\
  &D_{+t}N^{n}_{\ell,\jv}=\left.D_{4t}N_\ell\right|_{\jv}^n +\frac{\dt}{2} \left.D_{4tt}N_\ell\right|_{\jv}^n+\frac{\dt^2}{6} {\color{blue}\left.D_{2ttt}N_\ell\right|_{\jv}^n}+\frac{\dt^3}{24} {\color{blue}\left.D_{2tttt}N_\ell\right|_{\jv}^n},\label{eqn:n4}
\end{align}
\end{subequations}
for $m=1,\dots,\mathcal{N}_p$ and $\ell=0,1,\dots,\mathcal{N}_n-1$. Here the two {\em \color{red} red} terms are included in \eqref{eqn:mbe4E} because they naturally occur in the ME formulation, although they cancel and therefore need not appear in the final discretization. Further, $\Delta_{4h}$ denotes the fourth-order accurate discrete laplacian, $D_{4t}$ and $D_{4tt}$ indicate 4th-order accurate approximation to the continuous time derivatives $\partial_t$ and $\partial_t^2$ respectively, and $D_{2ttt}$ and $D_{2tttt}$ indicate 2nd-order accurate approximation to the continuous time derivatives $\partial_t^3$ and $\partial_t^4$ respectively (definitions of these terms are presented below). The various orders of accuracy for each term are consistent with the requirements for overall 4th-order accuracy of the scheme, as discussed for example in~\cite{max2006b,angel2019high,banks2020high}, and ultimately yield a fully 4th-order scheme in a compact spatial stencil using only three time levels. Definitions of approximations to the terms $\left.\partial^3_t\Pv_m\right|^n_{\jv}$, $\left.\partial^4_t\Pv_m\right|^n_{\jv}$, $\left.\partial_tN_{\ell}\right|^n_{\jv}$, $\left.\partial^2_tN_{\ell}\right|^n_{\jv}$, $\left.\partial^3_tN_{\ell}\right|^n_{\jv}$, and $\left.\partial^4_tN_{\ell}\right|^n_{\jv}$, as needed in \eqref{eqn:mbe4} appear in~\ref{sec:4thsupp}.


As in the case of the second-order discretization, the fourth-order HIME scheme permits a decoupled update of all quantities at the new time. The procedure is similar to second-order where first the polarization vectors are updated, then the fields, and finally the carrier population densities. Because this represents a significant advantage of the HIME versus traditional ME schemes for multilevel nonlinear electromagnetics, it is useful to describe this update in detail. After performing the predicted second-order update in a local stencil, the polarization vectors are updated as
\begin{align}
  \label{eqn:p-update-4}
  \Pv^{n+1}_{m,\jv} = \frac{1}{1+b_{1,m}\frac{\dt}{2}}
    \Big(
      2\Pv^n_{m,\jv} - \Pv^{n-1}_{m,\jv}
      +\frac{\dt^4}{12}{\color{blue}\left.D_{2tttt}\Pv^*_{m}\right|_{\jv}^n} 
      +\frac{\dt}{2}b_{1,m} \Pv^{n-1}_{m,\jv}\nonumber\\
      +\frac{\dt^4}{6}b_{1,m}{\color{blue}\left.D_{2ttt}\Pv^*_m\right|_{\jv}^n}
      -\dt^2b_{0,m}\Pv^n_{m,\jv} 
      +\dt^2\sum_{\ell}a_{m,\ell}N_{\ell,\jv}^n\Ev^n_{\jv}\Big).
\end{align}
Subsequently the electric fields can be updated as
\begin{align}
  \label{eqn:e-update-4}
  \Ev^{n+1}_{\jv}= 
    2\Ev^{n}_{\jv} - \Ev^{n-1}_{\jv}
    +\frac{\dt^4}{12}\left(
      c^4\Delta_{2h}^2\Ev^n_{\jv}
      -\epsilon^{-1}_0c^2{\color{blue}D_{+t}D_{-t}\Delta_{2h}\Pv^n_{\jv}}
    \right) \nonumber \\
    + \dt^2c^2\Delta_{4h}\Ev^n_{\jv}
    - \dt^2\epsilon^{-1}_0D_{+t}D_{-t}\Pv^n_{\jv}.
\end{align}
Finally the carrier populations are updated as 
\begin{align}
  \label{eqn:e-updatn-4}
  N^{n+1}_{\ell,\jv}=N^n_{\ell,\jv}+\dt \left.D_{4t}N_\ell\right|_{\jv}^n +\frac{\dt^2}{2} \left.D_{4tt}N_\ell\right|_{\jv}^n+\frac{\dt^3}{6} {\color{blue}\left.D_{2ttt}N_\ell\right|_{\jv}^n}+\frac{\dt^4}{24} {\color{blue}\left.D_{2tttt}N_\ell\right|_{\jv}^n}.
\end{align}



\subsection{Numerical interface approximations}\label{sec:interface}

We now proceed to a description of the numerical treatment of interface conditions \eqref{eqn:primary-ic} in the framework of HIME time-stepping. Throughout this section, we assume that 2nd- or 4th-order solution approximations have been time advanced to time level $t^{n}$ for all grid points on the domain interiors and along material interfaces. Furthermore, it is assumed that the grids across the interface are matched at point $\jv$ in the tangential direction, see Figure~\ref{fig:InterfaceCartoonOrder2Cartesian} and Figure~\ref{fig:InterfaceCartoonOrder4Cartesian}. The discrete interface conditions are then enforced using ghost cells, and the primary purpose of the present section is to describe how these ghost cells are determined at $t^{n}$. Once solution approximations in the ghost cells have been determined, subsequent time stepping using \eqref{eqn:mbe2} or \eqref{eqn:mbe4} will yield fully 2nd- or 4th-order accurate approximations. 
Note that to simplify presentation, the following discussion will be restricted to the case of Cartesian grids. 
The algorithms for the curvilinear case are very similiar and
are presented in~\ref{sec:InterfaceMLAOrder2Algorithm} and~\ref{sec:InterfaceMLAOrder4Algorithm}.

Determination of solution approximations in ghost cells naturally relies on interface conditions containing spatial derivatives, since undifferentiated terms constrain the solution directly on the interface. Primary interface conditions \eqref{eqn:ic-3} and \eqref{eqn:ic-4} already involve the requisite derivative operators, but \eqref{eqn:ic-1} and \eqref{eqn:ic-2} do not. Following the approach described in~\cite{max2006b,banks2020high}, the primary interface conditions \eqref{eqn:ic-1} and \eqref{eqn:ic-2} are therefore time differentiated, and the governing PDEs used to yield
\begin{subequations}
  \begin{align}
    [\nv\times \partial_t^2\Ev]_I &= [\nv\times(c^2\Delta \Ev-\epsilon^{-1}_0\partial^2_t\Pv)]_I=0,\\
    [\nv\cdot (\epsilon_0\partial_t^2\Ev+\partial_t^2\Pv)]_I &= [\nv\cdot (\epsilon_0c^2\Delta \Ev)]_I=0,
  \end{align}
\end{subequations}
which are used in place of \eqref{eqn:ic-1} and \eqref{eqn:ic-2}. The full set of primary interface conditions that are used to determine solution approximations in ghost cells are therefore
\begin{subequations}
  \label{eq:ic-primary}
  \begin{align}
    [\nv\times(c^2\Delta \Ev-\epsilon^{-1}_0\partial^2_t\Pv)]_I=0,\label{eqn:ic-1-primary}\\
    [\nv\cdot (\epsilon_0c^2\Delta \Ev)]_I=0,\label{eqn:ic-2-primary}\\
    \left[\mu_0^{-1}\nv\times\nabla\times\Ev\right]_I=0,\label{eqn:ic-3-primary}\\
    [\nabla\cdot\Ev]_I=0.\label{eqn:ic-4-primary}
  \end{align}
\end{subequations}

\subsubsection{Second-order accurate interface approximation} \label{sec:ic2}

{
\newcommand{\smallss}{\sffamily\small}
\newcommand{\mytrix}{(\x,-.15) -- ++(.3,0) -- ++(-.15,.26) -- (\x,-.15)}
\newcommand{\ghostMark}{circle (3pt)}
\newcommand{\tikzcircle}[2][red,fill=red]{\tikz[baseline=-0.5ex]\draw[#1,radius=#2] (0,0) circle ;}%
\newcommand{\tikzopencircle}[2][black,very thick]{\tikz[baseline=-0.5ex]\draw[#1,radius=#2] (0,0) circle ;}%
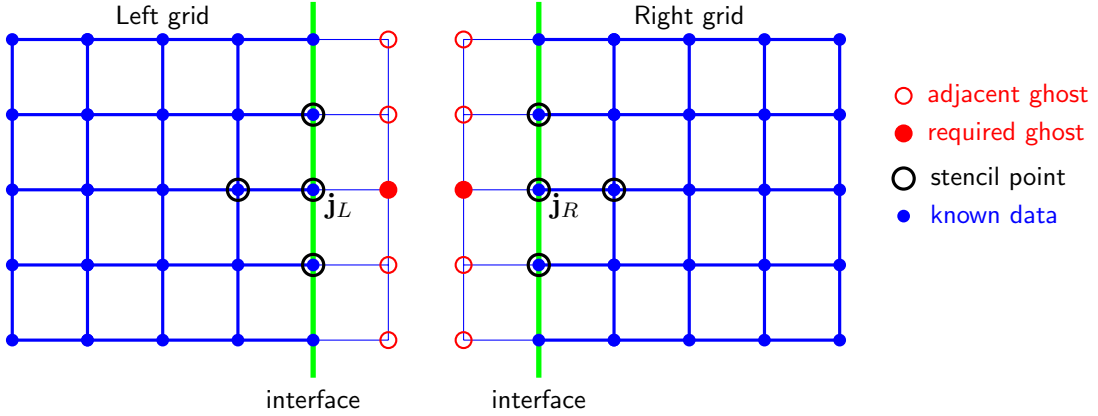
\begin{figure}[H]
\centering
\begin{tikzpicture}[scale=1]
   \useasboundingbox (.0,0.0) rectangle (13.5,5.5);  
  \begin{scope}[xshift=0cm,yshift=.75cm]

    \draw (2,4) node[above] {\smallss Left grid};
    \draw[step=1cm,blue,very thick] (0,0) grid (4,4);
    \draw[step=1cm,blue,thin] (4,0) grid (5,4);

    \draw[green,line width=2pt] (4,-.5) node[anchor=north,black] {\smallss interface} -- (4,4.5); 

    \foreach \x in {0,1,2,,3,4}
    {
       \foreach \y in {0,1,2,3,4}
       {
          \filldraw[blue,thick] (\x,\y) circle (2pt); 
       }
    }

    \foreach \y in {0,1,3,4}
      \draw[red,thick] (5,\y) \ghostMark;

     \filldraw[red,thick] (5,2) \ghostMark;

     \foreach \position in {(3,2),(4,2),(4,1),(4,3)}
        \draw[black,very thick] \position circle (4pt); 

     \draw (4,2) node[xshift=10pt,yshift=-6pt] {$\jv_L$};

  \end{scope}
  \begin{scope}[xshift=7cm,yshift=.75cm]

    \draw (2,4) node[above] {\smallss Right grid};
    \draw[step=1cm,blue,very thick] (0,0) grid (4,4);
    \draw[step=1cm,blue,thin] (-1,0) grid (0,4);
    \foreach \y in {0,1,3,4}
      \draw[red,thick] (-1,\y) \ghostMark;

     \filldraw[red,thick] (-1,2) \ghostMark;

     \draw[green,line width=2pt] (0,-.5) node[anchor=north,black] {\smallss interface} -- (0,4.5);   

      \foreach \x in {0,1,2,,3,4}
      {
         \foreach \y in {0,1,2,3,4}
         {
            \filldraw[blue,thick] (\x,\y) circle (2pt); 
         }
      }    

     \foreach \position in {(1,2),(0,2),(0,1),(0,3)}
        \draw[black,very thick] \position circle (4pt); 

     \draw (0,2) node[xshift=10pt,yshift=-6pt] {$\jv_R$};

  \end{scope}
  \begin{scope}[xshift=11.75cm,yshift=3cm]
     \draw[red,thick] (.1,1.) \ghostMark  node[anchor=west,xshift=5,yshift=0pt] {\smallss adjacent ghost};
     \filldraw[red,thick] (.1,.5) \ghostMark  node[anchor=west,xshift=5,yshift=0pt] {\smallss required ghost};
     \draw[black,very thick] (.1,-.1) circle (4pt)  node[anchor=west,xshift=6] {\smallss stencil point};
     \filldraw[blue,thick] (.1,-.6) circle (2pt)  node[anchor=west,xshift=6,yshift=1pt] {\smallss known data};
  \end{scope}
\end{tikzpicture}
\caption{Local stencil for filling in the interface ghost values for $\Ev$ for the second-order accurate scheme on Cartesian grids. 
   The points with indexes $\jv_L$ and $\jv_R$ correspond to a common physical point $\xv_{\jv_L}=\xv_{\jv_R}$ on the interface.
   The two required ghost values \tikzcircle{3pt} depend on the solution values at the stencil points \tikzopencircle{3pt}.
   There is no tangential coupling with adjacent ghost points and there are no nonlinear iterations required to update the ghost values.
}
\label{fig:InterfaceCartoonOrder2Cartesian}
\end{figure}
}

Discretization of the primary interface conditions \eqref{eq:ic-primary} to second-order accuracy is straightforward, and yields
\begin{subequations}
  \label{eqn:ic2}
  \begin{align}
    \left[\nv\times\left(c^2\Delta_{2h} \Ev^n_{\jv}-\epsilon^{-1}_0{\left.D_{+t}D_{-t}\Pv\right|^n_{\jv}}\right)\right]=0,\quad \jv\in\Gamma_h,\label{eqn:ic2-1}\\
    \left[\nv\cdot (\epsilon_0c^2\Delta_{2h} \Ev^n_{\jv}\right]=0,\quad \jv\in\Gamma_h,\label{eqn:ic2-2}\\
    \left[\mu_0^{-1}\nv\times\nabla_{2h}\times\Ev^n_{\jv}\right]=0,\quad \jv\in\Gamma_h,\label{eqn:ic2-3}\\
    \left[\nabla_{2h}\cdot\Ev^n_{\jv}\right]=0,\quad \jv\in\Gamma_h,\label{eqn:ic2-4}
  \end{align}
\end{subequations}
where $\Gamma_h$ indicate the set of indices along the material interface, and $\Delta_{2h}$ and $\nabla_{2h}$ denote second-order accurate finite differences.
 Note that the term ${\left.D_{+t}D_{-t}\Pv\right|^n_{\jv}}$ involves $\Pv_{\jv}^{n+1}$, i.e. $\Pv$ at a future time, which is not known. However, \eqref{eqn:p-update-2} can be used to determine $\Pv_{\jv}^{n+1}$ based on known information at $t_n$. Equivalently, the definition of $\Pv_{\jv}^{n+1}$ from \eqref{eqn:p-update-2} can be inserted directly in ${\left.D_{+t}D_{-t}\Pv\right|^n_{\jv}}$ to give
\begin{align}
  \label{eqn:Ptt2}
  {\left.D_{+t}D_{-t}\Pv\right|^n_{\jv}}=&\sum_{m=1}^{\mathcal{N}_p}\left(\frac{1}{\Delta t^2}+b_{1,m}\frac{1}{2\Delta t}\right)^{-1}\Big(2\Pv^n_{m,\jv}-\Pv^{n-1}_{m,\jv} +b_{1,m} \frac{\Delta t}{2}\Pv^{n-1}_{m,\jv}\nonumber\\
  &-\Delta t^2b_{0,m}\Pv^n_{m,\jv} + \Delta t^2a_{m,\ell}N_{\ell}^n\Ev^n_{\jv}\Big)-\frac{2}{\Delta t^2}\Pv^n_{\jv}+\frac{1}{\Delta t^2}\Pv^{n-1}_{\jv}.
\end{align} 
Clearly, the definition of $\left.D_{+t}D_{-t}\Pv\right|^n_{\jv}$ in \eqref{eqn:Ptt2} is not coupled with any unknown ghost values at time level $t^n$, and so the interface condition \eqref{eqn:ic2-1} can be expressed 
\begin{align}
  \left[\nv\times c^2\Delta_{2h} \Ev^n_{\jv}\right]_I=\left[\nv\times\epsilon^{-1}_0{\left.D_{+t}D_{-t}\Pv\right|^n_{\jv}}\right]_I,\quad \jv\in\Gamma_h\label{eqn:ic2-1-alt}\tag{\ref*{eqn:ic2-1}*},
\end{align} 
where the right-hand side is considered to be known data via \eqref{eqn:Ptt2}. 
The system of equations \eqref{eqn:ic2-1-alt} along with \eqref{eqn:ic2-2} -- \eqref{eqn:ic2-4} 
therefore defines values of $\Ev$ in the ghost cells (see Figure~\ref{fig:InterfaceCartoonOrder2Cartesian}).
 An additional advantageous property of this approach is that ghost values for $\Ev^n$ are decoupled from the ghost values for $\Pv^n$, a property that will be referred to as \texttt{EP-decoupling}. In fact no ghost values for $\Pv^n$ are needed in the actual update for this second-order scheme \eqref{eqn:mbe2}, although in practice extrapolation is used to define values in the ghost so that the entire grid function is defined.


\subsubsection{Fourth-order accurate interface approximation} \label{sec:ic4}

{
\newcommand{\smallss}{\sffamily\small}
\newcommand{\ghostMark}{circle (3pt)}
\newcommand{\tikzcircle}[2][red,fill=red]{\tikz[baseline=-0.5ex]\draw[#1,radius=#2] (0,0) circle ;}%
\newcommand{\tikzopencircle}[2][black,very thick]{\tikz[baseline=-0.5ex]\draw[#1,radius=#2] (0,0) circle ;}%
\begin{figure}[H]
\centering
\begin{tikzpicture}[scale=1]
   \useasboundingbox (.0,0.0) rectangle (13.5,5.5);  
  \begin{scope}[xshift=-1cm,yshift=.75cm]

    \draw (2,4) node[above] {\smallss Left grid};
    \draw[step=1cm,blue,very thick] (0,0) grid (4,4);
    \draw[step=1cm,blue,thin] (4,0) grid (6,4);

    \draw[green,line width=2pt] (4,-.5) node[anchor=north,black] {\smallss interface} -- (4,4.5); 

    \foreach \x in {0,1,2,,3,4}
    {
       \foreach \y in {0,1,2,3,4}
       {
          \filldraw[blue,thick] (\x,\y) circle (2pt); 
       }
    }

    \foreach \y in {0,1,3,4}
    {
      \draw[red,thick] (5,\y) \ghostMark;
      \draw[red,thick] (6,\y) \ghostMark;
    }
    \foreach \x in {5,6}
      \filldraw[red,thick] (\x,2) \ghostMark;
    \foreach \position in {(4,0),(4,1),(4,2),(4,3),(4,4), (3,0),(3,1),(3,2),(3,3),(3,4),(2,1),(2,2),(2,3)}
       \draw[black,very thick] \position circle (4pt); 
    \draw (4,2) node[xshift=10pt,yshift=-6pt] {$\jv_L$};

  \end{scope}
  \begin{scope}[xshift=8.25cm,yshift=.75cm]

    \draw (2,4) node[above] {\smallss Right grid};
    \draw[step=1cm,blue,very thick] (0,0) grid (4,4);
    \draw[step=1cm,blue,thin] (-2,0) grid (0,4);
    \foreach \y in {0,1,3,4}
    {
      \draw[red,thick] (-1,\y) \ghostMark;
      \draw[red,thick] (-2,\y) \ghostMark;
     }

    \draw[green,line width=2pt] (0,-.5) node[anchor=north,black] {\smallss interface} -- (0,4.5);   
    
     \foreach \x in {-2,-1}
       \filldraw[red,thick] (\x,2) \ghostMark;

      \foreach \x in {0,1,2,,3,4}
      {
         \foreach \y in {0,1,2,3,4}
         {
            \filldraw[blue,thick] (\x,\y) circle (2pt); 
         }
      }    

     \foreach \position in {(0,0),(0,1),(0,2),(0,3),(0,4), (1,0),(1,1),(1,2),(1,3),(1,4),   (2,1),(2,2),(2,3)}
        \draw[black,very thick] \position circle (4pt); 

     \draw (0,2) node[xshift=10pt,yshift=-6pt] {$\jv_R$};

  \end{scope}
  \begin{scope}[xshift=12.75cm,yshift=2.75cm]
     \draw[red,thick] (.1,1.) \ghostMark  node[anchor=west,xshift=5,yshift=0pt] {\smallss adjacent ghost};
     \filldraw[red,thick] (.1,.5) \ghostMark  node[anchor=west,xshift=5,yshift=0pt] {\smallss required ghost};
     \draw[black,very thick] (.1,-.1) circle (4pt)  node[anchor=west,xshift=6] {\smallss stencil point};
     \filldraw[blue,thick] (.1,-.6) circle (2.5pt)  node[anchor=west,xshift=6,yshift=1pt] {\smallss known data};     
  \end{scope}
\end{tikzpicture}
\caption{Local stencil for filling in the interface ghost values for $\Ev$ for the fourth-order accurate scheme on Cartesian grids. 
   The points with indexes $\jv_L$ and $\jv_R$ correspond to a common physical point $\xv_{\jv_L}=\xv_{\jv_R}$ on the interface.
   The four required ghost values \tikzcircle{3pt} depend on the solution values at the stencil points \tikzopencircle{3pt}.
   There is no tangential coupling with adajcent ghost points and there are no nonlinear iterations required to update the ghost values.
}
\label{fig:InterfaceCartoonOrder4Cartesian}
\end{figure}
}

As we employ fourth-order accurate central finite difference scheme, we have two ghost lines at both sides of the interface, which requires additional four jump conditions obtained by differentiating the four primary interface conditions \eqref{eqn:primary-ic} in time:
\begin{subequations}
\begin{align}
[\nv\times\partial^4_t\Ev]_I=0,\\
[\nv\cdot (\epsilon_0\partial^4_t\Ev+\partial^4_t\Pv)]_I=0,\\
\left[\mu_0^{-1}\nv\times\nabla\times\partial^2_t\Ev\right]_I=0,\\
[\nabla\cdot\partial^2_t\Ev]_I=0,
\end{align}
\end{subequations}
then we use the compatibility condition again and obtain the following additional four interface conditions
\begin{subequations}\label{eqn:add-ic}
\begin{align}
[\nv\times(c^4\Delta^2 \Ev-c^2\epsilon^{-1}_0\Delta \partial^2_t\Pv-\epsilon^{-1}_0\partial^4_t\Pv)]_I=0,\label{eqn:add-ic-1}\\
[\nv\cdot \epsilon_0(c^4\Delta^2 \Ev-c^2\Delta \partial^2_t\Pv)]_I=0,\label{eqn:add-ic-2}\\
\left[\mu_0^{-1}\nv\times\nabla\times(c^2\Delta \Ev-\epsilon^{-1}_0\partial^2_t\Pv)\right]_I=0,\label{eqn:add-ic-3}\\
[\nabla\cdot (c^2\Delta \Ev-\epsilon^{-1}_0\partial^2_t\Pv)]_I=0.\label{eqn:add-ic-4}
\end{align}
\end{subequations}

Fourth-order accurate approximation of the primary interface conditions \eqref{eqn:primary-ic} and the additional four jump conditions \eqref{eqn:add-ic} are given by
\begin{subequations}\label{eqn:ic4}
\begin{align}
  \left[\nv\times(c^2\Delta_{4h} \Ev^n_{\jv}-\epsilon^{-1}_0{\color{magenta}\left.D_{4tt}\Pv\right|^n_{\jv}})\right]_I=0,\quad \jv\in\Gamma_h,\label{eqn:ic4-1}\\
  \left[\nv\cdot (\epsilon_0c^2\Delta_{4h} \Ev^n_{\jv})\right]_I=0,\quad \jv\in\Gamma_h,\label{eqn:ic4-2}\\
  \left[\mu_0^{-1}\nv\times\nabla_{4h}\times\Ev^n_{\jv}\right]_I=0,\quad \jv\in\Gamma_h,\label{eqn:ic4-3}\\
  \left[\nabla_{4h}\cdot\Ev^n_{\jv}\right]_I=0,\quad \jv\in\Gamma_h,\label{eqn:ic4-4}\\
  \left[\nv\times(c^4\Delta^2_{2h} \Ev^n_{\jv}-c^2\epsilon^{-1}_0{\color{blue}\left.D_{+t}D_{-t}\Delta_{2h}\Pv\right|^n_{\jv}}-\epsilon^{-1}_0{\color{blue}\left.D_{2tttt}\Pv\right|^n_{\jv}})\right]_I=0,\quad \jv\in\Gamma_h,\label{eqn:ic4-5}\\
  \left[\nv\cdot \epsilon_0(c^4\Delta^2_{2h} \Ev^n_{\jv}-c^2{\color{blue}\left.D_{+t}D_{-t}\Delta_{2h}\Pv\right|^n_{\jv}})\right]_I=0,\quad \jv\in\Gamma_h,\label{eqn:ic4-6}\\
  \left[\mu_0^{-1}\nv\times(c^2\nabla_{2h}\times\Delta_{2h}\Ev^n_{\jv}-\epsilon^{-1}_0{\color{blue}\left.D_{+t}D_{-t}\nabla_{2h}\times\Pv\right|^n_{\jv}})\right]_I=0,\quad \jv\in\Gamma_h,\label{eqn:ic4-7}\\
  \left[\nabla_{2h} \cdot c^2\Delta_{2h}\Ev^n_{\jv}-\epsilon^{-1}_0{\color{blue}\left.D_{+t}D_{-t}\nabla_{2h} \cdot \Pv\right|^n_{\jv}}\right]_I=0,\quad \jv\in\Gamma_h,\label{eqn:ic4-8}
\end{align}
\end{subequations}
where the subscript $2h$ and $4h$ denote the second-order and fourth-order accurate discretizations respectively. The magenta term $\left.D_{4tt}\Pv\right|^n_{\jv}$ requires fourth-order accurate approximation, while the blue terms need only second-order accurate approximations. 

\begin{remark}
In 3D, the interface conditions \eqref{eqn:ic4} at point $\jv$ can be rewritten into the vector form:
\begin{subequations}\label{eqn:ic4-3d}
{\small
\begin{align}
\left[(\nabla_{4h}\cdot\Ev^n_{\jv})\nv+(\Iv-\nv\nv^T)\mu_0^{-1}\nabla_{4h}\times\Ev^n_{\jv}\right]_I&=0,\label{eqn:ic4-3d-1}\\
\left[\nv\nv^T(\epsilon_0c^2\Delta_{4h} \Ev^n_{\jv})+(\Iv-\nv\nv^T)(c^2\Delta_{4h} \Ev^n_{\jv}-\epsilon^{-1}_0{\color{magenta}\left.D_{4tt}\Pv\right|^n_{\jv}})\right]_I&=0,\label{eqn:ic4-3d-2}\\
\left[\nabla_{2h} \cdot (c^2\Delta_{2h}\Ev^n_{\jv}-\epsilon^{-1}_0{\color{blue}\left.D_{+t}D_{-t}\Pv\right|^n_{\jv}})\nv+(\Iv-\nv\nv^T)\mu_0^{-1}\nabla_{2h}\times(c^2\Delta_{2h}\Ev^n_{\jv}-\epsilon^{-1}_0{\color{blue}\left.D_{+t}D_{-t}\Pv\right|^n_{\jv}})\right]_I&=0,\label{eqn:ic4-3d-3}\\
\big[\nv\nv^T\epsilon_0\Delta_{2h}(c^4\Delta_{2h} \Ev^n_{\jv}-c^2{\color{blue}\left.D_{+t}D_{-t}\Pv\right|^n_{\jv}})&\nonumber\\
+(\Iv-\nv\nv^T)(c^4\Delta^2_{2h} \Ev^n_{\jv}-c^2\epsilon^{-1}_0{\color{blue}\left.D_{+t}D_{-t}\Delta_{2h}\Pv\right|^n_{\jv}}-\epsilon^{-1}_0{\color{blue}\left.D_{2tttt}\Pv\right|^n_{\jv}})\big]_I&=0,\label{eqn:ic4-3d-4}
\end{align}
}
\end{subequations}
by combining the tangential and normal components, which allows convenient implementation of the jump conditions. Here $\nv$ is the unit normal vector at point $\jv$, and $\Iv$ is an identity matrix.
\end{remark}

Note that, as written, the numerical interface conditions~\eqref{eqn:ic4} have some un-desirable characteristics.
Firstly, due to cross-derivative terms such as  $\p_x\p_y u$ in $\Delta^2u$, the discrete interface conditions couple ghost points
in the tangential directions; this would require solution of a system of equations along the entire interface.
Secondly, there is no \texttt{EP-decoupling}.
To see this, we take the example of the blue term $\left.D_{+t}D_{-t}\Delta_{2h}\Pv\right|^n_{\jv}$ and note that:
\begin{align}
\left.D_{+t}D_{-t}\Delta_{2h}\Pv\right|^n_{\jv} =& \sum_{m=1}^{\mathcal{N}_p}\left(\frac{1}{\Delta t^2}+b_{1,m}\frac{1}{2\Delta t}\right)^{-1}\Big(2\Delta_{2h}\Pv^n_{m,\jv}-\Delta_{2h}\Pv^{n-1}_{m,\jv} +b_{1,m} \frac{\Delta t}{2}\Delta_{2h}\Pv^{n-1}_{m,\jv}\nonumber\\
&-\Delta t^2b_{0,m}\Delta_{2h}\Pv^n_{m,\jv} + \Delta t^2a_{m,\ell}{\color{red}\Delta_{2h}(N_{\ell,\jv}^n\Ev^n_{\jv})}\Big)-\frac{2}{\Delta t^2}\Delta_{2h}\Pv^n_{\jv}+\frac{1}{\Delta t^2}\Delta_{2h}\Pv^{n-1}_{\jv}.
\end{align}
The red term $\Delta_{2h}(N_{\ell,\jv}^n\Ev^n_{\jv})$ obviously has time-dependent variable coefficients for the unknown first ghost point lines (ghost points that are one grid away from the interface), which implies that we need to formulate jump conditions \eqref{eqn:ic4-5} and \eqref{eqn:ic4-6} at each time step. 
Similarly, for terms $ \left.D_{+t}D_{-t}\nabla_{2h}\times\Pv\right|^n_{\jv}$ and $ \left.D_{+t}D_{-t}\nabla_{2h} \cdot \Pv\right|^n_{\jv}$, variable-coefficient terms $\nabla_{2h}\times(N_{\ell}^n\Ev^n_{\jv})$ and $\nabla_{2h}\cdot(N_{\ell}^n\Ev^n_{\jv})$ persists respectively, thus \eqref{eqn:ic4-7} and \eqref{eqn:ic4-8} need to be reformulated at each time step accordingly as well.

Moreover, the magenta term ${\color{magenta}\left.D_{4tt}\Pv\right|_{\jv}^n}$ needs to be fourth-order accurate:
\begin{align}
{\color{magenta}\left.D_{4tt}\Pv\right|^n_{\jv}} =& \sum_{m=1}^{\mathcal{N}_p}\frac{1}{\Delta t^2+b_{1,m}\frac{\Delta t^3}{2}}\Big(2\Pv^n_{m,\jv}-\Pv^{n-1}_{m,\jv}+\frac{\Delta t^4}{12}\left.D_{2tttt}\Pv_m\right|^n_{\jv} +\frac{\Delta t}{2}b_{1,m} \Pv^{n-1}_{m,\jv}\nonumber\\
&\qquad\qquad\qquad\qquad\qquad+\frac{\Delta t^4}{6}\left.D_{2ttt}\Pv_m\right|^n_{\jv}-\Delta t^2b_{0,m}\Pv^n_{m,\jv} + \Delta t^2a_{m,\ell}N_{\ell}^n\Ev^n_{\jv}\Big)\nonumber\\
&-\frac{2}{\Delta t^2}\Pv^n_{\jv}+\frac{1}{\Delta t^2}\Pv^{n-1}_{\jv}-\frac{\Delta t^2}{12}\sum_{m=1}^{\mathcal{N}_p}\left.D_{2tttt}\Pv_m\right|^{n}_{\jv},
\end{align}
where the second-order accurate terms $\left.D_{2ttt} \Pv_m\right|^n_{\jv}$ and $\left.D_{2tttt} \Pv_m\right|^n_{\jv}$ are approximated by:
\begin{align}
\left.D_{2ttt} \Pv_m\right|^n_{\jv} &= -b_{1,m}D_{+t}D_{-t} \Pv^{n,*}_{m,\jv} - b_{0,m}D_{0t}\Pv^{n,*}_{m,\jv} + a_{m,\ell}\left.D_{2t} N_{\ell}\right|^{n,*}_{\jv} \Ev^n_{\jv} +a_{m,\ell}N^n_{\ell,\jv} \Big(2\Ev^n_{\jv}-2\Ev^{n-1}_{\jv} \nonumber\\
&+\Delta t^2c^2{\color{red}\Delta_{2h}\Ev^n_{\jv}}-\epsilon^{-1}_0\sum_{m=1}^{\mathcal{N}_p}\left(\Pv^{n+1}_{m,\jv}-2\Pv^n_{m,\jv}+\Pv^{n-1}_{m,\jv}\right)\Big),
\end{align}
and
\begin{align}
\left.D_{2tttt} \Pv_m\right|^n_{\jv} = & b_{1,m}b_{0,m}D_{0t}\Pv^{n,*}_{m,\jv}+(b_{1,m}^2- b_{0,m})D_{+t}D_{-t} \Pv^{n,*}_{m,\jv}\nonumber\\
& + (a_{m,\ell}\left.D_{2tt} N_{\ell}\right|^{n,*}_{\jv}-b_{1,m} a_{m,\ell}\left.D_{2t} N_{\ell}\right|^{n,*}_{\jv})\Ev^n_{\jv} \nonumber\\
&+\left(\frac{a_{m,\ell}N^n_{\ell,\jv}}{\Delta t^2}+\frac{2a_{m,\ell}\left.D_{2t}N_{\ell}\right|^{n,*}_{\jv}-b_{1,m}a_{m,\ell}N^n_{\ell,\jv}}{2\Delta t}\right)\Delta t^2c^2{\color{red}\Delta_{2h}\Ev^n_{\jv}}\nonumber\\
&-\left(\frac{a_{m,\ell}N^n_{\ell,\jv}}{\Delta t^2}+\frac{2a_{m,\ell}\left.D_{2t}N_{\ell}\right|^{n,*}_{\jv}-b_{1,m}a_{m,\ell}N^n_{\ell,\jv}}{2\Delta t}\right)
\epsilon^{-1}_0(\Pv^{n+1}_{\jv}-2\Pv^n_{\jv}+\Pv^{n-1}_{\jv})\nonumber\\
&+\frac{(2a_{m,\ell}\left.D_{2t}N_{\ell}\right|^{n,*}_{\jv}-b_{1,m}a_{m,\ell}N^n_{\ell,\jv})}{2\Delta t}(2\Ev^n_{\jv}-2\Ev^{n-1}_{\jv}).
\end{align}
Again, we see that $\left.D_{2ttt} \Pv_m\right|^n_{\jv}$ and $\left.D_{2tttt} \Pv_m\right|^n_{\jv}$ have time dependent coefficients for $\Ev^n$ on the first ghost lines, which implies that the interface condition \eqref{eqn:ic4-1} needs to be reformulated at each time step.

Fortunately, the colored terms $\left.D_{+t}D_{-t}\Delta_{2h}\Pv\right|^n_{\jv}$, $\left.D_{+t}D_{-t}\nabla_{2h}\times \Pv\right|^n_{\jv}$, $\left.D_{+t}D_{-t}\nabla_{2h} \cdot \Pv\right|^n_{\jv}$, $\left.D_{2tttt}\Pv\right|_{\jv}^{n}$ and $\left.D_{4tt}\Pv\right|_{\jv}^{n}$ only nonlinearly depend on $\Ev^{n}$ on the first ghost lines. To this end, we approximate the red terms using second-order accurate predicted values on the first ghost lines that are obtained using the second-order numerical interface conditions \eqref{eqn:ic2}, which leads to the linear system:
\begin{subequations}\label{eqn:ic4-alt}
\begin{align}
\left[\nv\times(c^2\Delta_{4h} \Ev^n_{\jv})\right]_I=\left[\nv\times(\epsilon^{-1}_0{\color{magenta}\left.D_{4tt}\Pv\right|^n_{\jv}})\right]_I,\quad \jv\in\Gamma_h,\label{eqn:ic4-1-alt}\tag{\ref*{eqn:ic4-1}*}\\
[\nv\times(c^4\Delta^2_{2h} \Ev^n_{\jv})]_I=[\nv\times(c^2\epsilon^{-1}_0{\color{blue}\left.D_{+t}D_{-t}\Delta_{2h} \Pv\right|^n_{\jv}}+\epsilon^{-1}_0{\color{blue}\left.D_{2tttt}\Pv\right|^n_{\jv}})]_I,\quad \jv\in\Gamma_h,\label{eqn:ic4-5-alt}\tag{\ref*{eqn:ic4-5}*}\\
[\nv\cdot (\epsilon_0c^4\Delta^2_{2h} \Ev^n_{\jv})]_I=[\nv\cdot (\epsilon_0c^2{\color{blue}\left.D_{+t}D_{-t}\Delta_{2h}\Pv\right|^n_{\jv}})]_I,\quad \jv\in\Gamma_h,\label{eqn:ic4-6-alt}\tag{\ref*{eqn:ic4-6}*}\\
\left[\mu_0^{-1}\nv\times(c^2\nabla_{2h}\times\Delta_{2h}\Ev^n_{\jv})\right]_I=\left[\mu_0^{-1}\nv\times(\epsilon^{-1}_0{\color{blue}\left.D_{+t}D_{-t}\nabla_{2h}\times \Pv\right|^n_{\jv}})\right]_I,\quad \jv\in\Gamma_h,\label{eqn:ic4-7-alt}\tag{\ref*{eqn:ic4-7}*}\\
[\nabla_{2h} \cdot c^2\Delta_{2h}\Ev^n_{\jv}]_I=[\epsilon^{-1}_0{\color{blue}\left.D_{+t}D_{-t}\nabla_{2h} \cdot\Pv\right|^n_{\jv}}]_I,\quad \jv\in\Gamma_h.\label{eqn:ic4-8-alt}\tag{\ref*{eqn:ic4-8}*}
\end{align}
\end{subequations}
Thus, we obtain the \texttt{EP-decoupling} for the fourth-order accurate approximations of the interface conditions similarly as the second-order case.



Moreover, as discussed in the previous work \cite{banks2020high}, cross terms in $\Delta^2_{2h}\Ev_{\jv}^n,\nabla_{2h}\times\Delta_{2h}\Ev_{\jv}^n,\nabla_{2h}\times\Delta_{2h}\Ev_{\jv}^n$ will invoke instability if the grid size in the tangential direction is smaller than the grid size in the normal direction (the instability will be investigated and dedicated to a separate paper). To avoid such instability, we follow \cite{banks2020high} and approximate the cross terms by decoupling $\Ev^{n}$ in the tangential direction, i.e., the second-order accurate predicted values on the first ghost lines are also employed here to approximate the cross terms. 
As a result of \texttt{tangential decoupling} that leads the stencil localization, the interface conditions \eqref{eqn:ic4} can be formulated into an identical time-independent local small linear system at each point $\jv$ independently rather than a large global system.

To summarize, for the fourth-order accurate interface conditions \eqref{eqn:ic4}, we linearize the jump conditions via \texttt{EP-decoupling} and localize the stencils using \texttt{tangential decoupling} by employing second-order accurate predicted values of $\Ev^{n}$ on the first ghost lines, obtained using the second-order approximations of interface conditions \eqref{eqn:ic2}.


\section{Numerical results}\label{sec:numerics}

This section presents results to verify the accuracy and stability of the new Maxwell-MLA schemes and interface approximations.
Section~\ref{sec:planarInteraceResults} studies the accuracy for manufactured solutions for a two-domain problem with a planar interface between two nonlinear materials.
Section~\ref{sec:curvedInteraceResults} then studies the accuracy for manufactured solutions on a problem with a curved interface
and overset grids. Section~\ref{sec:soliton} considers the propagation of a soliton for which there is 
an approximate asymptotic solution. A grid self-convergence study is used to evaluate the convergence rate of the computed soliton.
A problem with multiple interfaces is studied in Section~\ref{sec:layeredDisk} and shows results for
an interface between a nonlinear medium and a linear dispersive medium.
Finally in Section~\ref{sec:ellipsoidArray}, results are shown from the scattering of a Gaussian
plane wave from an array of ellipsoids composed from active materials.



\subsection{Planar interface results}\label{sec:planarInteraceResults}




To verify the accuracy of the interface implementations we perform a grid-refinement convergence study for problems with a planar interface. 
The domain for the interface problem in $d$-dimensions consists of two squares (cubes), 
\ba 
  \Omega = [-1,0]\times[0,1]^{d-1} \, \cup \, [0,1]\times[0,1]^{d-1}. 
\ea
Each square (cube) is covered by a Cartesian grid with grid spacing $1/(10j)$, $j=1,2,\ldots$. 
Let $\Gc_{Id}^{(j)}$ denote the composite grid of resolution $j$ for this domain.

\newcommand{\freqx}{f_{m,x}}
\newcommand{\freqy}{f_{m,y}}
\newcommand{\freqz}{f_{m,z}}
\newcommand{\freqt}{f_{m,t}}
Manufactured solutions are used to generate a known solution. 
The manufactured solution is defined using trigonometric functions and takes the form
\ba \label{eqn:manufactured}
    U_m = a_m \, \cos(\freqx \, x + \phi_{m,x}) \, \cos( \freqy \, y + \phi_{m,y}) \, \cos( \freqz \, z + \phi_{m,z}) \, \cos(\freqt \,t ), 
\ea
where $U_m$ denotes any component of the solution (e.g. a component of $\Ev$, $\Pv_j$ or $N_\ell$).
The amplitudes $a_m$, frequencies $\freqx$, $\freqy$ and phases $\phi_{m,x}$ are chosen differently for different $m$. For convenience the
manufactured solution for $\Ev$ is chosen to be divergence free as this simplifies implementations of the boundary conditions.
The initial conditions and boundary conditions are set to the known solution. The discrete solutions at the interface
are obtained with the numerical interface conditions as discussed in Section~\ref{sec:interface}.
The MLA materials in the left and right domains are chosen to be \texttt{mlaMat2} and \texttt{mlaMat3} as
defined in~\ref{sec:materialDefinitions}. These materials have different numbers of polarization vectors and population densities.
 


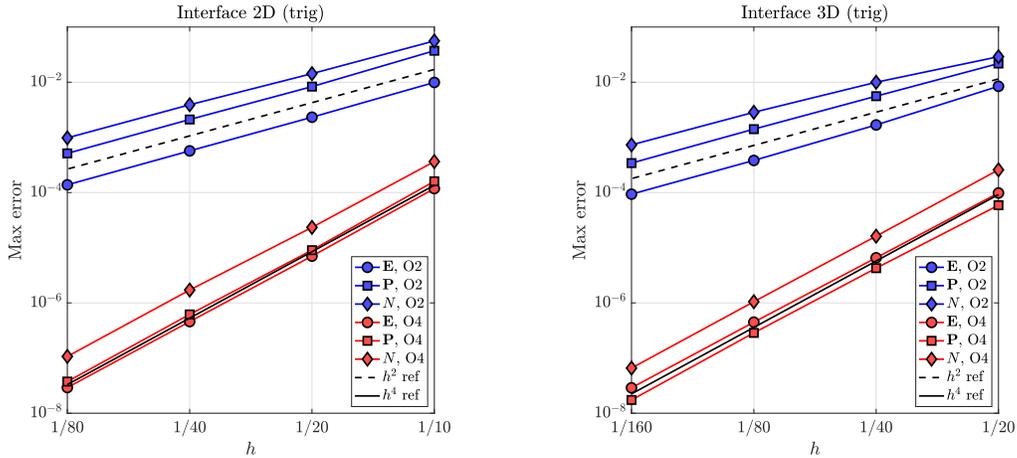
\begin{figure}[htb]
\begin{center}
\begin{tikzpicture}
   \useasboundingbox (0,.75) rectangle (14,6);  
  \figByHeight{  0}{0}{fig/interfaceMLA2d}{6cm}[0][0][0.][0]
  \figByHeight{7.5}{0}{fig/interfaceMLA3d}{6cm}[0][0][0.][0]
\end{tikzpicture}
\end{center}
\caption{Convergence for planar interfaces with manufactured solutions at $t=1.0$. Left: two dimensional results. Right: three dimensional results.}
\label{fig:interface2dTrig}
\end{figure}

Figure~\ref{fig:twoSquaresInterface} gives an example of the Cartesian meshes for 2D interface (left), approximated solution $E_y$ (middle) and the errors against the exact solution for $E_y$ component (right) at $t=1$ on meshes with grid spacing $h=1/160$. Figure~\ref{fig:interface2dTrig} shows the convergence results with manufactured solutions~\eqref{eqn:manufactured} for 2D (left) and 3D (right), with errors measured as the maximum norm in space.
The max-norm of a vector quantity is the maximum over all grid points of the maximum absolute value of all components of the vector.

{
\newcommand{\drawContourB}[7]{%
\begin{scope}[#1]
\draw(0.0,0) node[anchor=south west,xshift=-4pt,yshift=+0pt] {\trimfigb{fig/#2}{\figWidth}};
  \draw(.75,.5) node[draw,fill=white,anchor=west,xshift=2pt,yshift=1pt] {\scriptsize #3};
\begin{scope}[xshift=-.1cm,yshift=-3pt]
  \draw (\xcb,\ycb) node[anchor=south west,xshift=0.25cm,yshift=.5cm,rotate=-90] {\trimfigcb{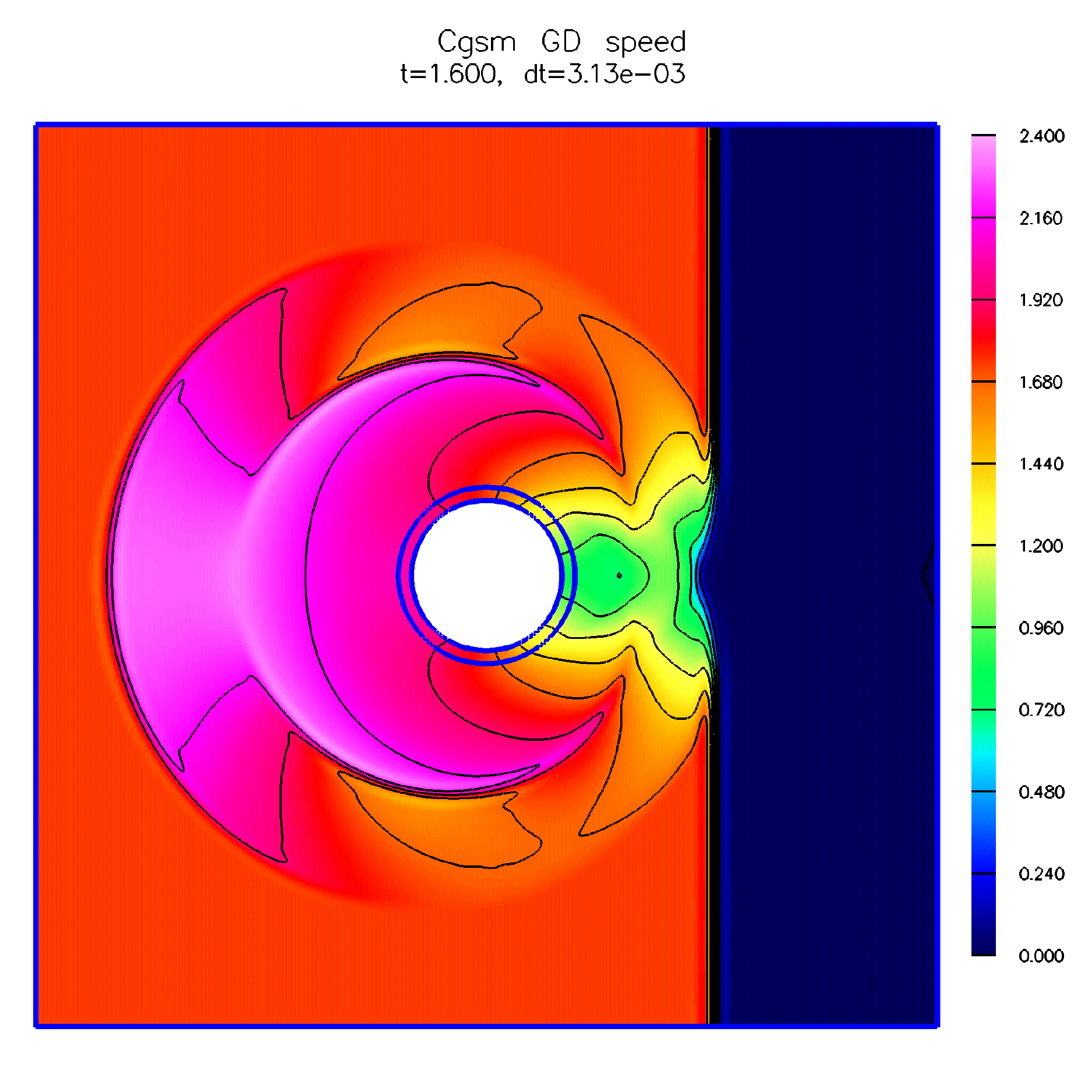}{\cbWidth}{\cbHeight}};
  \draw (.8,0) node[anchor=north,xshift=+3pt,yshift=+2pt] {\scriptsize $#6$};
  \draw (4.8,0) node[anchor=north,xshift=+0pt,yshift=+2pt] {\scriptsize $#7$};
\end{scope}
\end{scope}
}
\newcommand{\drawContourC}[7]{%
\begin{scope}[#1]
\draw(0.0,0) node[anchor=south west,xshift=-4pt,yshift=+0pt] {\trimfigc{fig/#2}{\figWidth}};
  \draw(.75,.5) node[draw,fill=white,anchor=west,xshift=2pt,yshift=1pt] {\scriptsize #3};
\begin{scope}[xshift=-.1cm,yshift=-3pt]
  \draw (\xcb,\ycb) node[anchor=south west,xshift=0.25cm,yshift=.5cm,rotate=-90] {\trimfigcb{fig/colourBarLines}{\cbWidth}{\cbHeight}};
  \draw (.8,0) node[anchor=north,xshift=+3pt,yshift=+2pt] {\scriptsize $#6$};
  \draw (4.8,0) node[anchor=north,xshift=+0pt,yshift=+2pt] {\scriptsize $#7$};
\end{scope}
\end{scope}
}
\newcommand{\cbWidth}{.2cm}
\newcommand{\cbHeight}{4cm}
\newcommand{\xcb}{.5cm}
\newcommand{\ycb}{-.2cm}
\setlength{\ycbTop}{\ycb+\cbHeight}
\setlength{\ycbMid}{\ycb+\cbHeight*\real{.5}}
\newcommand{\trimfigcb}[3]{\includegraphics[width=#2, height=#3, clip, trim=17cm 2.35cm 1.65cm 2.35cm]{#1}}
%
%
\def\rad{1.27}
\newcommand{\figWidth}{4.9cm}
\newcommand{\trimfig}[2]{\trimh{#1}{#2}{.02}{.02}{.02}{.02}}
\newcommand{\figWidtha}{5.75cm}
\newcommand{\trimfiga}[2]{\trimh{#1}{#2}{.04}{.12}{.07}{.09}}
\newcommand{\trimfigb}[2]{\trimh{#1}{#2}{.195}{.01}{.15}{.4}}
\newcommand{\trimfigc}[2]{\trimh{#1}{#2}{.18}{.18}{.15}{.4}}
\begin{figure}[htb]
\begin{center}
\resizebox{15.25cm}{!}{
\begin{tikzpicture}[scale=1]
  \useasboundingbox (.4,0.2) rectangle (16.5,4.75);  
   \draw(0.,0) node[anchor=south west,xshift=-4pt,yshift=0.5cm] {\includegraphics[width=6.25cm]{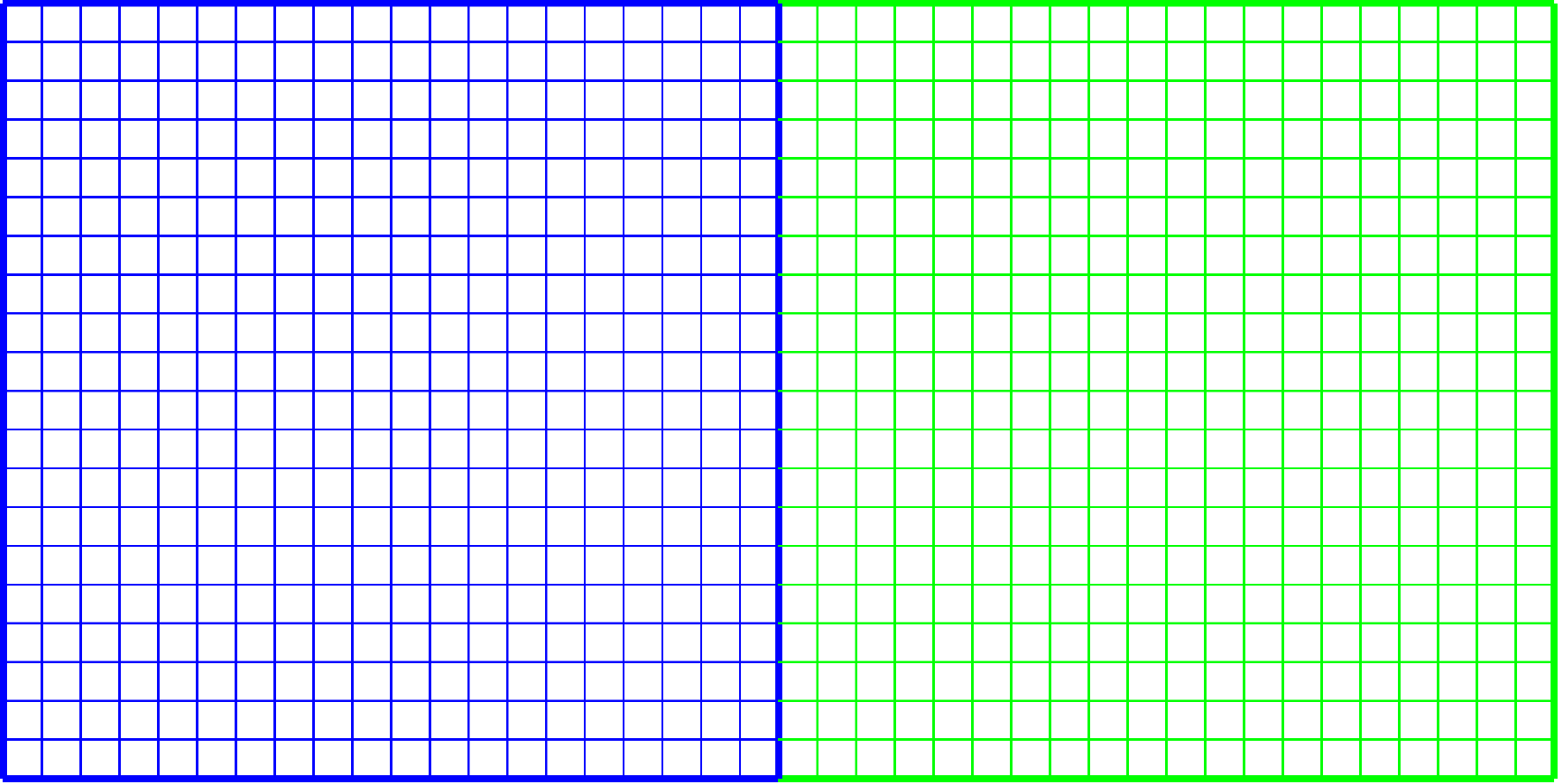}};
   \drawContourC{xshift= 5.9cm,yshift=-.08cm}{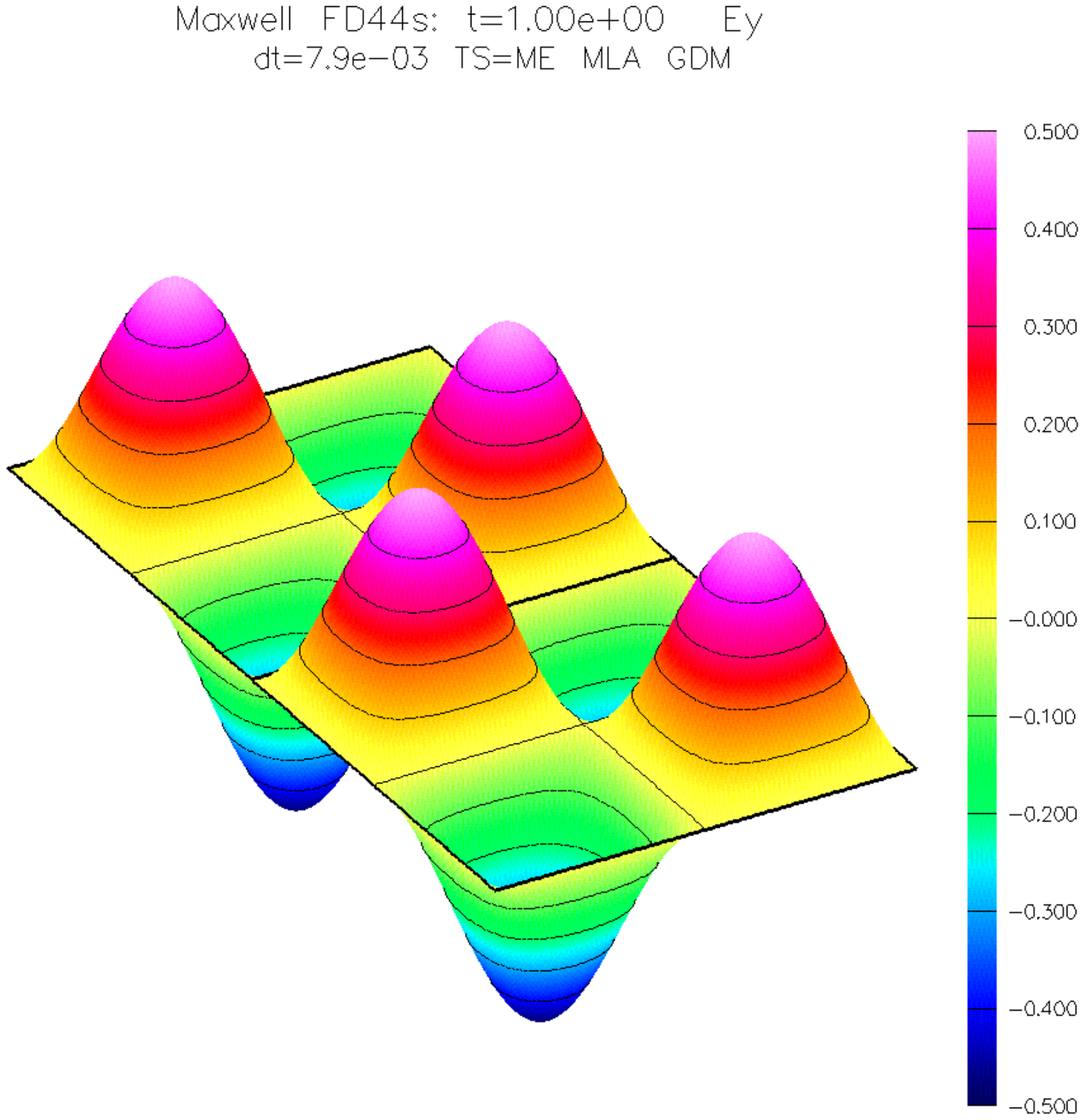}{$E_y$ at $t=1.0$}{$\Ev$}{$t=0.3$}{$-0.5$}{$0.5$};
   \drawContourC{xshift=11.cm,yshift=-.08cm}{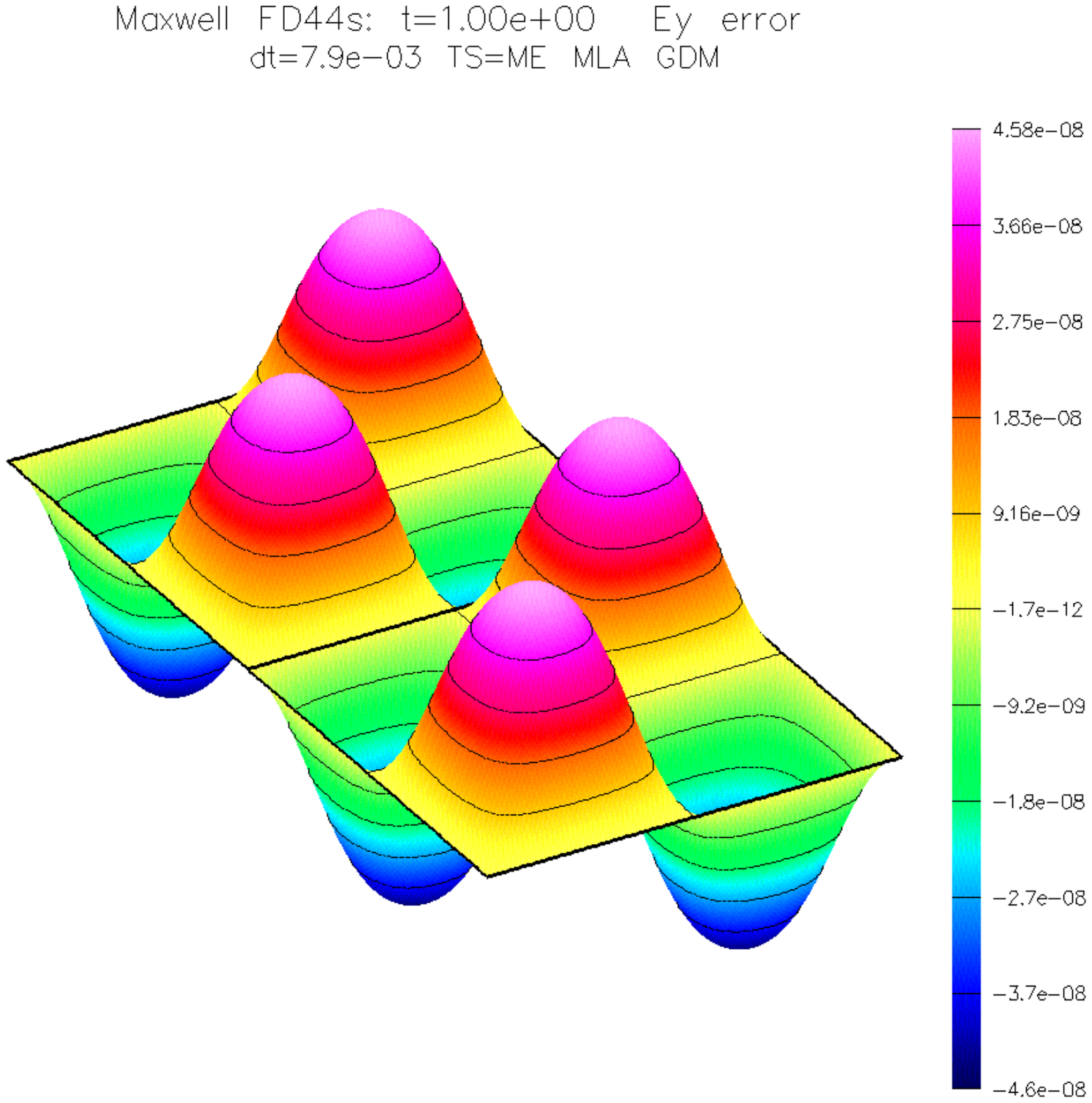}{$E_y$ error at $t=1.0$}{$p$}{$t=0.3$}{$-4.6e-8$}{$4.58e-8$};

  %


  %
  
%
\end{tikzpicture}
}
\end{center}
\caption{ Left: Grids for two rectangles;
          Middle: Solution $E_y$ at $t=1.0$;
          Right: Error of $E_y$ at $t=1.0$.
   }
  \label{fig:twoSquaresInterface}
\end{figure}
}

\subsection{Curved interface results}\label{sec:curvedInteraceResults}

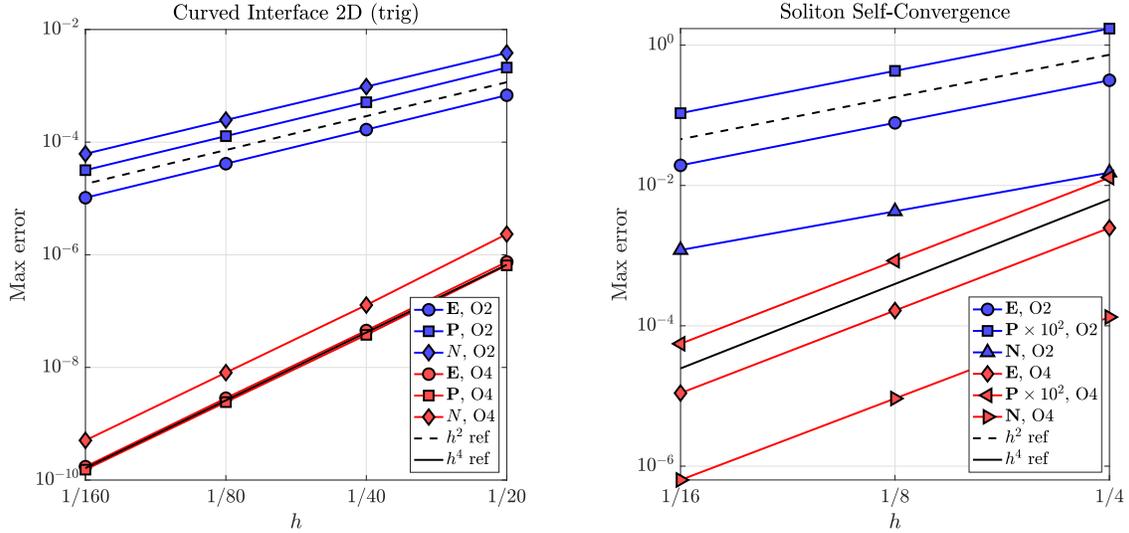
\begin{figure}[htb]
\begin{center}
\begin{tikzpicture}
   \useasboundingbox (0,.75) rectangle (15,7);  
  \figByHeight{  0}{0}{fig/curvedInterfaceMLA2d}{7cm}[0][0][0.][0];
  \figByHeight{  8}{0}{fig/solitonMLAConvergence}{7cm}[0][0][0.][0];
\end{tikzpicture}
\end{center}
\caption{Left: estimated max-norm errors for the curved interface problem at $t=1.0$ using a self-convergence grid refinement study.
  Right: Self-convergence estimated errors for the soliton solution at $t=100$.
  }
\label{fig:curvedInterfaceConvergence}
\end{figure}

In this section we verify the accuracy of the second-order and fourth-order accurate schemes for curved interfaces. The computational domain $[-1,1]\times[-1,1]$, the circular interface of radius $r_d=0.5$ and the overset grids are as shown in the left of Figure~\ref{fig:dielectricDisk}. We should note that the matched grids are used across the circular interfaces as is the square/cube case in Section~\ref{sec:planarInteraceResults}. For the numerical results, manufactured solutions in the form of \eqref{eqn:manufactured} are also employed, and their approximations and errors are shown in the middle and right of Figure~\ref{fig:dielectricDisk}. Grid refinement study is also performed in Figure~\ref{fig:curvedInterfaceConvergence}.
The MLA materials in the outer and inner domains are chosen to be \texttt{mlaMat2} and \texttt{mlaMat3} as
defined in~\ref{sec:materialDefinitions}.

{
\newcommand{\drawContour}[7]{%
\begin{scope}[#1]
\draw(0.0,0) node[anchor=south west,xshift=-4pt,yshift=+0pt] {\trimfiga{fig/#2}{\figWidtha}};
  \draw(.75,.5) node[draw,fill=white,anchor=west,xshift=2pt,yshift=1pt] {\scriptsize #3};
\begin{scope}[xshift=-.1cm,yshift=-3pt]
  \draw (\xcb,\ycb) node[anchor=south west,xshift=0.25cm,yshift=.5cm,rotate=-90] {\trimfigcb{fig/colourBarLines}{\cbWidth}{\cbHeight}};
  \draw (.8,0) node[anchor=north,xshift=+3pt,yshift=+2pt] {\scriptsize $#6$};
  \draw (4.8,0) node[anchor=north,xshift=+0pt,yshift=+2pt] {\scriptsize $#7$};
\end{scope}
\end{scope}
}
\newcommand{\drawContourC}[7]{%
\begin{scope}[#1]
\draw(0.0,0) node[anchor=south west,xshift=-4pt,yshift=+0pt] {\trimfigc{fig/#2}{\figWidth}};
  \draw(.75,.5) node[draw,fill=white,anchor=west,xshift=2pt,yshift=1pt] {\scriptsize #3};
\begin{scope}[xshift=-.1cm,yshift=-3pt]
  \draw (\xcb,\ycb) node[anchor=south west,xshift=0.25cm,yshift=.5cm,rotate=-90] {\trimfigcb{fig/colourBarLines}{\cbWidth}{\cbHeight}};
  \draw (.8,0) node[anchor=north,xshift=+3pt,yshift=+2pt] {\scriptsize $#6$};
  \draw (4.8,0) node[anchor=north,xshift=+0pt,yshift=+2pt] {\scriptsize $#7$};
\end{scope}
\end{scope}
}
\newcommand{\cbWidth}{.2cm}
\newcommand{\cbHeight}{4cm}
\newcommand{\xcb}{.5cm}
\newcommand{\ycb}{-.2cm}
\setlength{\ycbTop}{\ycb+\cbHeight}
\setlength{\ycbMid}{\ycb+\cbHeight*\real{.5}}
\newcommand{\trimfigcb}[3]{\includegraphics[width=#2, height=#3, clip, trim=17cm 2.35cm 1.65cm 2.35cm]{#1}}
%
%
\def\rad{1.27}
\newcommand{\plotDisk}{
\fill[fill=red!20,draw=red,line width=2pt] 
      plot[samples=100, domain=0.:360] ( {\rad*cos(\x)} , {\rad*sin(\x)} ) -- cycle ;
}
\newcommand{\figWidth}{5.cm}
\newcommand{\trimfig}[2]{\trimh{#1}{#2}{.02}{.02}{.02}{.02}}
\newcommand{\figWidtha}{5.75cm}
\newcommand{\trimfiga}[2]{\trimh{#1}{#2}{.04}{.12}{.07}{.09}}
\newcommand{\trimfigc}[2]{\trimh{#1}{#2}{.18}{.18}{.15}{.4}}
\begin{figure}[htb]
\begin{center}
\resizebox{15.25cm}{!}{
\begin{tikzpicture}[scale=1]
  \useasboundingbox (.4,0.2) rectangle (16.5,5.8);  
  \draw(0.25,0) node[anchor=south west,xshift=-4pt,yshift=-1pt] {\trimfig{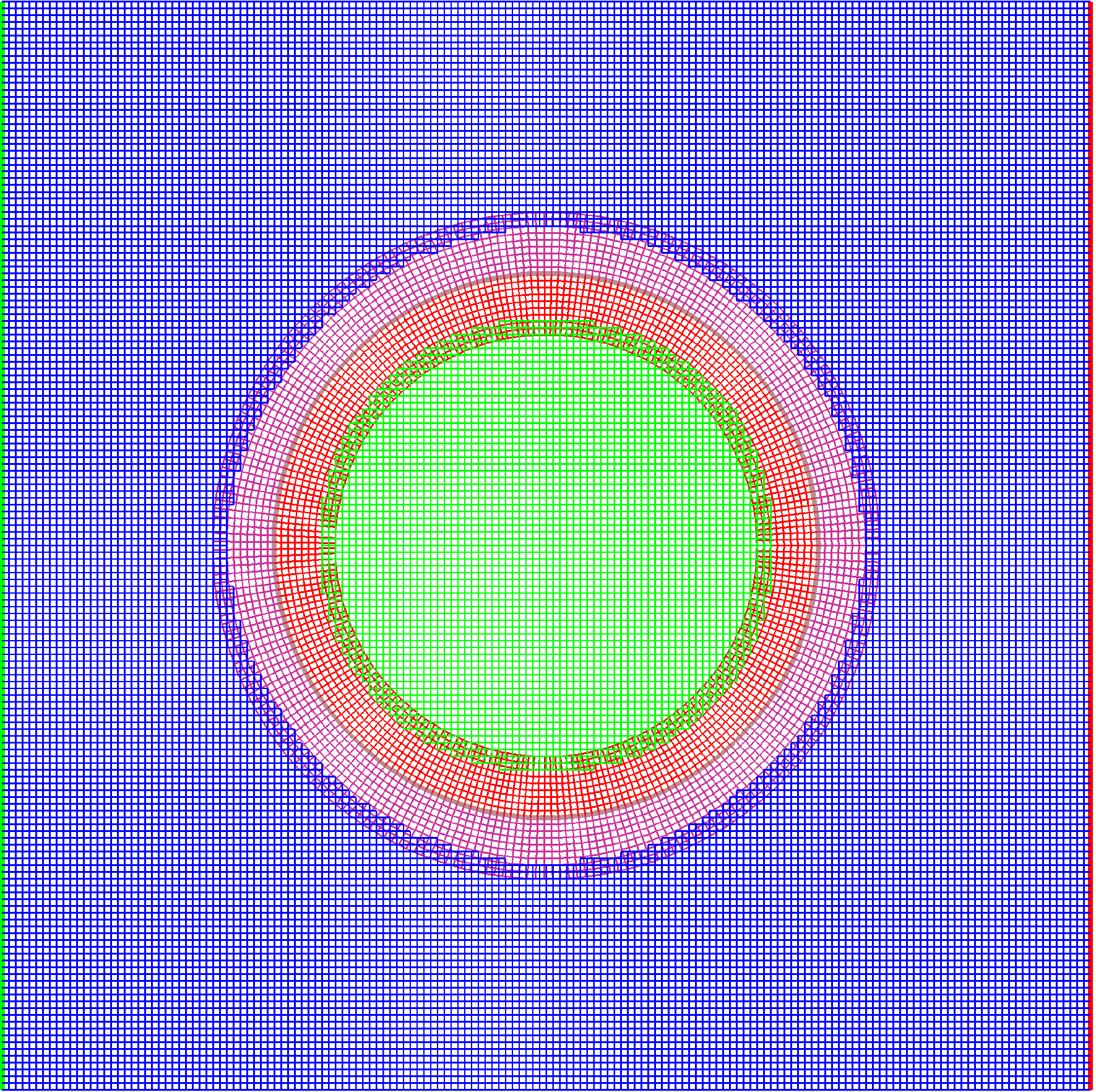}{\figWidth}};
   \drawContourC{xshift= 5.6cm,yshift=-.08cm}{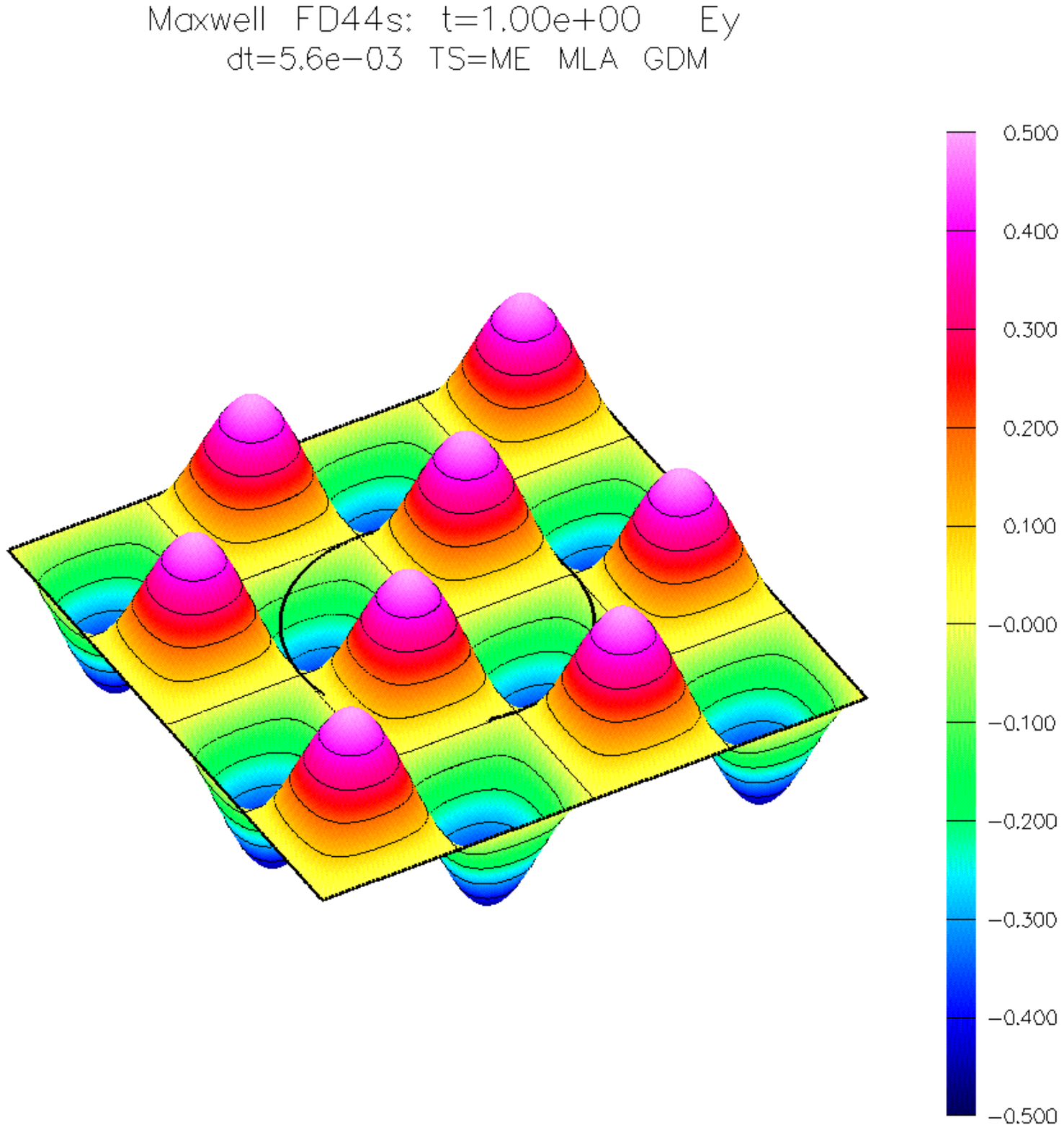}{$E_y$ $t=1.0$}{$\Ev$}{$t=0.3$}{$-0.5$}{$0.5$};
   \drawContourC{xshift=11.5cm,yshift=-.08cm}{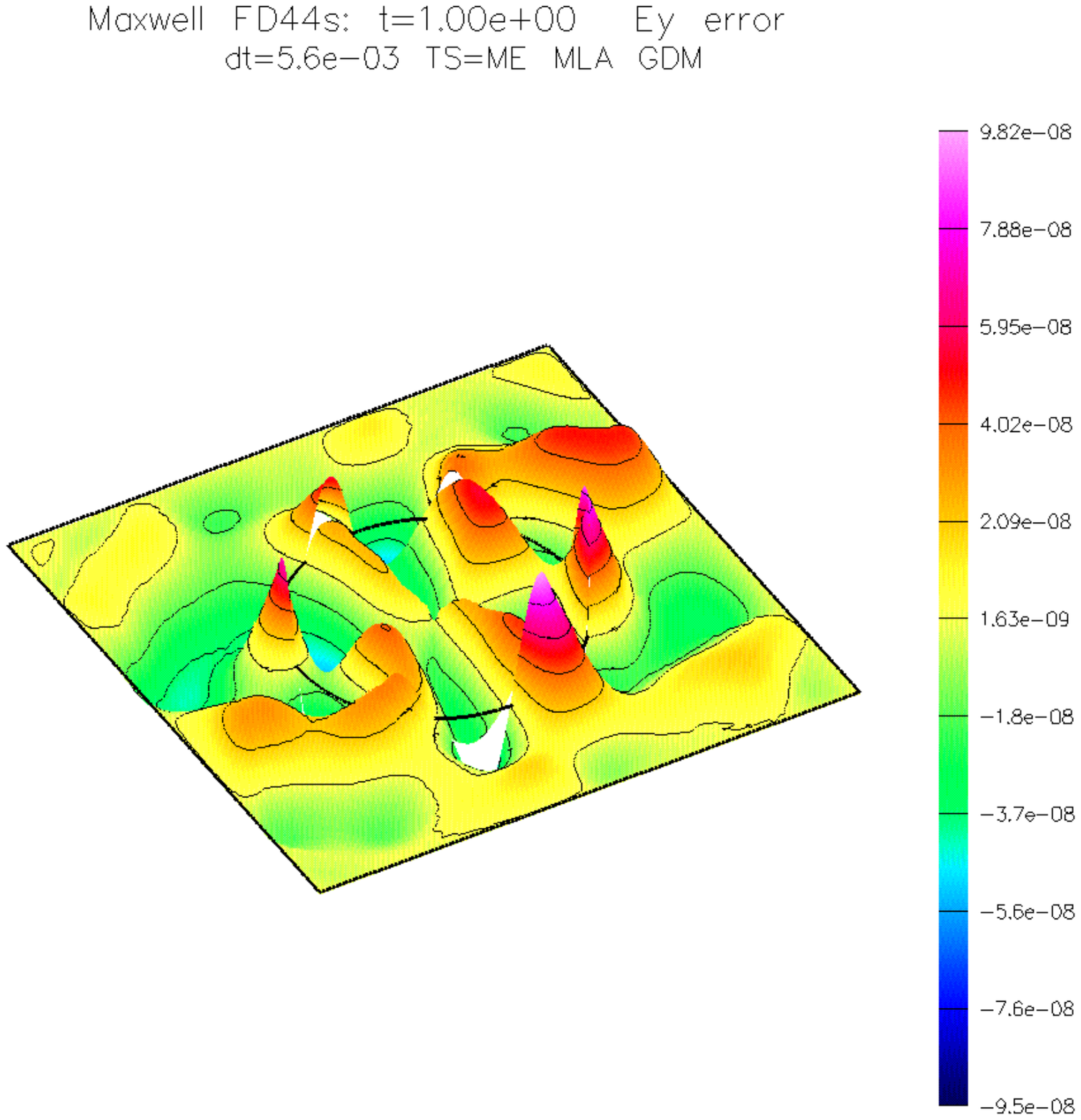}{$E_y$-error  $t=1.0$}{$p$}{$t=0.3$}{$-9.5e-8$}{$9.82e-8$};

  %
  \begin{scope}[xshift=.6cm,yshift=.05cm]
     \draw[|-|,very thick,xshift=-4pt] (0,0.1) node[anchor=east] {\footnotesize $-1$} -- (0,5.1) node[anchor=east] {\footnotesize $1$}; 
     \draw[|-|,very thick,yshift=-3pt] (0.08,0.0) node[anchor=north] {\footnotesize $-1$} -- (5.1,0) node[anchor=north,xshift=-3pt] {\footnotesize $1$}; 

     \begin{scope}[xshift=1.8pt,yshift=1.95pt]
       \draw[very thick] (2.5cm,2.5cm) circle(1.3cm);
       \draw[->,thick,xshift=2.5cm,yshift=2.5cm] (0,0) -- (.91,.91) node[anchor=east] {\scriptsize $r_d$};
      \end{scope}
  \end{scope}

  %
  
%
\end{tikzpicture}
}
\end{center}
\caption{ Left: A coarse grid representation of the composite grid $\Gc^{(2)}$ for the MLA disk.
 Middle: computed solution for $E_y$ at $t=1$ for a manufactured solution. Right: Errors in $E_y$.
   }
  \label{fig:dielectricDisk}
\end{figure}
}


\newcommand{\deltaHat}{\hat{\delta}}
\subsection{Soliton}\label{sec:soliton}




In this section, we present a soliton-like solution for the nonlinear system \eqref{eqn:nonlin_e}-\eqref{eqn:nonlin_n} in the following form:
\begin{subequations}\label{eqn:oneLevelSoliton}
\begin{align}
E_{tt} - c^2\Delta E = -\eta P_{tt},\label{eqn:2levelE}\\
P_{tt} + P = \deltaHat^2DE,\label{eqn:2levelP}\\
D_t = - EP_t,\label{eqn:2levelD}
\end{align}
\end{subequations}
where $D=N_0-N_1$ denotes the difference of carrier population density in the 2 atomic levels. The above system \eqref{eqn:oneLevelSoliton} can be cast into the Maxwell-MLA system \eqref{eqn:mbe} with 1 polarization and 1 level with parameters given in \ref{sec:mlaMat3}.



Multi-scale analysis in space and time of the 2-level system \eqref{eqn:2levelE}--\eqref{eqn:2levelD} would give the following asymptotic solutions
\begin{subequations}
\label{eq:solitionSolution}
\begin{align}
  E(x,t) &= 2\sqrt{\frac{\eta U}{1-U}}\sech(\deltaHat(x-x_0-Ut))\sin(x-t),\label{eqn:asymE1d}\\
  P(x,t) &= 2\deltaHat\tanh(\deltaHat(x-x_0-Ut))\sech( \deltaHat ( x-x_0-Ut))\cos(x-t),\label{eqn:asymP1d}\\
  D(x,t) &=1-2\sech^2(\deltaHat(x-x_0-Ut)),\label{eqn:asymD1d}
\end{align}
\end{subequations}
where $x_0$ is a free parameter for the center of the soliton solution. 

For the numerical simulations below, a thin rectangular domain $[0,1/2]\times[0,1000]$ is employed.
The boundary conditions on the left and right are set equal to the asymptotic solution~\eqref{eq:solitionSolution}, 
although this has negligible effect on the results since the solution is extremely small at these boundaries.
Periodic boundary conditions are imposed in the $y$-direction.
The initial conditions at $t=0$ are chosen to be $E(x,0)$, $P(x,0)$, $D(x,0)$, while Taylor's expansions are employed to obtain values of $E,P,D$ at time $t=-dt$, where low-order derivatives such as $E_t(x,0)$, and $P_t(x,0)$ are assumed known from the soliton solutions and high-order derivatives are computed recursively using the PDEs \eqref{eqn:2levelE}--\eqref{eqn:2levelD}, with parameters $x_0=0$, $U=1/2$, $\eta=1$, $c=1$, and $\deltaHat=0.1$. 

\plotTwoFigsByWidth{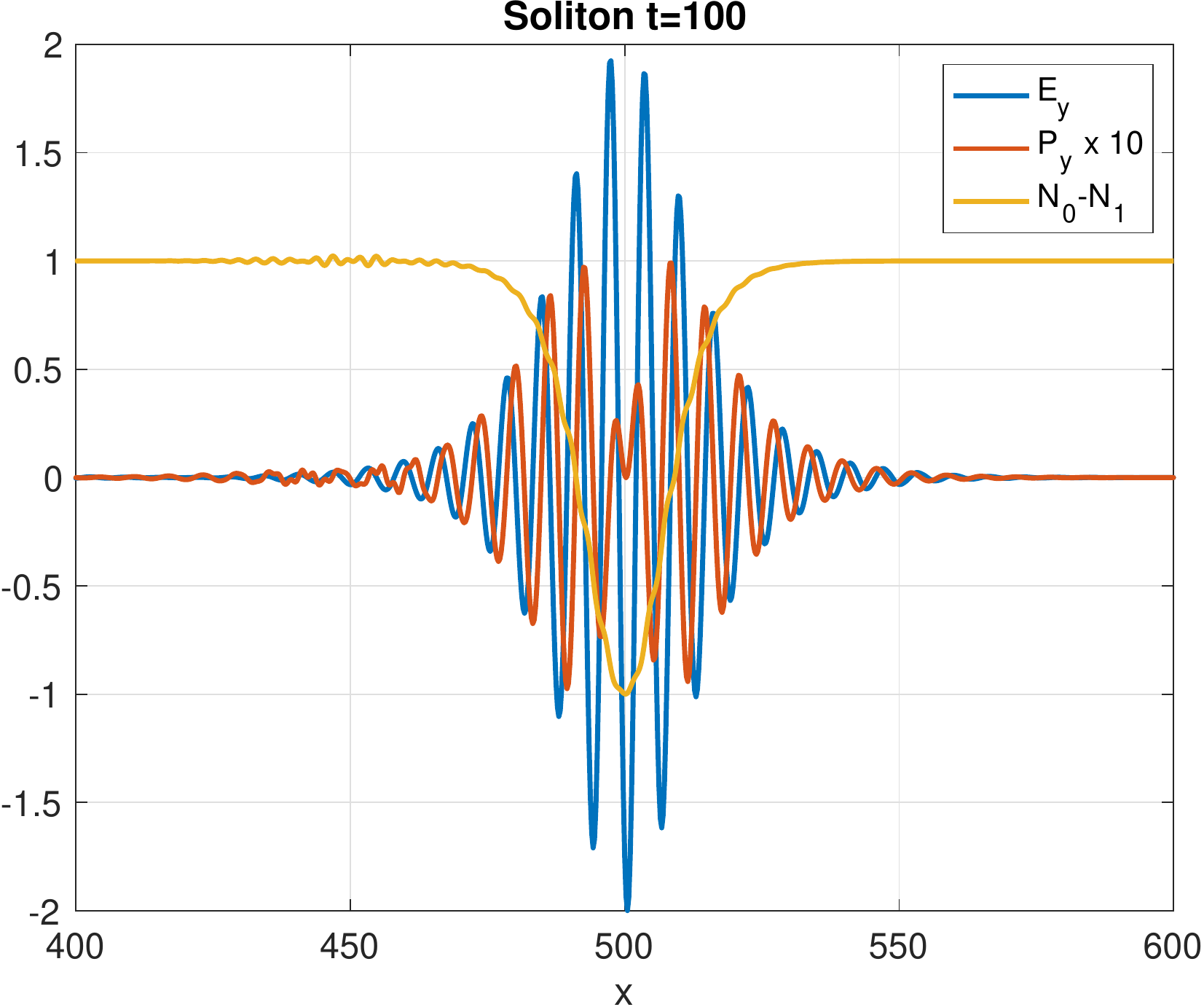}{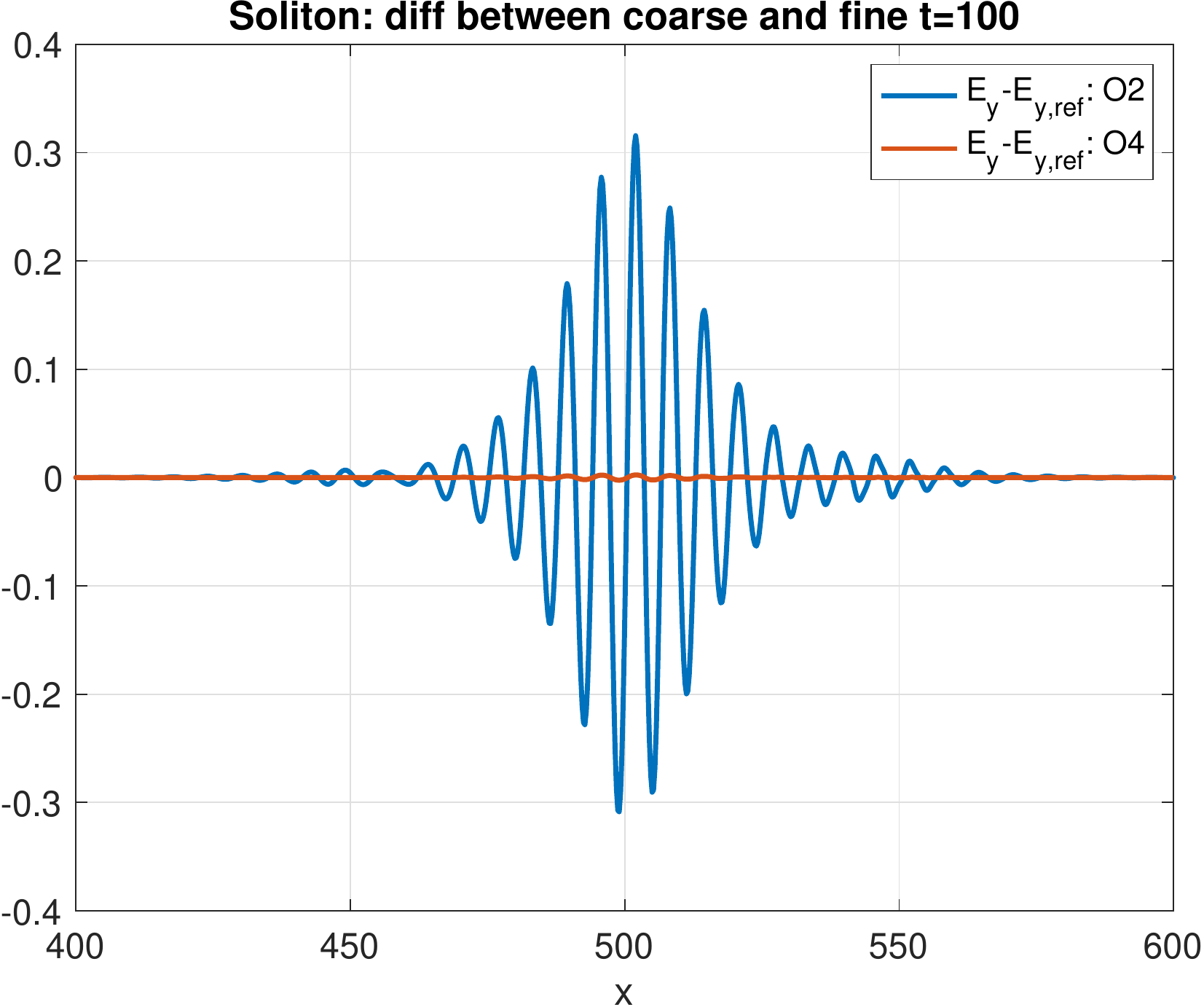}{Left: line plot along $y=0$ of the O4 soliton solutions ($h=1/4$). Right: difference between O2/O4 with a fine grid O4 solution $E_{y,ref}$ ($h=1/16$).  }{fig:soliton}{8cm}

In the left figure of Figure~\ref{fig:soliton}, extracted line plots of the approximated soliton solutions at $t=100$ along $y=0$ are presented, while the right figure illustrates the accuracy between second-order accurate (O2) simulations with fourth-order accurate (O4) simulations (both on coarse meshes with $h=1/4$) by comparing their differences with O4 simulations on a finer mesh with grid spacing $h=1/16$, which is denoted as $E_{y,ref}$ in the legend.  
Figure~\ref{fig:curvedInterfaceConvergence} (right)
compares the self-convergence study for the second- and fourth-order accurate simulations.

\subsection{Scattering from a layered disk} \label{sec:layeredDisk}



In this section, we provide an example of multiple types of material interfaces, i.e., interfaces between both linear/nonlinear and linear/linear materials. The computational domain is $[-1.75,1.75]\times[-1.5,1.5]$ with a disk of radius $r=0.4$ centered at origin and two layers of width 0.1 outside the disk. The background rectangle is assumed to be vacuum with normalized permittivity $\epsilon=1$, 
while the center disk is made of nonlinear 4-level active material with 2 polarization vectors with $\epsilon=2$ as depicted in Fig~\ref{fig:4-level}
(see~\ref{sec:mlaMat4levels} for the material parameters). The first layer is of linear material modeled by generalized dispersive model (GDM) with 1 polarization vector and $\epsilon=4$ as in \cite{angel2019high,banks2020high} and interface treatments of GDM materials can be found therein. The second layer is assumed to be a linear isotropic material with $\epsilon=3$. The permeability in all materials are set to $\mu=1$.



{
\newcommand{\drawContour}[7]{%
\begin{scope}[#1]
\draw(0.0,0) node[anchor=south west,xshift=-4pt,yshift=+0pt] {\trimfiga{fig/#2}{\figWidtha}};
\end{scope}
}
\newcommand{\drawContourC}[7]{%
\begin{scope}[#1]
\draw(0.0,0) node[anchor=south west,xshift=-4pt,yshift=+0pt] {\trimfigc{fig/#2}{\figWidth}};
\end{scope}
}
\newcommand{\cbWidth}{.2cm}
\newcommand{\cbHeight}{4cm}
\newcommand{\xcb}{.5cm}
\newcommand{\ycb}{-.2cm}
\setlength{\ycbTop}{\ycb+\cbHeight}
\setlength{\ycbMid}{\ycb+\cbHeight*\real{.5}}
\newcommand{\trimfigcb}[3]{\includegraphics[width=#2, height=#3, clip, trim=17cm 2.35cm 1.65cm 2.35cm]{#1}}
%
%
\newcommand{\figWidth}{5.cm}
\newcommand{\trimfig}[2]{\trimh{#1}{#2}{.02}{.02}{.02}{.02}}
\newcommand{\figWidtha}{5.75cm}
\newcommand{\trimfiga}[2]{\trimh{#1}{#2}{.04}{.12}{.07}{.09}}
\newcommand{\trimfigc}[2]{\trimh{#1}{#2}{.1}{.1}{.08}{.3}}
\begin{figure}[htb]
\begin{center}
\resizebox{17cm}{!}{
\begin{tikzpicture}[scale=1]
  \useasboundingbox (0.4,.65) rectangle (16.5,5.25);  
   \draw(0.,0) node[anchor=south west,xshift=-4pt,yshift=-1pt] {\trimfig{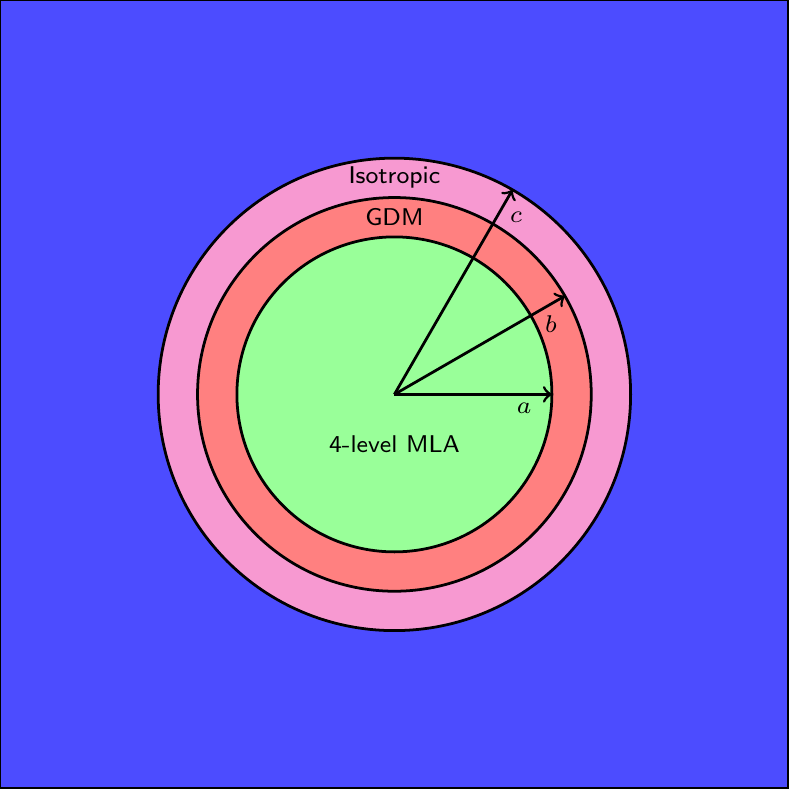}{\figWidth}};
   \drawContourC{xshift=5.3cm,yshift=0.cm}{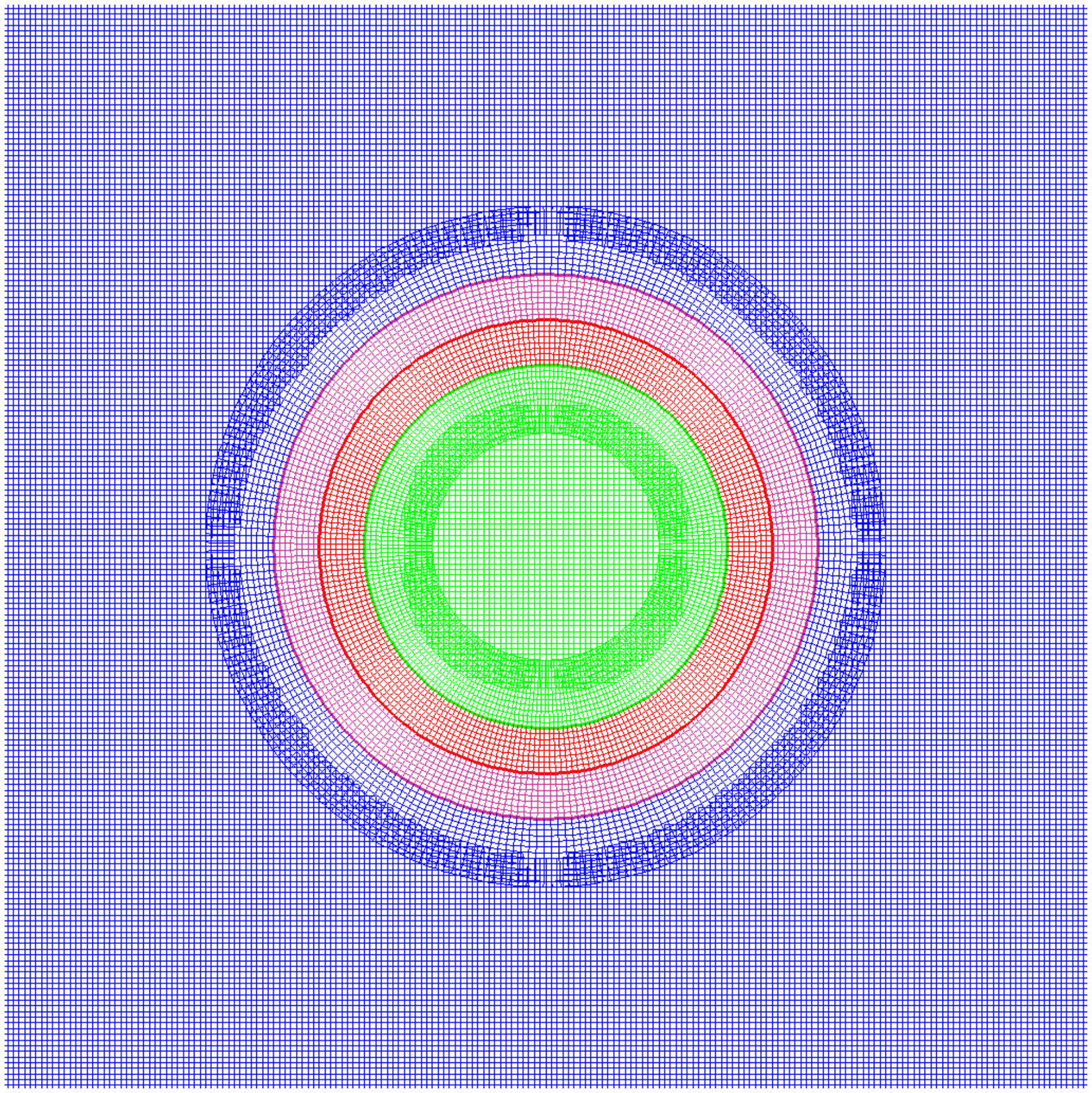}{}{$p$}{$t=0.3$}{$0$}{$1.04325$};
   \drawContourC{xshift=10.6cm,yshift=0.cm}{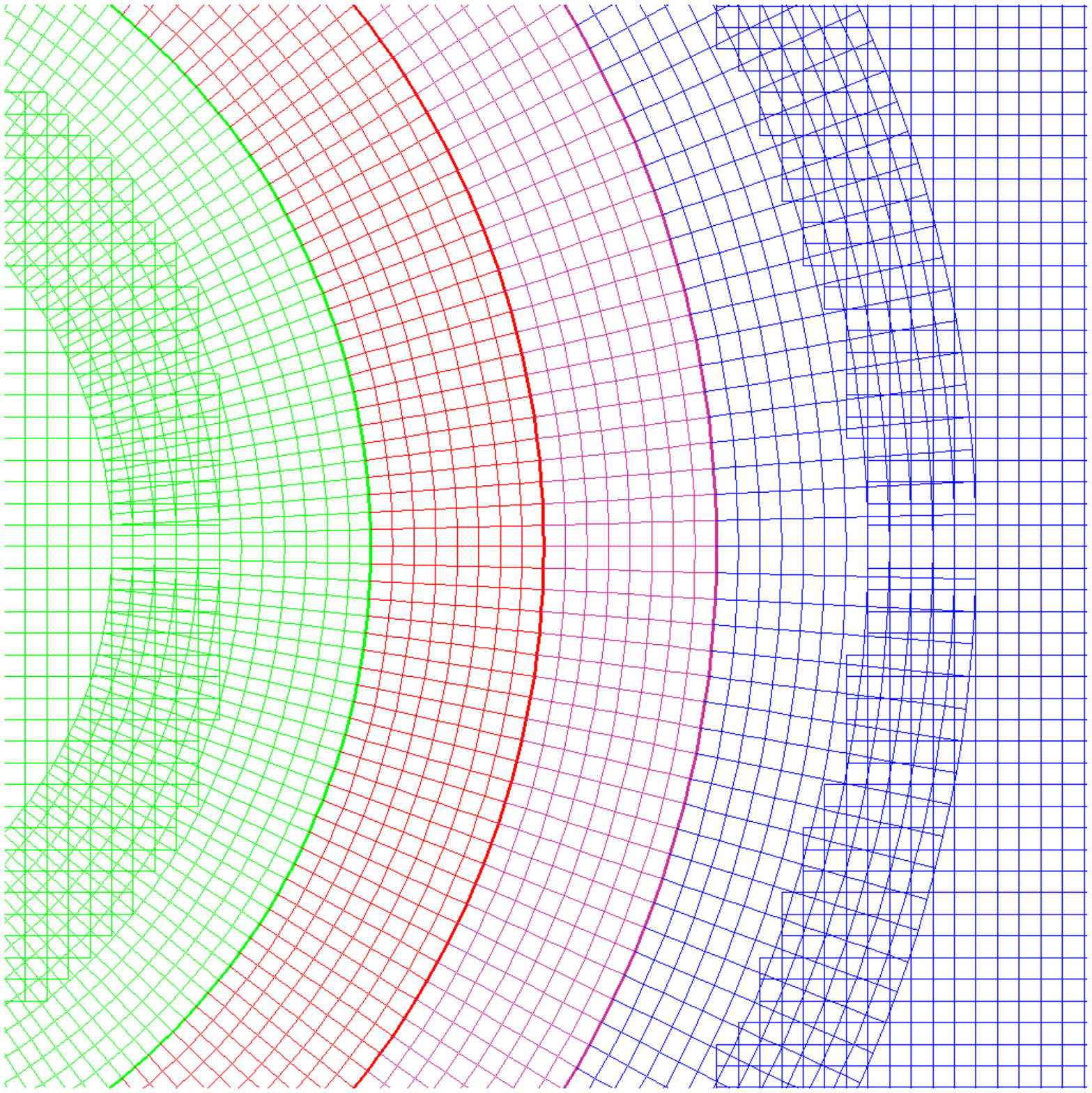}{}{$p$}{$t=0.3$}{$0$}{$1.04325$};

  %


  %
  
%
\end{tikzpicture}
}
\end{center}
\caption{ Left: material configurations; Middle: overset grids for material bodies; Right: zoomed-in view of the overset grids.
   }
  \label{fig:oneDiskTwoLayersGrids}
\end{figure}
}

{
\newcommand{\drawContour}[7]{%
\begin{scope}[#1]
\draw(0.0,0) node[anchor=south west,xshift=-4pt,yshift=+0pt] {\trimfiga{fig/#2}{\figWidtha}};
  \draw(.75,.5) node[draw,fill=white,anchor=west,xshift=2pt,yshift=1pt] {\scriptsize #3};
\begin{scope}[xshift=-.1cm,yshift=-3pt]
  \draw (\xcb,\ycb) node[anchor=south west,xshift=0.25cm,yshift=.5cm,rotate=-90] {\trimfigcb{fig/colourBarLines}{\cbWidth}{\cbHeight}};
  \draw (.8,0) node[anchor=north,xshift=+3pt,yshift=+2pt] {\scriptsize $#6$};
  \draw (4.8,0) node[anchor=north,xshift=+0pt,yshift=+2pt] {\scriptsize $#7$};
\end{scope}
\end{scope}
}
\newcommand{\drawContourC}[7]{%
\begin{scope}[#1]
\draw(0.0,0) node[anchor=south west,xshift=-4pt,yshift=+0pt] {\trimfigc{fig/#2}{\figWidth}};
  \draw(.75,.75) node[draw,fill=white,anchor=west,xshift=2pt,yshift=1pt] {\scriptsize #3};
\begin{scope}[xshift=-.1cm,yshift=-3pt]
  \draw (\xcb,\ycb) node[anchor=south west,xshift=0.25cm,yshift=.5cm,rotate=-90] {\trimfigcb{fig/colourBarLines}{\cbWidth}{\cbHeight}};
  \draw (.75,0) node[anchor=north,xshift=+3pt,yshift=+2pt] {\scriptsize $#6$};
  \draw (4.8,0) node[anchor=north,xshift=+0pt,yshift=+2pt] {\scriptsize $#7$};
\end{scope}
\end{scope}
}
\newcommand{\cbWidth}{.2cm}
\newcommand{\cbHeight}{4cm}
\newcommand{\xcb}{.5cm}
\newcommand{\ycb}{-.2cm}
\setlength{\ycbTop}{\ycb+\cbHeight}
\setlength{\ycbMid}{\ycb+\cbHeight*\real{.5}}
\newcommand{\trimfigcb}[3]{\includegraphics[width=#2, height=#3, clip, trim=17cm 2.35cm 1.65cm 2.35cm]{#1}}
%
%
\newcommand{\figWidth}{5.cm}
\newcommand{\trimfig}[2]{\trimh{#1}{#2}{.02}{.02}{.02}{.02}}
\newcommand{\figWidtha}{5.75cm}
\newcommand{\trimfiga}[2]{\trimh{#1}{#2}{.04}{.12}{.07}{.09}}
\newcommand{\trimfigc}[2]{\trimh{#1}{#2}{.18}{.18}{.15}{.4}}
\begin{figure}[htb]
\begin{center}
\resizebox{15.25cm}{!}{
\begin{tikzpicture}[scale=1]
  \useasboundingbox (.4,0.2) rectangle (16.5,5.8);  
   \drawContourC{xshift=0.cm,yshift=0.cm}{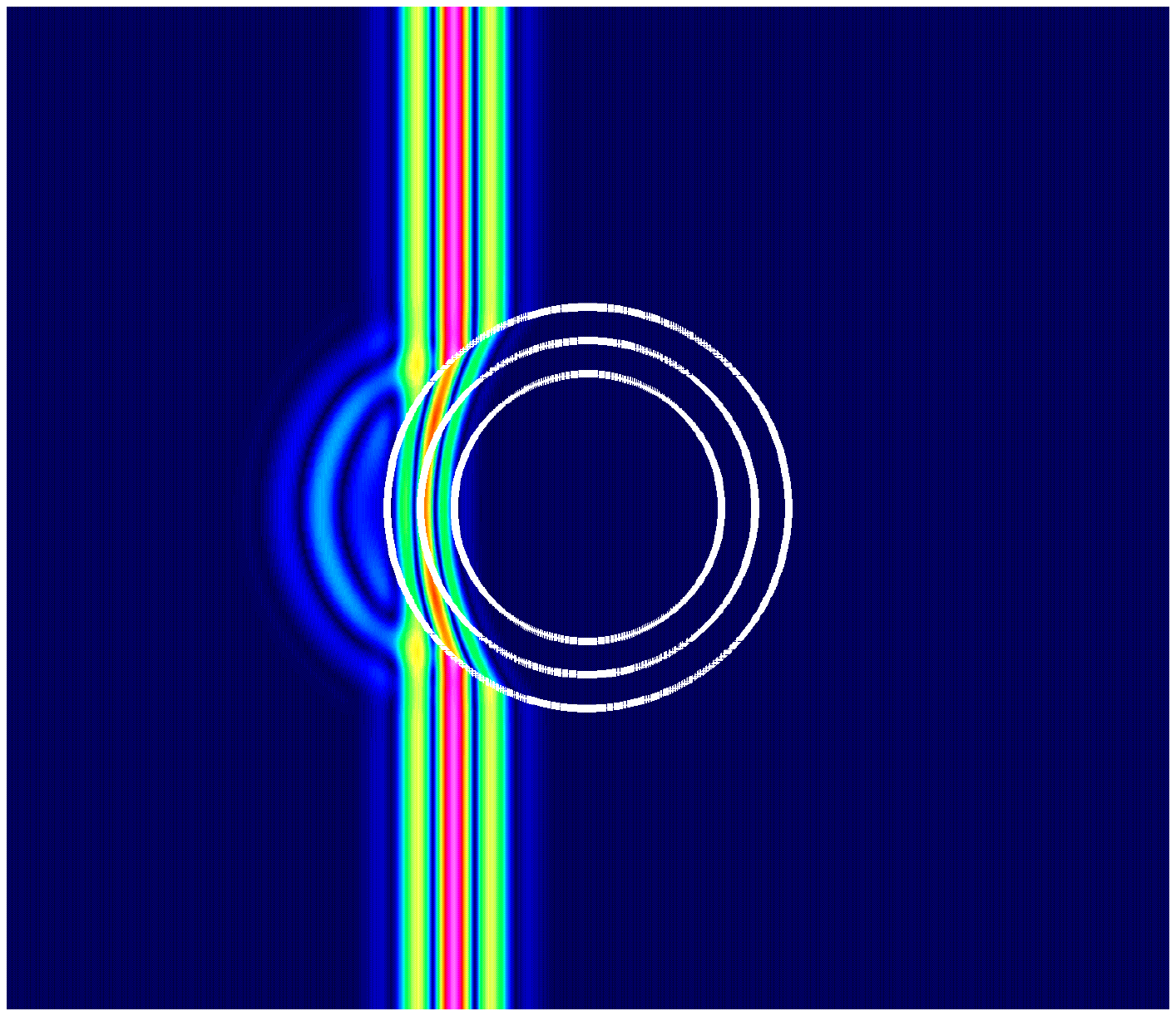}{$||\bm{E}||$ at $t=2.6$}{$p$}{$t=0.3$}{$0$}{$1.038$};
   \drawContourC{xshift=5.6cm,yshift=0.cm}{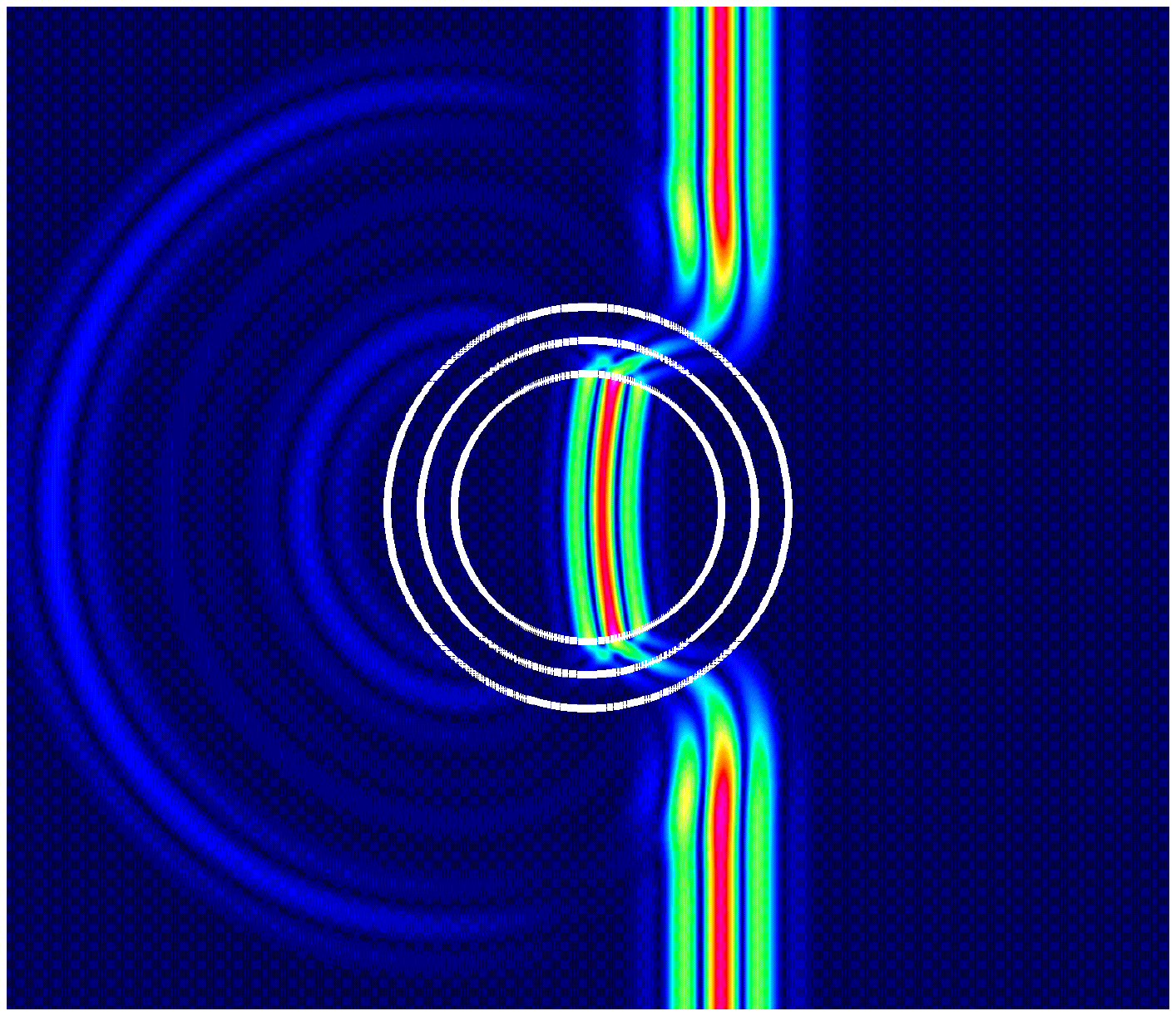}{$||\bm{E}||$ at $t=3.4$}{$p$}{$t=0.3$}{$0$}{$1.216$};
   \drawContourC{xshift=11.5cm,yshift=0.cm}{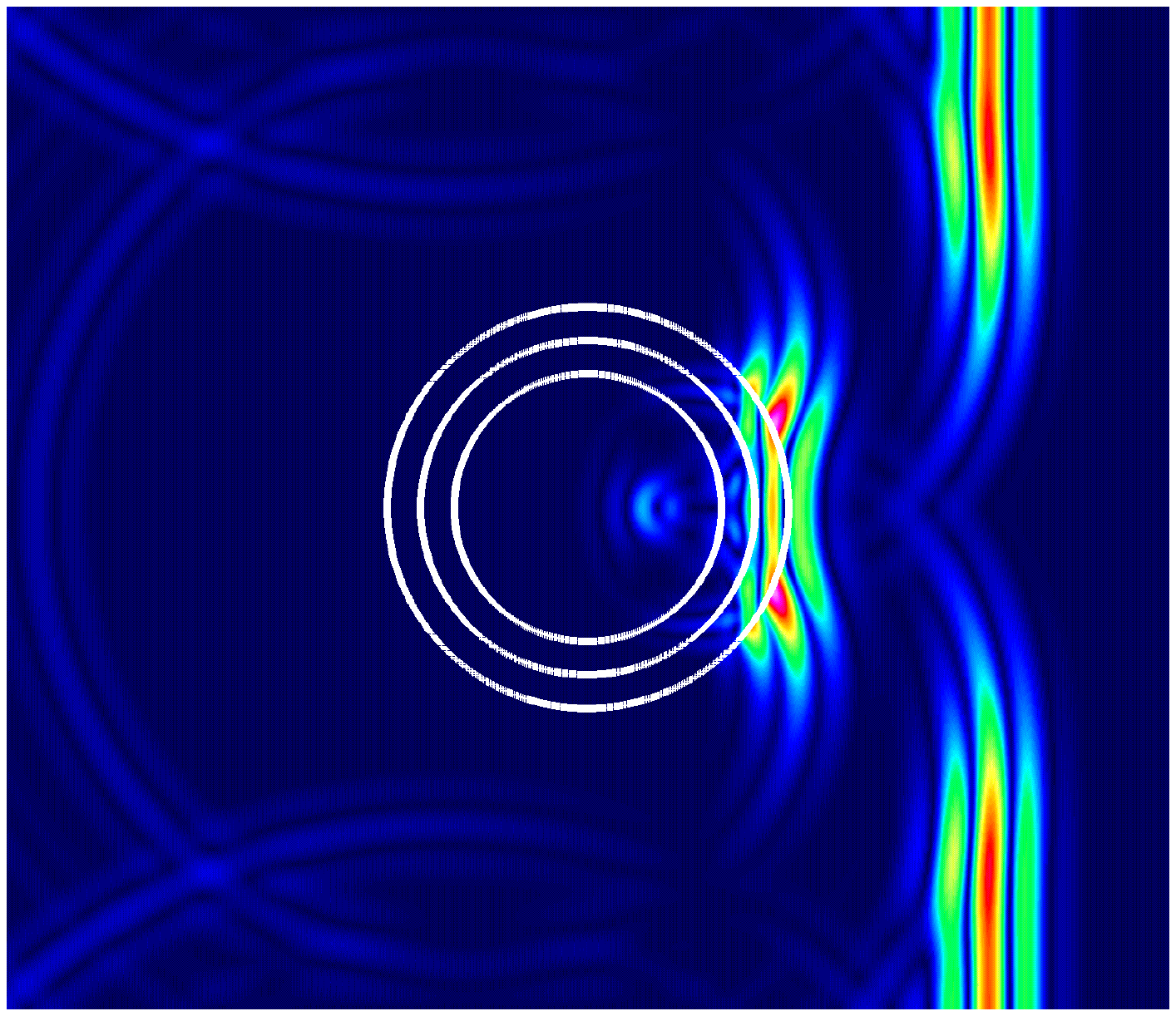}{$||\bm{E}||$ at $t=4.2$}{$p$}{$t=0.3$}{$0$}{$1.043$};

  %


  %
  
%
\end{tikzpicture}
}
\end{center}
\caption{ Snapshots of $\| \Ev\|$ at 3 different times.
   }
  \label{fig:oneDiskTwoLayersEfiledNorm}
\end{figure}
}


\begin{figure}[htb]
\begin{center}
\begin{tikzpicture}
   \useasboundingbox (0,.75) rectangle (15,7);  
  \figByHeight{  0}{0}{fig/oneDiskTwoLayersMLAConvergence}{7cm}[0][0][0.][0]
  \figByHeight{8.25}{0}{fig/oneDiskTwoLayersMLARuntime}{7cm}[0][0][0.][0]
\end{tikzpicture}
\end{center}
\caption{Left: Self-convergence for one disk with two layers at $t=3.8$. 
Right: Comparison of runtime between order 4 (O4) and order 2 (O2) simulations a final time $t=5$ with approximate speedup in purple.
}
\label{fig:layeredDiskResults}
\end{figure}
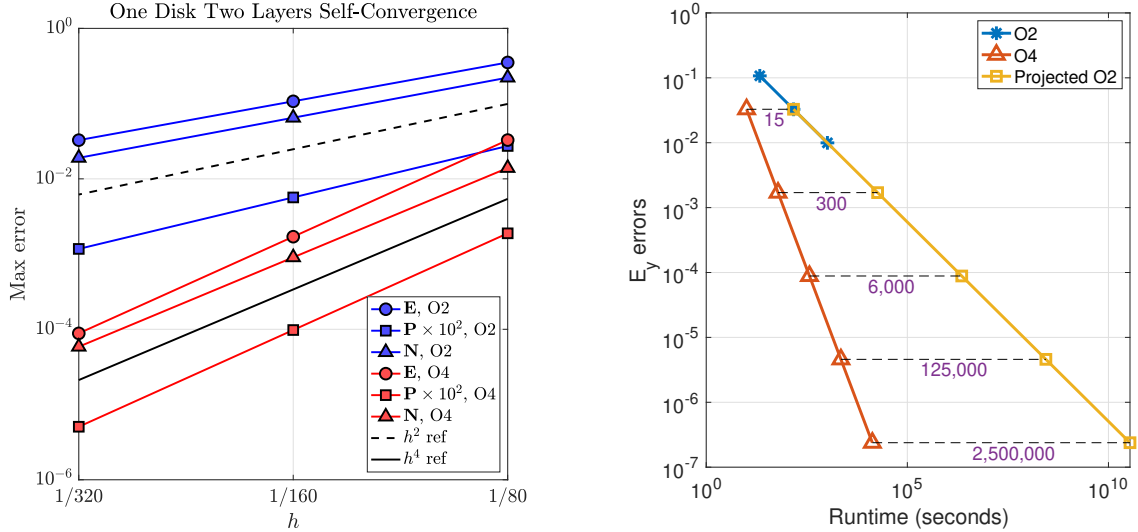

For the numerical simulations below, the initial conditions are from a modulated Gaussian plane wave with nonzero $E_y$ component, i.e.,
\begin{align}\label{eqn:gpw}
E_y = \exp(-50(x+3-t)^2)\cos(4\pi(x+3-t)),
\end{align}
with center $x=-3$. Nonlocal radiation boundary conditions\footnote{The non-local radiation boundary
conditions are based on the work in~\cite{AlpertGreengardHagstrom2002}.}
 are imposed on the left and right of the domain while periodic boundary conditions are imposed on the top and bottom of the domain. The overlapping grids (zoomed-in view) can be found in Figure~\ref{fig:oneDiskTwoLayersGrids}.



%




In Figure~\ref{fig:oneDiskTwoLayersEfiledNorm}, we show three snapshots of the norm of the electric fields as a result of interacting with the nonlinear 4 level active material, which also demonstrates a focusing effect from the circular geometry. 
Figure~\ref{fig:layeredDiskResults} (left) gives the results of self-convergence study from both the second-order and fourth-order accurate simulations at $t=3.8$. Furthermore, a comparison of runtime between order 2 (O2) and order 4 (O4) simulations with the same accuracy till final time $t=5$ are performed using the serial codes are shown in Figure~\ref{fig:layeredDiskResults} (right), where the purple digits indicate the approximate wall-clock runtime speedup between O2 and O4.
The projected O2 runtimes are predicted using the fact that the computational cost on a 2$\times$-refined mesh is 8$\times$ slower.


\subsection{An array of active material ellipsoids}\label{sec:ellipsoidArray}


As a final example we consider the scattering of a Gaussian plane wave from a collection
of 36 solid ellipsoid meta-atoms, each containing an active material.
The overset grid for the geometry in shown in Figure~\ref{fig:ellipsoidArrayFieldsAndGrid}.
The ellipsoids are enclosed in a rectangular box $\Bc=[x_a,x_b]\times[y_a,y_b]\times[z_a,z_b]$. 
The ellipsoids have different shapes and orientations. To avoid polar-type singularities in the grid mappings,
the surface of the ellipsoid is covered with three patches. 
The interior of each ellipsoid contains an active material defined by the active material \texttt{mlaMat4levels},
given in~\ref{sec:mlaMat4levels}.
The ellipsoids are surrounded by a vacuum region. A Gaussian plane wave enters the domain from the left 
at the face $x=x_a$. Radiation boundary conditions are used the faces $x_a$ and $x_b$ while the 
solution is periodic in $y$ and $z$. 
A $y$-polarized Gaussian plane wave travels in the $x$-direction and enters the domain through the left face at $x=x_a$.

Figure~\ref{fig:ellipsoidArrayFieldsAndGrid} shows the contours of the computed solution on selected contour cutting planes.
The magnitude of the electric field $\| \Ev\|$ is shown along with the $P_y$ component of total polarization and
the population density $N_3$. 
The incident wave is seen to excite the meta-atoms as it passes through.








{
\newcommand{\drawContour}[7]{%
\begin{scope}[#1]
\draw(0.0,0) node[anchor=south west,xshift=-4pt,yshift=+0pt] {\trimfiga{fig/#2}{\figWidtha}};
  \draw(.5,.5) node[draw,fill=white,anchor=west,xshift=2pt,yshift=1pt] {\scriptsize #3};
\begin{scope}[xshift=.2cm,yshift=-3pt]
  \draw (\xcb,\ycb) node[anchor=south west,xshift=0.25cm,yshift=.5cm,rotate=-90] {\trimfigcb{fig/colourBarLines}{\cbWidth}{\cbHeight}};
  \draw (.8,0) node[anchor=north,xshift=+3pt,yshift=+2pt] {\scriptsize $#6$};
  \draw (4.8,0) node[anchor=north,xshift=+0pt,yshift=+2pt] {\scriptsize $#7$};
\end{scope}
\end{scope}
}
\newcommand{\cbWidth}{.2cm}
\newcommand{\cbHeight}{4cm}
\newcommand{\xcb}{.5cm}
\newcommand{\ycb}{-.2cm}
\setlength{\ycbTop}{\ycb+\cbHeight}
\setlength{\ycbMid}{\ycb+\cbHeight*\real{.5}}
\newcommand{\trimfigcb}[3]{\includegraphics[width=#2, height=#3, clip, trim=17cm 2.35cm 1.65cm 2.35cm]{#1}}
\newcommand{\figWidtha}{5.75cm}
\newcommand{\trimfiga}[2]{\trimwb{#1}{#2}{.0}{.0}{.0}{.0}}
\begin{figure}[htb]
\begin{center}
\begin{tikzpicture}
   \useasboundingbox (0,.25) rectangle (11.5,12);  

   \figByWidth{   0}{6.00}{fig/36SolidEllipsoidsGrid}{6cm}[0][0][0.][0]

   \drawContour{xshift=6.cm,yshift=6.25cm}{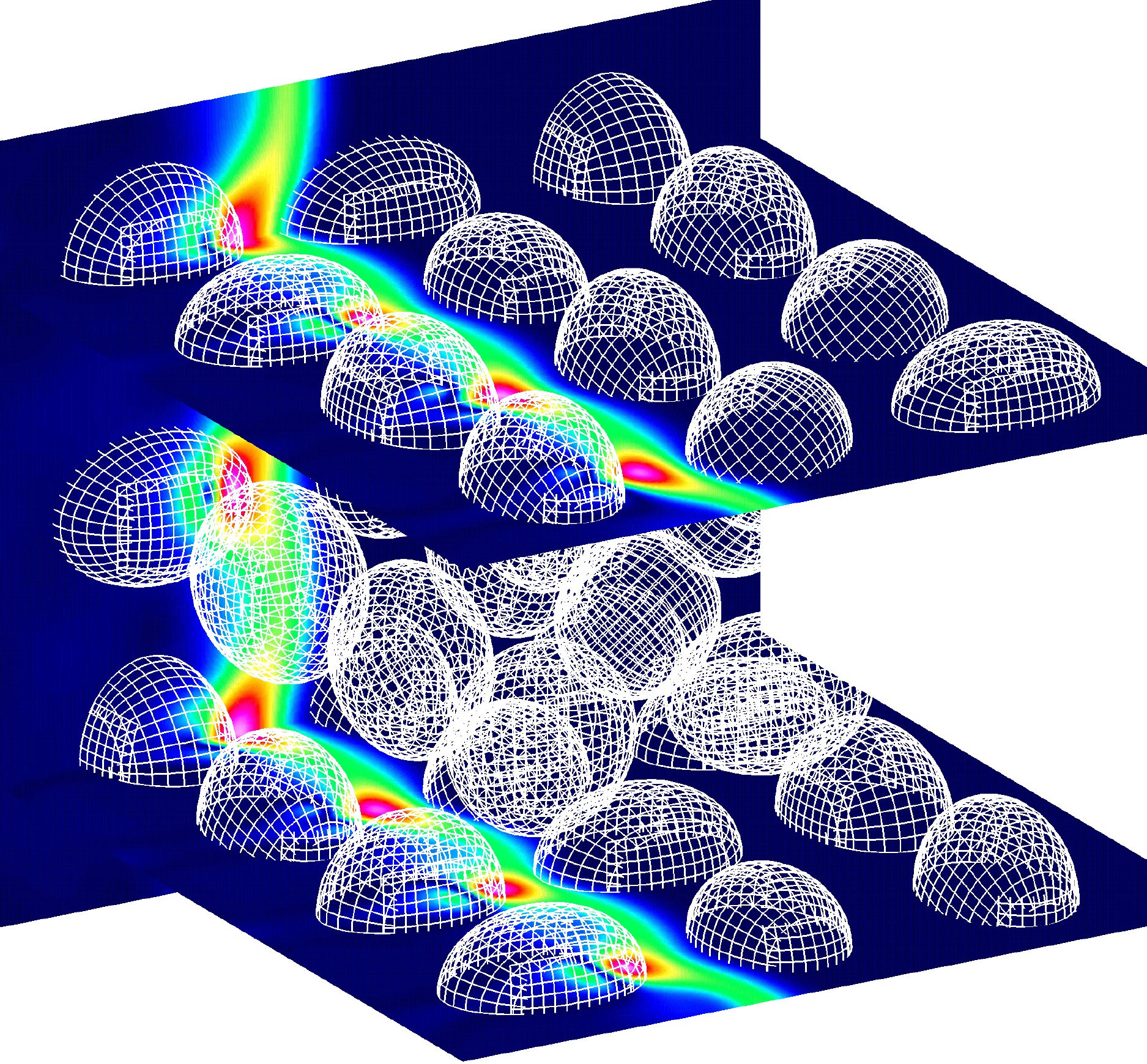}{$||\Ev||$ at $t=2.5$}{$p$}{$t=0.3$}{$0$}{$1.82$};
   \drawContour{xshift=0.cm,yshift=0.00cm}{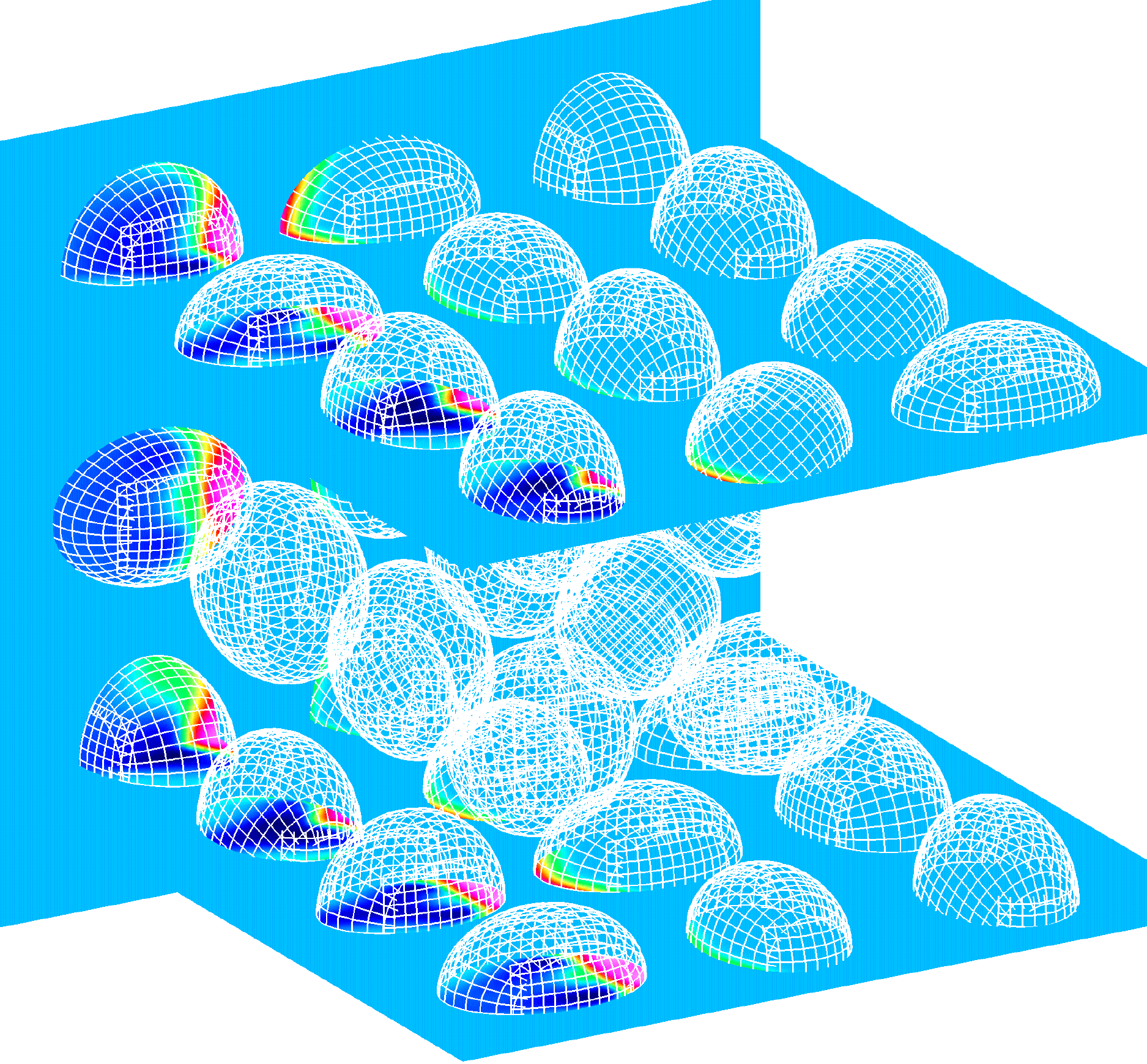}{$P_y$ at $t=2.5$}{$p$}{$t=0.3$}{$-1.9e-4$}{$7.1e-4$};
   \drawContour{xshift=6.cm,yshift=0.00cm}{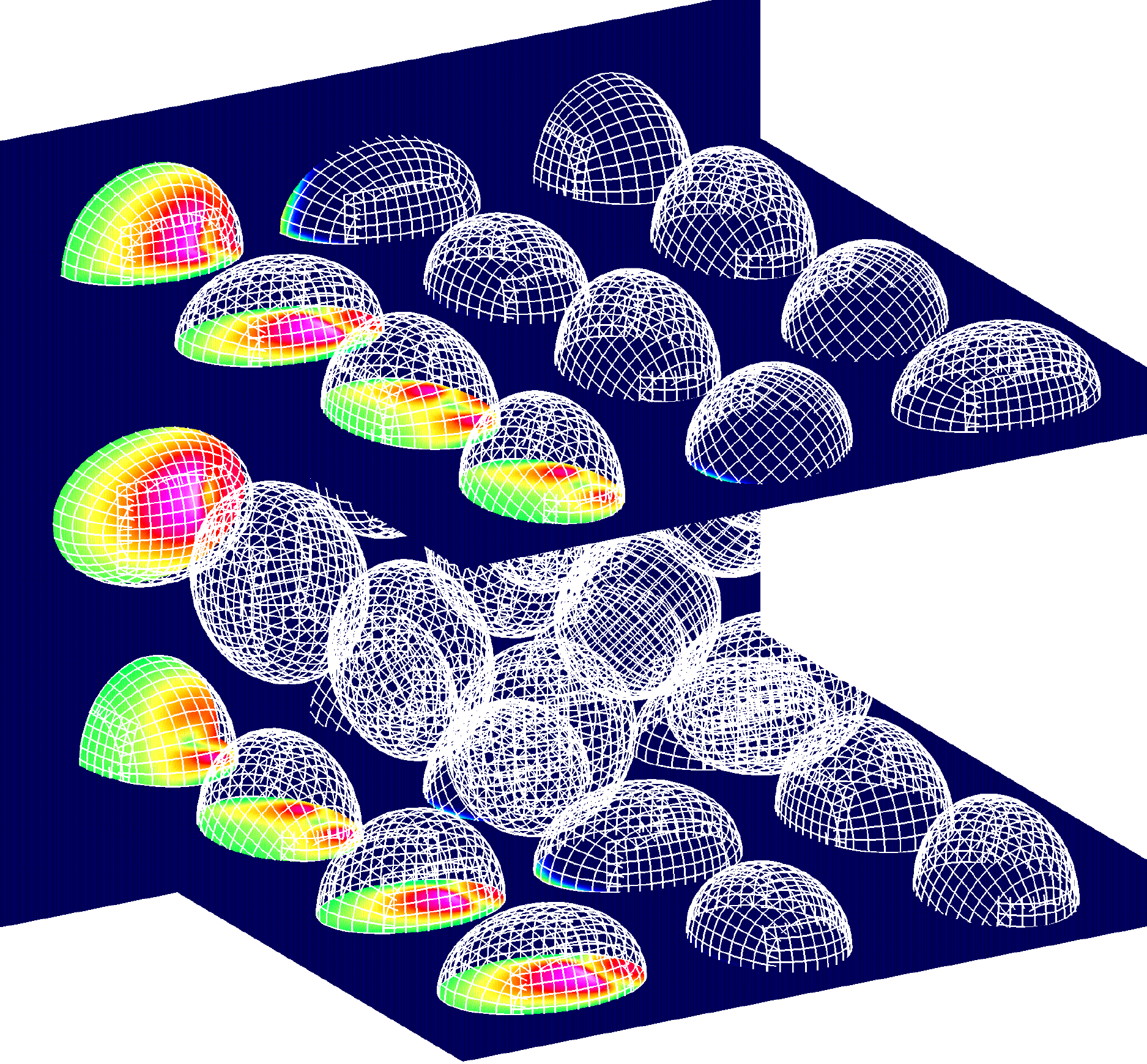}{$N_3$ at $t=2.5$}{$p$}{$t=0.3$}{$0$}{$.56$};


   %
   \draw( .00,6.5) node[draw,fill=white,anchor=west,xshift=2pt,yshift=1pt] {\scriptsize Overset grid.}; 
\end{tikzpicture}
\end{center}
\caption{Gaussian plane wave hitting thirty-six solid-ellipsoid meta-atoms. Each ellipsoid contains a four-level MLA active material.
   Left: overset grid showing the grids on the ellipsoid surfaces.
   Right and bottom: contours of the electric field-norm $\|\Ev\|$, the $y$-component of the total polarization $P_y$,
     and the level-three population density, $N_3$, are shown.
     }
\label{fig:ellipsoidArrayFieldsAndGrid}
\end{figure}
}







\section{Conclusions}\label{sec:conclusion}

\
  A high-order accurate finite-difference time-domain scheme for solving Maxwell's equations coupled to multi-level carrier kinetics models was developed.
The Maxwell-MLA method uses an efficient single-step three-level modified-equation for time-stepping approach Maxwell's equations in the
form of a second-order vector wave equation. Nonlinear effects for active materials is treated with a fairly general class of multi-level atomic models
involving ODEs for any number of polarization vectors and population densities. Complex geometries with curved boundary and interfaces
are accurately treated using conforming and overset grids.
One key property of the scheme is that through a hierarchical modified equation (HIME) approach no nonlinear solves
are required to time-step the equations at high-order accuracy. 
Another key property is that through the use of a hierarchical approach that couples low-order accurate and
high-order accurate approximations, the update of the ghost values at the interface are local with no nonlinear
solves and no coupling with adjacent ghost values.
Stability on overset grids was maintained using a novel high-order upwind scheme that applies to wave equations in second-order form.
The initial-boundary value problem for the Maxwell-MLA equations was shown to be well posed. 
An $L_2$-energy estimate was derived for a restricted class of commonly used MLA models 
which showed long-time existence. Numerical results in two and three dimensions were presented that verified the accuracy and stability of
the new Maxwell-MLA schemes. Results were shown for problems with interfaces between two active materials and also for 
active materials adjacent to linear dispersive materials.
Verification was performed using manufactured solutions and an asymptotic soliton solution. For problems
without exact solutions, a self-convergence grid refinement procedure was used to estimate the errors and convergence rates.
%
%
   Some possible future steps include extending the current numerical schemes to handle junctions of three or materials,
   adding support for adaptive mesh refinement, and adding support for a changing time-step.

\appendix

\section{Maxwell-MLA algorithms}\label{sec:algorithm}

Pseudo-code algorithms for the second-order accurate and fourth-order accurate Maxwell-MLA
time-stepping schemes are given in this section.

-------------------------------------------------------------------------------------
\subsection{Second-order accurate Maxwell-MLA algorithm}\label{sec:algorithmOrder2}

\begin{algorithm}
\algFontSize 
\caption{Second-order accurate Maxwell-MLA algorithm}
\begin{algorithmic}[1]
    \State Initialize; 
    \While{$t<T_{final}$} \Comment Begin time-stepping loop
        \For{$i=1,\dots, nDomains$} \Comment Loop over each subdomain
          \For{$m=1,\dots, \mathcal{N}_p$}
              \State $\Pv^{n+1}_{m,\jv}=\frac{1}{1+b_{1,m}\frac{\Delta t}{2}}\left(2\Pv^n_{m,\jv}-\Pv^{n-1}_{m,\jv} +b_{1,m} \frac{\Delta t}{2}\Pv^{n-1}_{m,\jv}-\Delta t^2b_{0,m}\Pv^n_{m,\jv} + \Delta t^2a_{m,\ell}N_{\ell}^n\Ev^n_{\jv}\right)$; 
              \State $D_{0t} \Pv^n_{m,\jv} = ({\Pv^{n+1}_{m,{\jv}}-\Pv^{n-1}_{m,{\jv}}})/(2\Delta t)$;
              \State $D_{+t}D_{-t} \Pv^n_{m,\jv} = ({\Pv^{n+1}_{m,{\jv}}-2\Pv^{n+1}_{m,{\jv}}+\Pv^{n-1}_{m,{\jv}}})/{\Delta t^2}$;
          \EndFor
          \State $D_{+t}D_{-t}\Ev^n_{\jv}=({\Ev^{n+1}_{\jv}-2\Ev^n_{\jv}+\Ev^{n-1}_{\jv}})/{\Delta t^2}$;
          
          \State $\Ev^{n+1}_{\jv}=2\Ev^n_{\jv}-\Ev^{n-1}_{\jv} + \Delta t^2c^2\Delta_{2h}\Ev^n_{\jv}-\epsilon^{-1}_0\sum_{m=1}^{\mathcal{N}_p}\left(\Pv^{n+1}_{m,\jv}-2\Pv^n_{m,\jv}+\Pv^{n-1}_{m,\jv}\right)$
          \For{$\ell=0,\dots, \mathcal{N}_\ell-1$}
              \State $\left.D_{2t}N_\ell\right|^n_{\jv} = \alpha_{\ell,k}N^n_{k,\jv} + \beta_{\ell,m}\Ev^n_{\jv}\cdot D_{0t} \Pv^n_{m,\jv}$;
              \State $\left.D_{2tt}N^{n}_\ell\right|_{\jv}^n = \alpha_{\ell,k}\left.D_{2t}N_k\right|^n_{\jv} + \beta_{\ell,m}D_{0t}\Ev^n_{\jv}\cdot D_{0t} \Pv^n_{m,\jv}+ \beta_{\ell,m}\Ev^n_{\jv}\cdot D_{+t}D_{-t} \Pv^n_{m,\jv}$;
              \State $N_{\ell,\jv}^{n+1}=N_{\ell,\jv}^{n}+\Delta t\left.D_{2t}N_{\ell}\right|^n_{\jv}+\frac{\Delta t^2}{2}\left.D_{2tt}N_{\ell}\right|^n_{\jv}$;
          \EndFor
        \EndFor
        \State Apply boundary and interface conditions;

        \State $t = t + \Delta t$, $n = n + 1$;
    \EndWhile    \Comment End time-stepping loop
\end{algorithmic} 
\label{alg:mbe2}
\end{algorithm}

\begin{algorithm}
\algFontSize 
\caption{Second-order scheme for filling ghost point values at interface}
\begin{algorithmic}[1]
    \State Extrapolate first ghost lines to second-order accuracy;
    \State Assemble and store the coefficient matrix of the linear system arising from the {\color{blue}second-order} discretizations \eqref{eqn:ic2-1-alt},\eqref{eqn:ic2-2},\eqref{eqn:ic2-3} and \eqref{eqn:ic2-4};
    \For{$\bm{j}\in\Gamma_h$} 
        \State Evaluate dispersive forcing functions, i.e., $D_+D_-\Pv|^n_{\bm{j}}$ by taking a fictitious forward step;
        \State Evaluate derivatives on the left hand side of \eqref{eqn:ic2-1-alt},\eqref{eqn:ic2-2},\eqref{eqn:ic2-3} and \eqref{eqn:ic2-4} using {\color{blue}second-order} accurate schemes with ``wrong'' values in the ghost points;
        \State Evaluate the {\color{blue}second-order} accurate residual of the jump conditions;
        \State Adjust right hand sides of interface conditions by subtracting off the ``wrong'' ghost values;
        \State Solve local linear system and fill in ghost values along the normal line;
    \EndFor
\end{algorithmic} 
\label{alg:ic2}
\end{algorithm}

\subsection{Fourth-order accurate Maxwell-MLA algorithm}\label{sec:algorithmOrder4}

\begin{algorithm}
\algFontSize 
\caption{Algorithm for fourth order approximations}
\begin{algorithmic}[1]
    \State Initialize; 
    \While{$t<T_{final}$} \Comment Begin time-stepping loop
        \For{$i=1,\dots, nDomains$} \Comment Loop over each subdomain
          \For{$m=1,\dots, \mathcal{N}_p$} 
              \State $\Pv^{n+1,*}_{m,\jv}=\frac{1}{1+b_{1,m}\frac{\Delta t}{2}}\left(2\Pv^n_{m,\jv}-\Pv^{n-1}_{m,\jv} +b_{1,m} \frac{\Delta t}{2}\Pv^{n-1}_{m,\jv}-\Delta t^2b_{0,m}\Pv^n_{m,\jv} + \Delta t^2a_{m,\ell}N_{\ell}^n\Ev^n_{\jv}\right)$; 
          \EndFor
          \State $\Ev^{n+1,*}_{\jv}=2\Ev^n_{\jv}-\Ev^{n-1}_{\jv} + \Delta t^2c^2\Delta_{2h}\Ev^n_{\jv}-\epsilon^{-1}_0\sum_{m=1}^{\mathcal{N}_p}\left(\Pv^{n+1,*}_{m,\jv}-2\Pv^n_{m,\jv}+\Pv^{n-1}_{m,\jv}\right)$
          \For{$\ell=0,\dots, \mathcal{N}_\ell-1$}
              \State ${\blue \left.D_{2t}N^*_\ell\right|^{n}_{\jv}} = \alpha_{\ell,k}N^n_{k,\jv} + \beta_{\ell,m}\Ev^n_{\jv}\cdot {\blue D_{0t} \Pv^{n,*}_{m,\jv}}$;
              \State ${\blue \left.D_{2tt}N^*_\ell\right|_{\jv}^{n}} = \alpha_{\ell,k}{\blue \left.D_{2t}N^*_k\right|^n_{\jv}} + \beta_{\ell,m}{\blue D_{0t}\Ev^n_{\jv}}\cdot {\blue D_{0t} \Pv^{n,*}_{m,\jv}}+ \beta_{\ell,m}\Ev^n_{\jv}\cdot {\blue D_{+t}D_{-t} \Pv^{n,*}_{m,\jv}}$;
          \EndFor
          \For{$m=1,\dots, \mathcal{N}_p$} 

              \State ${\blue \left.D_{2ttt} \Pv^*_m\right|_{\jv}^{n}} = -b_{1,m}{\blue D_{+t}D_{-t} \Pv^{n,*}_{m,\jv}} - b_{0,m}{\blue D_{0t}\Pv^{n,*}_{m,\jv}} + a_{m,\ell}{\blue \left.D_{2t} N^{*}_{\ell}\right|_{\jv}^n} \Ev^n_{\jv} +a_{k,\ell}N^n_{\ell,\jv} {\blue D_{0t}\Ev^{n,*}_{\jv}}$;

              \State ${\blue \left.D_{2tttt} \Pv^*_m\right|_{\jv}^{n}} = -b_{1,m}{\blue \left.D_{2ttt} \Pv^*_m\right|_{\jv}^{n}} - b_{0,m}{\blue D_{+t}D_{-t} \Pv^{n,*}_{m,\jv}} + a_{m,\ell}{\blue \left.D_{2tt} N^{*}_\ell\right|_{\jv}^n} \Ev^n_{\jv} $
              \State \qquad\qquad $+ 2a_{m,\ell}{\blue \left.D_{2t} N^{*}_{\ell}\right|_{\jv}^n} {\blue D_{0t}\Ev^{n,*}_{\jv}}+a_{m,\ell}N^n_{\ell,\jv} {\blue D_{+t}D_{-t}\Ev^{n,*}_{\jv}}$;

              \State $\Pv^{n+1}_{m,\jv}=\frac{1}{1+b_{1,m}\frac{\Delta t}{2}}\Big(2\Pv^n_{m,\jv}-\Pv^{n-1}_{m,\jv}+\frac{\Delta t^4}{12}{\blue \left.D_{2tttt}\Pv^*_{m}\right|_{\jv}^n} +\frac{\Delta t}{2}b_{1,m} \Pv^{n-1}_{m,\jv} $
              \State \qquad\qquad $+\frac{\Delta t^4}{6}b_{1,m}{\blue \left.D_{2ttt}\Pv^*_m\right|_{\jv}^n}-\Delta t^2b_{0,m}\Pv^n_{m,\jv} + \Delta t^2a_{m,\ell}N_{\ell,\jv}^n\Ev^n_{\jv}\Big)$; 

              \State ${\blue D_{+t}D_{-t}\Delta_{2h}\Pv^n_{m,\jv}}=({\Delta_{2h}\Pv^{n+1}_{m,\jv}-2\Delta_{2h}\Pv^n_{m,\jv}+\Delta_{2h}\Pv^{n-1}_{m,\jv}})/{\Delta t^2}$; 

              \State $\left.D_{4t}\Pv_m\right|_{\jv}^{n}=D_{0t} \Pv^n_{m,\jv}-\frac{\Delta t^2}{6}{\blue \left.D_{2ttt}\Pv^*_m\right|_{\jv}^{n}}$;

              \State $\left.D_{4tt}\Pv_m\right|_{\jv}^{n}=-b_{1,m} \left.D_{4t}\Pv_m\right|_{\jv}^{n} - b_{0,m}\Pv^n_{m,\jv}+a_{m,\ell}N_{\ell,\jv}^n\Ev^n_{\jv}$; 
          \EndFor

          \State $\Ev^{n+1}_{\jv}= 2\Ev^{n}_{\jv}-\Ev^{n-1}_{\jv}+\frac{\Delta t^4}{12}\left(c^4\Delta_{2h}^2\Ev^n_{\jv}-\sum_{m=1}^{\mathcal{N}_p}\epsilon^{-1}_0c^2{\blue D_{+t}D_{-t}\Delta_{2h}\Pv^n_{m,\jv}}\right) $
          \State \qquad\qquad  $+ \Delta t^2c^2\Delta_{4h}\Ev^n_{\jv}-\epsilon^{-1}_0\sum_{m=1}^{\mathcal{N}_p}\left(\Pv^{n+1}_{m,\jv}-2\Pv^{n}_{m,\jv}+\Pv^{n-1}_{m,\jv}\right)$;

          \State ${\blue \left.D_{2ttt}\Ev\right|_{\jv}^n} = \frac{1}{2\Delta t}\left[c^2\Delta_{2h}\Ev^{n+1}_{\jv}-c^2\Delta_{2h}\Ev^{n-1}_{\jv}\right]-\epsilon_0^{-1}\sum_{m=1}^{\mathcal{N}_p}{\blue \left.D_{2ttt}\Pv\right|^{n+1}_{m,\jv}}$;

          \State $\left.D_{4t}\Ev\right|_{\jv}^{n}=D_{0t} \Ev^n_{\jv}-\frac{\Delta t^2}{6}{\blue \left.D_{2ttt}\Ev\right|_{\jv}^{n}}$; 

          \For{$\ell=0,\dots, \mathcal{N}_\ell-1$} 
              \State $\left.D_{4t}N_\ell\right|_{\jv}^{n} = \alpha_{\ell,k}N^n_k + \beta_{\ell,m}\Ev^n_{\jv}\cdot \left.D_{4t}\Pv_m\right|_{\jv}^{n}$;
              \State $\left.D_{4tt}N_\ell\right|_{\jv}^{n} = \alpha_{\ell,k}\left.D_{4t}N_\ell\right|_{\jv}^{n} + \beta_{\ell,m}\left.D_{4t}\Ev\right|_{\jv}^{n}\cdot \left.D_{4t}\Pv_m\right|_{\jv}^{n}+ \beta_{\ell,m}\Ev^n\cdot \left.D_{4tt}\Pv_m\right|_{\jv}^{n}$;
          \EndFor

          \For{$m=1,\dots,\mathcal{N}_p$} 
              \State ${\blue \left.D_{2ttt} \Pv_m\right|_{\jv}^n} = -b_{1,m}D_{4tt} \Pv^{n}_{m,\jv} - b_{0,m}D_{4t}\Pv^{n}_{m,\jv} + a_{m,\ell}\left.D_{4t} N_{\ell}\right|_{\jv}^n \Ev^n_{\jv} +a_{m,\ell}N^n_{\ell,\jv} D_{4t}\Ev^{n}_{\jv}$;
              \State ${\blue \left.D_{2tttt} \Pv_m\right|_{\jv}^n} = -b_{1,m}{\blue \left.D_{2ttt} \Pv_m\right|_{\jv}^n} - b_{0,m}D_{4tt} \Pv^{n}_{m,\jv} + a_{m,\ell}\left.D_{4tt} N_\ell\right|_{\jv}^n \Ev^n_{\jv}$ 
              \State \qquad\qquad $+2a_{m,\ell}\left.D_{4t} N_{\ell}\right|_{\jv}^n D_{4t}\Ev^{n}_{\jv}+a_{m,\ell}N^n_{\ell,\jv} D_{+t}D_{-t}\Ev^{n}_{\jv}$;
          \EndFor

          \For{$\ell=0,\dots, \mathcal{N}_\ell-1$} 
              \State ${\blue \left.D_{2ttt}N_\ell\right|_{\jv}^{n}} = \alpha_{\ell,k}\left.D_{4tt}N_k\right|_{\jv}^{n} + \beta_{\ell,m}D_{+t}D_{-t}\Ev^{n}_{\jv}\cdot \left.D_{4t}\Pv_m\right|_{\jv}^{n}+2\beta_{\ell,m}\left.D_{4t}\Ev\right|_{\jv}^{n}\cdot \left.D_{4tt}\Pv_m\right|_{\jv}^{n}$
              \State \qquad \quad $+ \beta_{\ell,m}\Ev^n_{\jv}\cdot {\blue \left.D_{2ttt}\Pv_m\right|_{\jv}^{n}}$;

              \State ${\blue \left.D_{2tttt}N_\ell\right|_{\jv}^{n}} = \alpha_{\ell,k}{\blue \left.D_{2ttt}N_k\right|_{\jv}^{n}}+ \beta_{\ell,m}{\blue \left.D_{2ttt}\Ev\right|_{\jv}^{n}}\cdot \left.D_{4t}\Pv_m\right|_{\jv}^{n} $ 
              \State \qquad\quad $+ 3\beta_{\ell,m}D_{+t}D_{-t}\Ev^n_{\jv}\cdot \left.D_{4tt}\Pv_m\right|_{\jv}^{n}+3\beta_{\ell,m}\left.D_{4t}\Ev\right|_{\jv}^{n}\cdot {\blue \left.D_{2ttt}\Pv_m\right|_{\jv}^{n}}+\beta_{\ell,m}\Ev^n_{\jv}\cdot {\blue \left.D_{2tttt}\Pv_m\right|_{\jv}^{n}}$;

              \State $N^{n+1}_{\ell,\jv}=N^n_{\ell,\jv}+\Delta t \left.D_{4t}N_\ell\right|_{\jv}^n +\frac{\Delta t^2}{2} \left.D_{4tt}N_\ell\right|_{\jv}^n+\frac{\Delta t^3}{6} {\blue \left.D_{2ttt}N_\ell\right|_{\jv}^n}+\frac{\Delta t^4}{24} {\blue \left.D_{2tttt}N_\ell\right|_{\jv}^n}$;
          \EndFor
        \EndFor
        \State Apply boundary and interface conditions;

        \State $t = t + \Delta t$, $n = n + 1$;
    \EndWhile    \Comment End time-stepping loop
\end{algorithmic} 
\label{alg:mbe4}
\end{algorithm}

\begin{algorithm}
\algFontSize 
\caption{Fourth-order scheme for filling ghost point values at interface}
\begin{algorithmic}[1]
    \State Extrapolate first and second ghost lines to fourth-order accuracy;
    \State Call Algorithm~\ref{alg:ic2} for points $\bm{j}\in\Gamma_h$ plus extra points, except to use {\color{magenta}fourth-order} residuals for evaluation;
    \State Assemble and store the coefficient matrix of the linear system arising from the {\color{magenta}fourth-order} discretizations \eqref{eqn:ic4-1-alt},\eqref{eqn:ic4-2}--\eqref{eqn:ic4-5},\eqref{eqn:ic4-6-alt}--\eqref{eqn:ic4-8-alt};
    \For{$\bm{j}\in\Gamma_h$} 
        \State Evaluate dispersive forcing functions, i.e., $D_+D_-\Pv|^n_{\bm{j}}$, etc. to their respective accuracy as in the discretized jump conditions by taking a fictitious forward step;
        \State Evaluate derivatives on the left hand side of \eqref{eqn:ic4-1-alt},\eqref{eqn:ic4-2}--\eqref{eqn:ic4-5},\eqref{eqn:ic4-6-alt}--\eqref{eqn:ic4-8-alt} using {\color{magenta}fourth-order} accurate schemes with ``wrong'' values in the ghost points;
        \State Evaluate the {\color{magenta}fourth-order} accurate residuals of the jump conditions;
        \State Adjust right hand sides of interface conditions by subtracting off the ``wrong'' ghost values;
        \State Solve local linear system and fill in first and second ghost values along the normal line;
    \EndFor
\end{algorithmic} 
\label{alg:ic4}
\end{algorithm}

\clearpage
\section{Well-posedness and long time stability of a restricted Maxwell-MLA system}\label{sec:stability}

In this section we consider the well-posedness and long-time stability of the Maxwell-MLA equations~\eqref{eqn:mbe}. 
If lower-order terms are dropped from the equations~\eqref{eqn:mbe}, Maxwell's equations for $\Ev$ decouples
from the polarization equations~\eqref{eqn:nonlin_p} and rate equations~\eqref{eqn:nonlin_n}. These equations are thus well posed 
with the appropriate initial conditions and boundary conditions. The solution to the IBVP will exist for
at least short times. To study the long-time existence of the nonlinear equations we restrict ourselves
to a class of Maxwell-MLA equations that are of common interest, see Figure~\ref{fig:atomicLevelsCartoon}. To simplify the discussion we consider
the equations written using the first-order form for $\Ev$ and $\Hv$
\bse\label{eqn:paired-system}
\ba
    &\eps_0\Ev_t = \nabla\times \Hv - \Pv_t,        \label{eqn:ps_e}\\
    &\mu_0 \Hv_t = -\nabla\times\Ev\label{eqn:ps_h} ,          
\ea
\ese
with
an MLA system of $\Nc_n$ energy levels, each with population density $N_\ell$, $\ell=0,1,\ldots,\Nc_n-1$
(see Figure~\ref{fig:4-level} for an example of a 4-level system). Polarization states $\Pv_{ji}$, with $i$ and $j$ integers in the range $0$ and $\Nc_n-1$ and $i < j$, 
may exist between any two levels 
with governing equation given by
\ba
   \partial_t^2 \Pv_{ji}+\gamma_{ji}\partial_t\Pv_{ji}+\omega_{ji}^2\Pv_{ji}=\kappa_{ji}(N_i-N_j)\Ev,  \qquad ji \in \Tc,  \label{eq:Pjieqn}
\ea
 where $ji$ belongs to the set of active transitions pairs, $\Tc$. For example, in the four-level system in Figure~\ref{fig:4-level}, $\Tc=\{ 30, 21 \}$.
 The parameters in~\eqref{eq:Pjieqn} are assumed to satisfy $\gamma_{ji}\ge 0$, $\omega_{ji} > 0 $ and  
 $\kappa_{ji} > 0$. 

{
\newcommand{\smallss}{\sffamily\small}
\newcommand{\xa}{-1.5}
\newcommand{\xb}{9}
\definecolor{ghostColour}{named}{DodgerBlue}
\begin{figure}[H]
\centering
\begin{tikzpicture}[scale=1]
   \useasboundingbox (.0,0.0) rectangle (9,6);  

    \foreach \y in {0,1,2,3,4,5,6}
        \draw[blue,dashed,very thick] (\xa,\y) -- (\xb,\y) node[anchor=west,xshift=5pt] {\smallss $N_{\y}$};

    \draw[<->,red,very thick] (0,0) -- (0,6) node[anchor=east,xshift=1pt,yshift=-10pt] {\smallss $P_{60}$};
    \draw[<->,red,very thick,xshift=1cm] (0,0) -- (0,3) node[anchor=east,xshift=1pt,yshift=-10pt] {\smallss $P_{30}$};

    \draw[<->,red,very thick,xshift=3cm] (0,2) -- (0,5) node[anchor=east,xshift=1pt,yshift=-10pt] {\smallss $P_{52}$};
    \draw[<->,red,very thick,xshift=4cm] (0,1) -- (0,5) node[anchor=east,xshift=1pt,yshift=-10pt] {\smallss $P_{51}$};
    \draw[<->,red,very thick,xshift=5.5cm] (0,1) -- (0,4) node[anchor=east,xshift=1pt,yshift=-10pt] {\smallss $P_{41}$};

    \draw[->,violet,dashed,very thick] (.25,6) -- (.5,4) node[anchor=south,xshift=10pt,yshift=12pt] {\smallss $\alpha_{6,4}$};
    \draw[->,violet,dashed,very thick,xshift=.5cm] (.25,6) -- (1.7,2) node[anchor=south,xshift=6pt,yshift=12pt] {\smallss $\alpha_{6,2}$};
    \draw[->,violet,dashed,very thick,xshift=1cm] (.25,6) -- (.75,5) node[anchor=south,xshift=10pt,yshift=12pt] {\smallss $\alpha_{6,5}$};

    \draw[->,violet,dashed,very thick,xshift=4cm] (.25,5) -- (.75,1) node[anchor=south,xshift=10pt,yshift=12pt] {\smallss $\alpha_{5,1}$};

    \draw[->,violet,dashed,very thick,xshift=1cm] (.25,0) -- (.75,2) node[anchor=east,xshift=20pt,yshift=-16pt] {\smallss $\alpha_{0,2}$};

    \filldraw[red] (-.5,3) node[draw=white,fill=white,rotate=90] {\smallss $\kappa_{60}(N_0-N_6) \Ev$};

    \begin{scope}[xshift=1cm]
      \draw[<->,blue,very thick,xshift=6.5cm] (0,1) -- (0,4);
      \draw[<->,blue,very thick,xshift=7.0cm] (0,4) -- (0,1) ;
      \draw[blue,xshift=6.75cm] (0,1) node[anchor=north]  {\smallss $\f{1}{\hbar\omega_{4,1}} \Ev\cdot\p_t\Pv_{4,1}$};
    \end{scope}

\end{tikzpicture}
\caption{Jablonski diagram for a common class of MLA systems showing the energy levels, population densities $N_i$ and selected transitions. 
Polarization states $P_{ji}$ can exist between any two atomic levels with $i<j$.
The $\alpha_{ji}$ are relaxation time-constants for relaxation from state $j$ to state $i$.
The rate of change of state $N_4$ depends on the source term $\f{1}{\hbar\omega_{4,1}} \Ev\cdot\p_t\Pv_{4,1}$
while the rate of change of state $N_1$ depends on the same source term with opposite sign.
The source term $\kappa_{60}(N_0-N_6) \Ev$ for the $P_{60}$ ODE is proportional to $\Ev$ and the difference $N_0-N_6$.
}
\label{fig:atomicLevelsCartoon}
\end{figure}
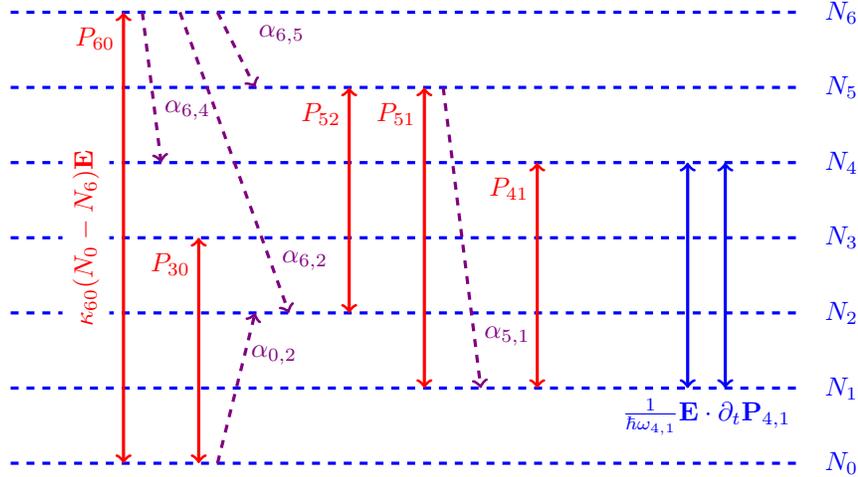
}

The population density $N_\ell$ satisfies 
\ba
 & \partial_tN_\ell = \sum_{k=0}^{\Nc_n-1} \alpha_{\ell k}N_k  + 
        \sum_{ji\in\Tc}  \sigma_{ji}\, \frac{1}{\hbar\omega_{ji}}\Ev\cdot\partial_t\Pv_{ji}, \label{eq:Nelleqn}\\
  & \sigma_{ij} = \begin{cases} +1, & \text{if $\ell=j$,} \\
                                -1, & \text{if $\ell=i$,} \\
                                 0, & \text{otherwise.} \\
                  \end{cases}
\ea 
Our goal is to derive an $L_2$-energy estimate to show that the solution remains bounded in time. 
Let $(f,g)$ denote the $L_2$-inner product over $\Omega$,
\ba 
   (f,g) = \int_\Omega f(\xv) \, g(\xv) \, d\xv,
\ea
for scalar functions $f$ and $g$ (all functions are assume to be real valued).
Let $\| \cdot \|$ denote the corresponding norm. For vector functions we use 
\ba 
   ( \fv , \gv ) = \int_\Omega \fv(\xv) \cdot  \, \gv(\xv) \, d\xv .
\ea 
 In the usual way we take inner products of the various equations with the corresponding
variable or its time-derivative, 
\bat
   & ( \Ev, \eps_0 \Ev_t)  = ( \Ev, \nabla\times \Hv - \p_t \Pv ),           \label{eq:EE-E}   \\
   & ( \Hv, \mu_0  \Hv_t )    = (\Hv, -\nabla\times\Ev),                     \label{eq:EE-H}   \\
   & ( \p_t \Pv_{ji}, \p_t^2 \Pv_{ji}+\gamma_{ji}\p_t\Pv_{ji}+\omega_{ji}^2\Pv_{ji} ) = (\p_t\Pv_{ji},\kappa_{ji}(N_i-N_j)\Ev),
              \quad ij \in \Tc, \label{eq:EE-P} \\
   & ( N_\ell, \partial_t{N_\ell} ) =  \SumkNn (N_\ell,\alpha_{\ell k}N_k)  
               + \sum_{ij \in\Tc}  (N_\ell ,\sigma_{ji}\, \frac{1}{\hbar\omega_{ji}}\Ev\cdot \p_t \Pv_{ji}) 
              \quad \ell =0,1,2,\ldots,\Nc_n-1 .\label{eq:EE-N} 
\eat
Integrating by parts the right hand sides of~\eqref{eq:EE-E} and~\eqref{eq:EE-H} and then adding these equations gives
an equation for the time-derivative of isotropic energy $\Ec_0 \eqdef \f{\eps_0}{2} \| \Ev\|^2 + \f{\mu_0}{2} \| \Hv\|^2$,
\ba
   \p_t\Big( \f{\eps_0}{2} \| \Ev\|^2 + \f{\mu_0}{2} \| \Hv\|^2 \Big) & = - (\Ev,\p_t \Pv) + BT.s, \label{eq:EE-EH}
\ea
where $BT.s$ denotes the usual boundary terms for isotropic Maxwell's equations. We assume the boundary conditions are chosen
to make the boundary terms to vanish or be negative. 
Let us now focus on equations~\eqref{eq:EE-P} and~\eqref{eq:EE-N} which can be written as
\ba
   & \half \p_t \| \p_t \Pv_{ji} \|^2  + \gamma_{ji} \| \p_t \Pv_{ji} \|^2 + \half\omega_{ji}^2 \, \p_t \| \Pv_{ji}\|^2  
             = (\p_t\Pv_{ji}, \kappa_{ji}(N_i-N_j)\Ev), \label{eq:EE-Pb} \\
   & \half \p_t \| N_\ell \|^2 = \SumkNn (N_\ell,\alpha_{\ell k}N_k) + \sum_{ji\in\Tc} (N_\ell ,\sigma_{ji}\, \frac{1}{\hbar\omega_{ji}}\Ev\cdot \p_t \Pv_{ji})    \label{eq:EE-Nb} 
\ea
\newcommand{\khoi}{\delta_{ji}}
\newcommand{\khoSum}{\Kc}
Define the quantities $\khoSum$ and $\khoi$ by 
\ba
   &  \khoSum \eqdef  \sum_{ij\in\Tc} \kappa_{ji} \hbar\omega_{ji}, \\
   &  \khoi \eqdef \f{\khoSum}{\kappa_{ji} \hbar\omega_{ji}} 
\ea
Scaling equation~\eqref{eq:EE-Pb} by $\khoi$ and equation~\eqref{eq:EE-Nb} by $\khoSum$ leads to 
\ba
   & \khoi\left\{ \half \p_t \| \p_t \Pv_{ji} \|^2  + \gamma_{ji} \| \p_t \Pv_{ji} \|^2 + \half\omega_{ji}^2 \, \p_t \| \Pv_{ji}\|^2 \right\} 
             =    ( N_i-N_j , \frac{\khoSum}{\hbar\omega_{ji}} \Ev \cdot \p_t\Pv_{ji}), \label{eq:EE-PbScaled} \\
   & \khoSum \half \p_t \| N_\ell \|^2 = 
      \khoSum  \SumkNn (N_\ell,\alpha_{\ell k}N_k) + \sum_{ji\in\Tc} ( \sigma_{ji}\, N_\ell , \frac{\khoSum}{\hbar\omega_{ji}}\Ev\cdot \p_t\Pv_{ji})   
            \label{eq:EE-NbScaled}
\ea
Adding~\eqref{eq:EE-NbScaled} to the sum over $ji\in\Tc$ of~\eqref{eq:EE-PbScaled} eliminates the inner products containing the
nonlinear terms $\Ev\cdot\p_t \Pv_{ji}$ to give

\ba
 & \sum_{ji \in \Tc} \khoi \left\{ \half \p_t \| \p_t \Pv_{ji} \|^2  + \gamma_{ji} \| \p_t \Pv_{ji} \|^2 + \half\omega_{ji}^2 \, \p_t \| \Pv_{ji}\|^2 \right\}
    + \khoSum \SumiNn \half \p_t \| N_\ell \|^2     \label{eq:EE-PN}  \\
    &\qquad\qquad = \khoSum\SumiNn \SumkNn (N_\ell,\alpha_{\ell k}N_k)  \nonumber
\ea
Let $\Ec_{PN}$ be defined from the terms on the left of~\eqref{eq:EE-PN}
\ba 
   \Ec_{PN} &\eqdef \sum_{ji \in \Tc} \khoi \left\{ \half \| \p_t \Pv_{ji} \|^2  + \half\omega_{ji}^2 \, \| \Pv_{ji}\|^2 \right\}
    + \khoSum  \half \| \Nv \|^2,                        \label{eq:EPNdef}
\ea
where
\ba
  \| \Nv \|^2 \eqdef \sum_{\ell=0}^{\mathcal{N}_n-1}  \| N_\ell \|^2.
\ea  
\newcommand{\Kpt}{K_p}
Note that 
\ba
  & \|\p_t \Pv \|^2 = \sum_{ji \in \Tc}  \| \p_t \Pv_{ji} \|^2    \le \Kpt \Ec_{PN},   \label{eq:PtvBound} \\
  &    \Kpt \eqdef  \min_{ji\in\Tc} \f{1}{\khoi}. 
\ea
Equation~\eqref{eq:EE-PN} becomes
\ba
  \p_t \Ec_{PN} = - \sum_{ji \in \Tc} \khoi \gamma_{ji} \| \p_t \Pv_{ji} \|^2 + \Kc \, \sum_{\ell=0}^{\mathcal{N}_n-1} \SumkNn (N_\ell,\alpha_{\ell k}N_k) \label{eq:E-PNbound}
\ea
Using $|N_i N_j| \le \half N_i^2 + \half N_j^2$ implies
\ba
   \Kc \, \Big|  \sum_\ell \sum_k (N_\ell,\alpha_{\ell k}N_k) \Big| \le \Ca \, \Kc \, \half \| \Nv \|^2 \le \Ca \Ec_{PN},  \label{eq:Nbound}
\ea
for some constant $\Ca$ which depends on $\alphaMax = \max_{ji} |\alpha_{ji}|$, $\Nc_n$ and $\Nc_p$. 
Using~\eqref{eq:Nbound} in~\eqref{eq:E-PNbound} implies
\ba
    & \p_t \Ec_{PN}  \le  \Ca \Ec_{PN}, 
\ea
Integrating this last expression in time implies $\Ec_{PN}$ is bounded in time,
\ba
    \Ec_{PN}(t) \le \Ec_{PN}(0) e^{\Ca t} ,
\ea 
which also, from the definition~\eqref{eq:EPNdef} for $\Ec_{PN}$ and~\eqref{eq:PtvBound} implies a bound on $\| \p_t \Pv \|^2$ and $\| \Nv\|^2$
\ba
   & \| \p_t \Pv \|^2 \le \Kpt \, \Ec_{PN}(0) e^{\Ca t}, \\
   &  \f{\Kc}{2} \| \Nv \|^2 \le \Ec_{PN}(0) e^{\Ca t}.
\ea

We are now prepared to form the final energy estimate.
Adding equations~\eqref{eq:EE-EH} and~\eqref{eq:E-PNbound} gives an equation for the total energy $\Ec \eqdef \Ec_0 + \Ec_{PN}$
\ba
  \p_t \Ec = - (\Ev,\p_t \Pv) - \sum_{ji \in \Tc} \khoi \,\gamma_{ji} \| \p_t \Pv_{ji} \|^2 + \Kc \SumiNn \SumkNn (N_i,\alpha_{ik}N_k)  + BT.s  , \label{eq:EE-total}
\ea
Using
\ba
    \big| (\Ev,\p_t \Pv) \big| \le \f{\eps_0}{2} \| \Ev \|^2 + \f{1}{2\eps_0} \| \p_t \Pv \|^2 
\ea
implies 
\ba
   \p_t \Ec & \le \f{\eps_0}{2} \| \Ev\|^2 + \f{1}{2\eps_0} \|\p_t \Pv\|^2 + C \| \Nv\|^2 + BT.s, \\
            & \le   \Ec  + C_2 \, \Ec_{PN}(0) e^{\Ca t}    \label{eq:EcBoundI}
\ea
for some constant $C_2$, where we have assumed the boundary terms are non-positive. 
To get a bound for $\Ec$ we integrate the inequality~\eqref{eq:EcBoundI} 
to give 
\ba
  \Ec(t) \le e^{t} \Ec(0) + C_2 \, e^t \int_0^t e^{ (\Ca-1)\tau}\, \Ec_{PN}(0) \, d\tau, \\
        =  e^{t} \Ec(0) + \f{C_2}{1-\Ca} \Big(  e^t - e^{\Ca t} \Big) \, \Ec_{PN}(0) . 
\ea
where the case $\Ca=1$ can be found with the appropriate limit. 

\medskip\noindent
We have therefore proved the following theorem.
\begin{theorem}
Given appropriate boundary conditions, 
the following $L_2$-``energy'' of the Maxwell-MLA system~\eqref{eqn:paired-system},\eqref{eq:Pjieqn},\eqref{eq:Nelleqn}
 \ba
   \Ec = \f{\eps_0}{2} \| \Ev\|^2 + \f{\mu_0}{2} \| \Hv\|^2 + 
    \sum_{ji \in \Tc} \khoi \left\{ \half \| \p_t \Pv_{ji} \|^2 +  \half\omega_{ji}^2 \, \| \Pv_{ji}\|^2 \right\}
      + \f{\Kc}{2} \| \Nv \|^2,  
\ea
has bounded exponential growth in time,
\ba
  \Ec(t) \le K_1 e^{K_2 t} ,
\ea
for some constants $K_1$ and $K_2$.  
\end{theorem}

\section{Supplemental Equations for Fourth-Order}
\label{sec:4thsupp}


The following expressions are used in Section~\ref{sec:4th} in the description
of the fourth-order accurate scheme.
\begin{subequations}
\begin{align}
  {\color{blue}\left.D_{2ttt} \Pv^*_m\right|_{\jv}^n} & \eqdef
    - b_{1,m}{\color{blue}D_{+t}D_{-t} \Pv^{n,*}_{m,\jv}} 
    - b_{0,m}{\color{blue}D_{0t}\Pv^{n,*}_{m,\jv}} \nonumber \\
  &\qquad + \sum_{\ell}a_{m,\ell}{\color{blue}\left.D_{2t} N^{*}_{\ell}\right|_{\jv}^n} \Ev^n_{\jv} 
    + \sum_{\ell}a_{k,\ell}N^n_{\ell,\jv} {\color{blue}D_{0t}\Ev^{n,*}_{\jv}}\label{eqn:Pttt}\\
  {\color{blue}\left.D_{2tttt} \Pv^*_m\right|_{\jv}^n} & \eqdef
    - b_{1,m}{\color{blue}\left.D_{2ttt} \Pv^*_m\right|_{\jv}^n}  
    - b_{0,m}{\color{blue}D_{+t}D_{-t} \Pv^{n,*}_{m,\jv}} 
    + \sum_{\ell}a_{m,\ell}{\color{blue}\left.D_{2tt} N^{*}_\ell\right|_{\jv}^n} \Ev^n_{\jv} \nonumber \\
  &\qquad +2\sum_{\ell}a_{m,\ell}{\color{blue}\left.D_{2t} N^{*}_{\ell}\right|_{\jv}^n} {\color{blue}D_{0t}\Ev^{n,*}_{\jv}}
    + \sum_{\ell}a_{m,\ell}N^n_{\ell,\jv} {\color{blue}D_{+t}D_{-t}\Ev^{n,*}_{\jv}}\label{eqn:Ptttt}\\
  \left.D_{4t}N_\ell\right|_{\jv}^{n} & \eqdef
    \sum_{\hat{\ell}}\alpha_{\ell,\hat{\ell}}N^n_{\hat{\ell}} 
    + \beta_{\ell,m}\Ev^n_{\jv}\cdot \left.D_{4t}\Pv_m\right|_{\jv}^{n}\\
  \left.D_{4tt}N_\ell\right|_{\jv}^{n} & \eqdef
    \sum_{\hat{\ell}}\alpha_{\ell,\hat{\ell}}\left.D_{4t}N_\ell\right|_{\jv}^{n} 
    + \beta_{\ell,m}\left.D_{4t}\Ev\right|_{\jv}^{n}\cdot \left.D_{4t}\Pv_m\right|_{\jv}^{n}
    + \beta_{\ell,m}\Ev^n\cdot \left.D_{4tt}\Pv_m\right|_{\jv}^{n}\\
  {\color{blue}\left.D_{2ttt}N_\ell\right|_{\jv}^{n}} & \eqdef
    \sum_{\hat{\ell}}\alpha_{\ell,\hat{\ell}}\left.D_{4tt}N_{\hat{\ell}}\right|_{\jv}^{n} 
    + \beta_{\ell,m}D_{+t}D_{-t}\Ev^{n}_{\jv}\cdot \left.D_{4t}\Pv_m\right|_{\jv}^{n}\nonumber\\
  &\qquad+ 2\beta_{\ell,m}\left.D_{4t}\Ev\right|_{\jv}^{n}\cdot \left.D_{4tt}\Pv_m\right|_{\jv}^{n}
    + \beta_{\ell,m}\Ev^n_{\jv}\cdot {\color{blue}\left.D_{2ttt}\Pv_m\right|_{\jv}^{n}}\\
  {\color{blue}\left.D_{2tttt}N_\ell\right|_{\jv}^{n}} & \eqdef
    \sum_{\hat{\ell}}\alpha_{\ell,\hat{\ell}}{\color{blue}\left.D_{2ttt}N_{\hat{\ell}}\right|_{\jv}^{n}}
    + \beta_{\ell,m}{\color{blue}\left.D_{2ttt}\Ev\right|_{\jv}^{n}}\cdot \left.D_{4t}\Pv_m\right|_{\jv}^{n}
    + 3\beta_{\ell,m}D_{+t}D_{-t}\Ev^n_{\jv}\cdot \left.D_{4tt}\Pv_m\right|_{\jv}^{n}\nonumber\\
  &\qquad +3\beta_{\ell,m}\left.D_{4t}\Ev\right|_{\jv}^{n}\cdot {\color{blue}\left.D_{2ttt}\Pv_m\right|_{\jv}^{n}}
    +\beta_{\ell,m}\Ev^n_{\jv}\cdot {\color{blue}\left.D_{2tttt}\Pv_m\right|_{\jv}^{n}}
\end{align}
\end{subequations}



\section{Material definitions}\label{sec:materialDefinitions}

Here are the definitions of the MLA materials used in this article.

\subsection{Active material mlaMat2}  \label{sec:mlaMat2}
\bni
Active material mlaMat2 has two polarization vectors and four atomic levels with MLA material parameters
given by
\bse 
\ba
  & \Nc_p =2, \quad \Nc_n=4, \\
  & \eps_0 = 1, \quad \mu_0=1, \\
  & \am= \begin{bmatrix*}[r]
         2.3418 &    0.0 &     0.0 &  2.3418   \\
            0.0 & 11.666 &  11.666 &     0.0   
       \end{bmatrix*}, \\
  & \bmm= \begin{bmatrix*}[r]
            1 & .1 \\
            1 & .1
         \end{bmatrix*} , \\       
  & \alpham = \begin{bmatrix*}[r]
          0.0 &   0.0010542 &            0.0 &   0.000000012723 & \\
          0.0 &  -0.0010542 &   0.0000014641 &              0.0 &  \\
          0.0 &         0.0 &  -0.0000014641 &        0.0012299 & \\
          0.0 &         0.0 &            0.0 &       -0.0012299 &
             \end{bmatrix*}, 
  \\
 &   \betam = \begin{bmatrix*}[r]
            -2.3418 &     0.0 \\
                0.0 & -2.4362 \\
                0.0 &  2.4362 \\
             2.3418 &     0.0   
              \end{bmatrix*} .
\ea
\ese
These parameters are adapted from \cite{Gain0DNanohubTool} for tests with manufactured solutions.

\subsection{Active material mlaMat3} \label{sec:mlaMat3}
\bni
Active material mlaMat3 has one polarization vector and one atomic level with MLA material parameters
given by
\bse 
\ba
  & \Nc_p =1, \quad \Nc_n=1, \\
  & \eps_0 = 2, \quad \mu_0=1, \\
  & \am= \begin{bmatrix*}[r]
         10 
       \end{bmatrix*},     \\
  & \bmm= \begin{bmatrix*}[r]
            1 & 0 \\
            1 & 0
         \end{bmatrix*} , \\    
  & \alpham = \begin{bmatrix*}[r]
          0.01
             \end{bmatrix*}, 
  \\
 &   \betam = \begin{bmatrix*}[r]
               1 
              \end{bmatrix*} .
\ea
\ese


\subsection{Active material mlaMat4levels} \label{sec:mlaMat4levels}
\bni
Active material mlaMat4levels has two polarization vectors and four atomic levels with MLA material parameters
given by 
\bse 
\ba
  & \Nc_p =2, \quad \Nc_n=4, \\
  & \eps_0 = 2, \quad \mu_0=1, \\
  & \am= \begin{bmatrix*}[r]
            2.3418 &    0.0 &      0.0 &  -2.3418 \\
               0.0 & 11.666 &  -11.666 &      0.0 
              \end{bmatrix*},   \\
  & \bmm= \begin{bmatrix*}[r]
           769.2308 & 64.0180  \\
           710.8037 & 152.1820
         \end{bmatrix*} , \\    
  & \alpham = \begin{bmatrix*}[r]
                  0.0 &   0.0010542 &            0.0 &   0.000000012723 \\
                  0.0 &  -0.0010542 &   0.0000014641 &              0.0 \\
                  0.0 &         0.0 &  -0.0000014641 &        0.0012299 \\
                  0.0 &         0.0 &            0.0 &       -0.0012299  
             \end{bmatrix*}, 
  \\
 &   \betam = \begin{bmatrix*}[r]
                  -1801.421965974931 &                0.0 \\
                                 0.0 & -1873.997239424281 \\
                                 0.0 &  1873.997239424281 \\
                   1801.421965974931 &                0.0 
              \end{bmatrix*} .
\ea
\ese
These parameters are adapted from \cite{Gain0DNanohubTool} with proper undimensionalization to match the form of \eqref{eqn:mbe}.

\section{Maxwell-MLA interface algorithm, second-order accuracy} \label{sec:InterfaceMLAOrder2Algorithm}

\mni
Guide to some variables:
\begin{flushleft}
   \texttt{u1(i1,i2,i3,0:2)}     : holds $E_x$, $E_y$ and $E_z$ at point $\iv=(i_1,i_2,i_3)$ on side 1, \\
   \texttt{u2(j1,j2,j3,0:2)}     : holds $E_x$, $E_y$ and $E_z$ at point $\jv=(j_1,j_2,j_3)$ on side 2, \\
   \texttt{rsxy1(i1,i2,i3,m,n)}  : $\displaystyle\f{\p r_m}{\p x_n}$ for a point $\iv=(i_1,i_2,i_3)$ on side 1, \\
   \texttt{rsxy2(j1,j2,j3,m,n)}  : $\displaystyle\f{\p r_m}{\p x_n}$ for a point $\jv=(j_1,j_2,j_3)$ on side 2, \\
   \texttt{axis1=0,1,2,3}        : interface side 1 is located on $r_{\rm axis1}=side1$ constant, side1=0,1\\
   \texttt{axis2=0,1,2,3}        : interface side 2 is located on $r_{\rm axis2}=side2$ constant, side2=0,1\\
   \texttt{is1 = 1 -2*side1}     : index shift vector on side1=0 (left) or side1=1 (right), \\
   \texttt{is2 = 1 -2*side2}     : index shift vector on side1=0 (left) or side1=1 (right), \\
   \texttt{dx1(0:2), dx2(0:2)}   : grid spacing for rectangular (Cartesian) grids on side 1 and side 2, \\
   \texttt{dr1(0:2), dr2(0:2)}   : unit-square grid spacing for curvilinear grids on side 1 and side 2, \\
\end{flushleft}

\subsection{Second-order accurate interfaces, Cartesian grids}

Here is a listing of the fortran code that defines the Maxwell-MLA interface algorithm
for second-order accuracy on Cartesian grids. The code is written with the bpp macro preprocessor from Overture.

\begin{lstlisting}[numbers=left, stepnumber=1, firstline=1]
! ---------------------------------------------------------------------------------------
! Macro: Assign nonlinear DISPERSIVE interface ghost values, DIM=2, ORDER=2, GRID=Rectangular
! 
! Here are the jump conditions (See notes in DMX_ADE)
!   [ u.x + v.y ] = 0
!   [ (1/mu)* tv,.( curl(E) ) ]
!   [ tv.( c^2*Delta(E) -alphaP*P_tt) ] = 0  --> [ tv.( beta*c^2*Delta(E) - alphaP* F) ]=0 
!   [ (1/mu)* nv.( Delta(E) ) ]=0
! 
! -------------------------------------------------------------------------------------------
#beginMacro assignNonlinearInterfaceGhost22r()

 ! ****************************************************
 ! ***********  2D, ORDER=2, RECTANGULAR **************
 ! ****************************************************

INFO("22rectangle-nonlinear-MLA") ! nonlinear multilevel atomic system

! For rectangular, both sides must axis axis1==axis2: 
if( axis1.ne.axis2 )then
  stop 8826
end if

! 
! Solve for the unknowns u1(-1),v1(-1),u2(-1),v2(-1)
!     
!       A [ U ] = A [ U(old) ] - [ f ]
!
!               [ u1(-1) ]
!       [ U ] = [ v1(-1) ]
!               [ u2(-1) ]
!               [ v2(-1) ]
!             

! --- initialize some forcing functions ---
! forcing functions for E and P
do n=0,nd-1
  fev1(n)=0.
  fev2(n)=0.
  if (dispersionModel1 .ne. noDispersion) then
    do jv=0,numberOfPolarizationVectors1-1
      fpv1(n,jv)=0.
    end do
  endif
  if (dispersionModel2 .ne. noDispersion) then
    do jv=0,numberOfPolarizationVectors2-1
      fpv2(n,jv)=0.
    end do
  endif
end do
! forcing functions for N
if (nonlinearModel1 .ne. noNonlinearModel) then
  do jv = 0,numberOfAtomicLevels1-1
      fnv1(jv) = 0.
      fntv1(jv) = 0.
  enddo
endif
if (nonlinearModel2 .ne. noNonlinearModel) then
  do jv = 0,numberOfAtomicLevels2-1
      fnv2(jv) = 0.
      fntv2(jv) = 0.
  enddo
endif

! print *, "-----------Now using MLA (RECTANGULAR)---------------"

! ----------------- START LOOP OVER INTERFACE -------------------------
beginLoopsMask2d()

  ! u1 = Ex, v1 = Ey on side 1
  ! u2 = Ex, v2 = Ey on side 2

  ! Evaluate derivatives of the solution, u1x, u1y, v1x, v1y using the wrong values at the ghost points:
  evalInterfaceDerivatives2d()
 
  ! Evaluate TZ forcing for dispersive equations in 2D 
  getTZForcingMLA(fpv1,fpv2,fev1,fev2,fnv1,fntv1,fnv2,fntv2)

  ! eval dispersive forcings for domain 1
  getMLAForcingOrder2(i1,i2,i3, fp1, fpv1,fev1,fnv1,fntv1,p1,p1n,p1m,q1,q1n,q1m, u1,u1n,u1m, dispersionModel1,nonlinearModel1,numberOfPolarizationVectors1,numberOfAtomicLevels1,alphaP1,beta1,pnec1,prc1,peptc1,b0v1,b1v1)

  ! eval dispersive forcings for domain 2
  getMLAForcingOrder2(j1,j2,j3, fp2, fpv2,fev2,fnv2,fntv2,p2,p2n,p2m,q2,q2n,q2m, u2,u2n,u2m, dispersionModel2,nonlinearModel2,numberOfPolarizationVectors2,numberOfAtomicLevels2,alphaP2,beta2,pnec2,prc2,peptc2,b0v2,b1v2)


  if( axis1.eq.0 )then ! vertical interfaces
    ! Interface equations for a boundary at x = 0 or x=1

    ! ---- EQUATION 0 -----
    ! [ u.x + v.y ] = 0
    ! NOTE: if mu==mu2 then we do not need TZ forcing for this eqn:
    f(0)=(u1x+v1y) - \
         (u2x+v2y)
    a4(0,0) = -is1/(2.*dx1(axis1))    ! coeff of u1(-1) from [u.x+v.y] 
    a4(0,1) = 0.                      ! coeff of v1(-1) from [u.x+v.y] 
    a4(0,2) =  js1/(2.*dx2(axis2))    ! coeff of u2(-1) from [u.x+v.y] 
    a4(0,3) = 0.                      ! coeff of v2(-1) from [u.x+v.y]
  
    ! ---- EQUATION 1 -----
    ! [ (1/mu)* tv,.( curl(E) ) ] = 0  
    ! NOTE: if mu==mu2 then we do not need TZ forcing for this eqn:
    f(1)=(v1x-u1y)/mu1 - \
         (v2x-u2y)/mu2
    a4(1,0) = 0.
    a4(1,1) = -is1/(2.*dx1(axis1))    ! coeff of v1(-1) from [v.x - u.y] 
    a4(1,2) = 0.
    a4(1,3) =  js1/(2.*dx2(axis2))    ! coeff of v2(-1) from [v.x - u.y]
   
    ! ---- EQUATION 2 -----
    ! [ (1/mu)* nv.( Delta(E) ) ]=0 (normal component)
    ! NOTE: if mu1==mu2 then we do not need TZ forcing for this eqn (TZ forcing canceled due to nonzero jump conditions)
    f(2)=( (u1xx+u1yy)/mu1  ) - \
         ( (u2xx+u2yy)/mu2  )
    a4(2,0) = 1./(dx1(axis1)**2)/mu1   ! coeff of u1(-1) from [(u.xx + u.yy)/mu]
    a4(2,1) = 0. 
    a4(2,2) =-1./(dx2(axis2)**2)/mu2   ! coeff of u2(-1) from [(u.xx + u.yy)/mu]
    a4(2,3) = 0. 
  
    ! ---- EQUATION 3 -----    
    ! [ tv.( c^2*Delta(E) -alphaP*P_tt) ] = 0 (tangential component)
    ! The coefficient of Delta(E) in this equation is altered due to Ptt term (not true for MLA)
    f(3)=( (v1xx+v1yy)*beta1/epsmu1 -alphaP1*fp1(1) + fev1(1)) - \
         ( (v2xx+v2yy)*beta2/epsmu2 -alphaP2*fp2(1) + fev2(1))


    a4(3,0) = 0.                      
    a4(3,1) = (beta1/epsmu1)/(dx1(axis1)**2) ! coeff of v1(-1) from [beta*c^2*(v.xx+v.yy)]
    a4(3,2) = 0. 
    a4(3,3) =-(beta2/epsmu2)/(dx2(axis2)**2) ! coeff of v2(-1) from [beta*c^2*(v.xx+v.yy)]

    ! print *, 'E TZ forcing (x)',fev1(0),fev2(0),'E TZ forcing (y)',fev1(1),fev2(1)

    ! print *, '============eps:', eps1,eps2, 'mu',mu1,mu2, 'epsmu',epsmu1,epsmu2,'beta',beta1,beta2,'alphaP',alphaP1,alphaP2

  else ! ---------- horizontal interfaces ---------------

    ! Interface equations for a boundary at y = 0 or y=1
    ! Switch u <-> v,  x<-> y in above equations 

    ! ---- EQUATION 0 -----
    f(0)=(v1y+u1x) - \
         (v2y+u2x)
    a4(0,0) = 0.                      ! coeff of u1(-1) from [u.x+v.y] 
    a4(0,1) = -is1/(2.*dx1(axis1))    ! coeff of v1(-1) from [u.x+v.y] 

    a4(0,2) = 0.                      ! coeff of u2(-1) from [u.x+v.y] 
    a4(0,3) = js1/(2.*dx2(axis2))     ! coeff of v2(-1) from [u.x+v.y]
  
    ! ---- EQUATION 1 -----
    f(1)=(u1y-v1x)/mu1 - \
         (u2y-v2x)/mu2
    a4(1,0) = -is1/(2.*dx1(axis1))
    a4(1,1) = 0.
    a4(1,2) =  js1/(2.*dx2(axis2))  
    a4(1,3) = 0.
   
    ! ---- EQUATION 2 -----    
    f(2)=( (v1xx+v1yy)/mu1 ) - \
         ( (v2xx+v2yy)/mu2 )
    a4(2,0) = 0.
    a4(2,1) = 1./(dx1(axis1)**2)/mu1  
    a4(2,2) = 0.
    a4(2,3) =-1./(dx2(axis2)**2)/mu2 
  
    ! ---- EQUATION 3 -----    
    ! The coefficient of Delta(E) in this equation is altered due to Ptt term 
    f(3)=( (u1xx+u1yy)*beta1/epsmu1 -alphaP1*fp1(0) +fev1(0) ) - \
         ( (u2xx+u2yy)*beta2/epsmu2 -alphaP2*fp2(0) +fev2(0) )
    a4(3,0) = (beta1/epsmu1)/(dx1(axis1)**2)
    a4(3,1) = 0.
    a4(3,2) =-(beta2/epsmu2)/(dx2(axis2)**2) 
    a4(3,3) = 0.


  end if


   q(0) = u1(i1-is1,i2-is2,i3,ex)
   q(1) = u1(i1-is1,i2-is2,i3,ey)
   q(2) = u2(j1-js1,j2-js2,j3,ex)
   q(3) = u2(j1-js1,j2-js2,j3,ey)


   ! subtract off the contributions from the wrong values at the ghost points:
   do n=0,3
     f(n) = (a4(n,0)*q(0)+a4(n,1)*q(1)+a4(n,2)*q(2)+a4(n,3)*q(3)) - f(n)
   end do

   ! solve A Q = F
   ! factor the matrix
   numberOfEquations=4
   call dgeco( a4(0,0), numberOfEquations, numberOfEquations, ipvt(0),rcond,work(0))

   ! solve
   job=0
   call dgesl( a4(0,0), numberOfEquations, numberOfEquations, ipvt(0), f(0), job)

   u1(i1-is1,i2-is2,i3,ex)=f(0)
   u1(i1-is1,i2-is2,i3,ey)=f(1)
   u2(j1-js1,j2-js2,j3,ex)=f(2)
   u2(j1-js1,j2-js2,j3,ey)=f(3)


 endLoopsMask2d()
 
 
#endMacro
\end{lstlisting}

\subsection{Second-order accurate interfaces, curvilinear grids} \label{sec:mlaInterface22c}

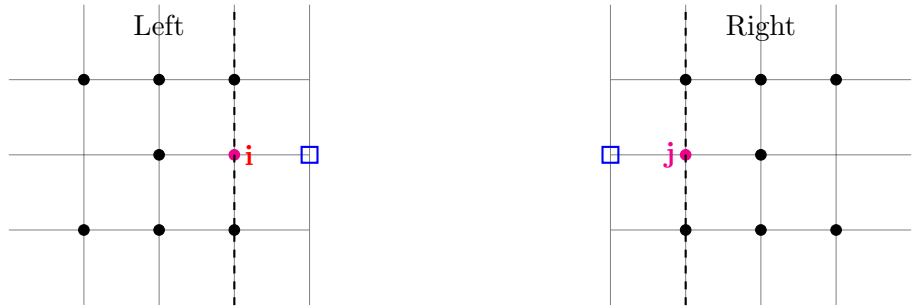
\begin{figure}[H]
\centering
\begin{tikzpicture}[scale=2]
\draw[step=0.5cm,gray,very thin] (0+0.001,0+0.001) grid (2,2-0.001);
\filldraw[magenta] (1.5,1) circle (1pt);
\filldraw (1.5,1.5) circle (1pt);
\filldraw (1,1.5) circle (1pt);
\filldraw (0.5,1.5) circle (1pt);
\filldraw (1.5,0.5) circle (1pt);
\filldraw (1,0.5) circle (1pt);
\filldraw (0.5,0.5) circle (1pt);
\filldraw (1.,1) circle (1pt);
\draw[blue,thick] (2.,1) +(-1.5pt,-1.5pt) rectangle +(1.5pt,1.5pt);
\draw[dashed,thick] (1.5,0) -- (1.5,2);
\draw[magenta] (1.5,1) node[right] { \color{red}$\iv$};
\draw (1,2) node[below] {Left};
\draw[step=0.5cm,gray,very thin] (4-0.001,0+0.001) grid (6-0.001,2-0.001);
\filldraw[magenta] (4.5,1) circle (1pt);
\filldraw (4.5,1.5) circle (1pt);
\filldraw (5,1.5) circle (1pt);
\filldraw (5.5,1.5) circle (1pt);
\filldraw (4.5,0.5) circle (1pt);
\filldraw (5,0.5) circle (1pt);
\filldraw (5.5,0.5) circle (1pt);
\filldraw (5,1) circle (1pt);
\draw[blue,thick] (4,1) +(-1.5pt,-1.5pt) rectangle +(1.5pt,1.5pt);
\draw[dashed,thick] (4.5,0) -- (4.5,2);
\draw[magenta] (4.5,1) node[left] { \color{magenta}$\jv$};
\draw (5,2) node[below] {Right};
\end{tikzpicture}
\caption{Local stencil at point $\iv$ (left) and $\jv$ (right) on the interface for the second order accurate numerical interface approximations (dashed line -- interface, dots -- known interior values, square -- ghost points to be filled)}
\label{fig:ic2-stencil}
\end{figure}

Here is a listing of the fortran code that defines the Maxwell-MLA interface algorithm
for second-order accuracy on curvilinear grids. The code is written with the bpp macro preprocessor from Overture.
The stencil for the code appears in Figure~\ref{fig:ic2-stencil}.

\begin{lstlisting}[numbers=left, stepnumber=1, firstline=1]
! --------------------------------------------------------------------
! Macro: Assign NONLINEAR interface ghost values, DIM=2, ORDER=2, GRID=Curvilinear
! 
! Here are the jump conditions (See notes in DMX_ADE)
!   [ u.x + v.y ] = 0
!   [ (1/mu)* tv,.( curl(E) ) ]
!   [ tv.( c^2*Delta(E) -alphaP*P_tt) ] = 0  --> [ tv.( beta*c^2*Delta(E) - alphaP* F) ]=0 
!   [ (1/mu)* nv.( Delta(E) ) ]=0
! 
! -------------------------------------------------------------------------------------------
#beginMacro assignNonlinearInterfaceGhost22c()

  ! ****************************************************
  ! ***********  2D, ORDER=2, CURVILINEAR **************
  ! ****************************************************

INFO("22curvilinear-nonlinear-MLA")

! --- initialize some forcing functions ---
do n=0,nd-1
  fev1(n)=0.
  fev2(n)=0.
  if (dispersionModel1 .ne. noDispersion) then
    do jv=0,numberOfPolarizationVectors1-1
      fpv1(n,jv)=0.
    end do
  endif
  if (dispersionModel2 .ne. noDispersion) then
    do jv=0,numberOfPolarizationVectors2-1
      fpv2(n,jv)=0.
    end do
  endif
end do
! forcing functions for N
if (nonlinearModel1 .ne. noNonlinearModel) then
  do jv = 0,numberOfAtomicLevels1-1
      fnv1(jv) = 0.
      fntv1(jv) = 0.
  enddo
endif
if (nonlinearModel2 .ne. noNonlinearModel) then
  do jv = 0,numberOfAtomicLevels2-1
      fnv2(jv) = 0.
      fntv2(jv) = 0.
  enddo
endif

! ----------------- START LOOP OVER INTERFACE -------------------------
beginLoopsMask2d()

  ! here is the normal (assumed to be the same on both sides)
  an1=rsxy1(i1,i2,i3,axis1,0)   ! normal (an1,an2)
  an2=rsxy1(i1,i2,i3,axis1,1)
  aNorm=max(epsx,sqrt(an1**2+an2**2))
  an1=an1/aNorm
  an2=an2/aNorm
  tau1=-an2
  tau2= an1

  ! first evaluate the derivatives of the solution using the wrong values at the ghost points:
  evalInterfaceDerivatives2d()
 
  ! Evaluate TZ forcing for dispersive equations in 2D 
  getTZForcingMLA(fpv1,fpv2,fev1,fev2,fnv1,fntv1,fnv2,fntv2)

  ! eval dispersive forcings for domain 1
  getMLAForcingOrder2(i1,i2,i3, fp1, fpv1,fev1,fnv1,fntv1,p1,p1n,p1m,q1,q1n,q1m, u1,u1n,u1m, dispersionModel1,nonlinearModel1,numberOfPolarizationVectors1,numberOfAtomicLevels1,alphaP1,beta1,pnec1,prc1,peptc1,b0v1,b1v1)

  ! eval dispersive forcings for domain 2
  getMLAForcingOrder2(j1,j2,j3, fp2, fpv2,fev2,fnv2,fntv2,p2,p2n,p2m,q2,q2n,q2m, u2,u2n,u2m, dispersionModel2,nonlinearModel2,numberOfPolarizationVectors2,numberOfAtomicLevels2,alphaP2,beta2,pnec2,prc2,peptc2,b0v2,b1v2)

  ! Evaulate RHS, f(n),n=0,1,2,3 using current ghost values: 
  eval2dJumpMLAOrder2()
  
  ! here is the matrix of coefficients for the unknowns u1(-1),v1(-1),u2(-1),v2(-1)
  ! Solve:
  !     
  !       A [ U ] = A [ U(old) ] - [ f ]
  ! ---- EQUATION 0 ----- 
  a4(0,0) = -is*rsxy1(i1,i2,i3,axis1,0)/(2.*dr1(axis1))    ! coeff of u1(-1) from [u.x+v.y] 
  a4(0,1) = -is*rsxy1(i1,i2,i3,axis1,1)/(2.*dr1(axis1))    ! coeff of v1(-1) from [u.x+v.y] 
  a4(0,2) =  js*rsxy2(j1,j2,j3,axis2,0)/(2.*dr2(axis2))    ! coeff of u2(-1) from [u.x+v.y] 
  a4(0,3) =  js*rsxy2(j1,j2,j3,axis2,1)/(2.*dr2(axis2))    ! coeff of v2(-1) from [u.x+v.y] 

  ! ---- EQUATION 2 ----- 
  a4(2,0) =  is*rsxy1(i1,i2,i3,axis1,1)/(2.*dr1(axis1))/mu1   ! coeff of u1(-1) from [(v.x - u.y)/mu] 
  a4(2,1) = -is*rsxy1(i1,i2,i3,axis1,0)/(2.*dr1(axis1))/mu1   ! coeff of v1(-1) from [(v.x - u.y)/mu] 

  a4(2,2) = -js*rsxy2(j1,j2,j3,axis2,1)/(2.*dr2(axis2))/mu2   ! coeff of u2(-1) from [(v.x - u.y)/mu] 
  a4(2,3) =  js*rsxy2(j1,j2,j3,axis2,0)/(2.*dr2(axis2))/mu2   ! coeff of v2(-1) from [(v.x - u.y)/mu] 


  ! coeff of u(-1) from lap = u.xx + u.yy
  rxx1(0,0,0)=aj1rxx
  rxx1(1,0,0)=aj1sxx
  rxx1(0,1,1)=aj1ryy
  rxx1(1,1,1)=aj1syy

  rxx2(0,0,0)=aj2rxx
  rxx2(1,0,0)=aj2sxx
  rxx2(0,1,1)=aj2ryy
  rxx2(1,1,1)=aj2syy

  
  clap1=(rsxy1(i1,i2,i3,axis1,0)**2+rsxy1(i1,i2,i3,axis1,1)**2)/(dr1(axis1)**2) \
            -is*(rxx1(axis1,0,0)+rxx1(axis1,1,1))/(2.*dr1(axis1))
  clap2=(rsxy2(j1,j2,j3,axis2,0)**2+rsxy2(j1,j2,j3,axis2,1)**2)/(dr2(axis2)**2) \
            -js*(rxx2(axis2,0,0)+rxx2(axis2,1,1))/(2.*dr2(axis2)) 

  ! ---- EQUATION 1 ----- 
  !   [ n.(uv.xx + u.yy)/mu ] = 0
  a4(1,0) = an1*clap1/mu1
  a4(1,1) = an2*clap1/mu1
  a4(1,2) =-an1*clap2/mu2
  a4(1,3) =-an2*clap2/mu2

  ! ---- EQUATION 3 ----- 
  !   [ tau.(uv.xx+uv.yy)*beta/(eps*mu) + ... ] = 0
  a4(3,0) = tau1*clap1*beta1/epsmu1
  a4(3,1) = tau2*clap1*beta1/epsmu1
  a4(3,2) =-tau1*clap2*beta2/epsmu2
  a4(3,3) =-tau2*clap2*beta2/epsmu2
    

  q(0) = u1(i1-is1,i2-is2,i3,ex)
  q(1) = u1(i1-is1,i2-is2,i3,ey)
  q(2) = u2(j1-js1,j2-js2,j3,ex)
  q(3) = u2(j1-js1,j2-js2,j3,ey)


  ! subtract off the contributions from the wrong values at the ghost points:
  do n=0,3
    f(n) = (a4(n,0)*q(0)+a4(n,1)*q(1)+a4(n,2)*q(2)+a4(n,3)*q(3)) - f(n)
  end do

  ! solve A Q = F
  ! factor the matrix
  numberOfEquations=4
  call dgeco( a4(0,0), numberOfEquations, numberOfEquations, ipvt(0),rcond,work(0))

  ! solve
  job=0
  call dgesl( a4(0,0), numberOfEquations, numberOfEquations, ipvt(0), f(0), job)
 
  u1(i1-is1,i2-is2,i3,ex)=f(0)
  u1(i1-is1,i2-is2,i3,ey)=f(1)
  u2(j1-js1,j2-js2,j3,ex)=f(2)
  u2(j1-js1,j2-js2,j3,ey)=f(3)

endLoopsMask2d()



#endMacro

\end{lstlisting}

\section{Maxwell-MLA interface algorithm, fourth-order accuracy} \label{sec:InterfaceMLAOrder4Algorithm}

\begin{figure}[H]
\centering
\begin{tikzpicture}[scale=1.5]
\draw[step=0.5cm,gray,very thin] (-1+0.001,-1+0.001) grid (2.5,2-0.001);
\filldraw[magenta] (1.5,0.5) circle (1pt);
\filldraw (1.5,1) circle (1pt);
\filldraw (1,1) circle (1pt);
\filldraw (0.5,1) circle (1pt);
\filldraw (0,1) circle (1pt);
\filldraw (-0.5,1) circle (1pt);
\filldraw (1.5,0) circle (1pt);
\filldraw (1,0) circle (1pt);
\filldraw (0.5,0) circle (1pt);
\filldraw (0,0) circle (1pt);
\filldraw (-0.5,0) circle (1pt);
\filldraw (1.,0.5) circle (1pt);
\filldraw (0.5,0.5) circle (1pt);
\filldraw (1.5,-0.5) circle (1pt);
\filldraw (1,-0.5) circle (1pt);
\filldraw (0.5,-0.5) circle (1pt);
\filldraw (0,-0.5) circle (1pt);
\filldraw (-0.5,-0.5) circle (1pt);
\filldraw (1.5,1.5) circle (1pt);
\filldraw (1,1.5) circle (1pt);
\filldraw (0.5,1.5) circle (1pt);
\filldraw (0,1.5) circle (1pt);
\filldraw (-0.5,1.5) circle (1pt);
\draw[blue,thick] (2.,0.5) +(-1.5pt,-1.5pt) rectangle +(1.5pt,1.5pt);
\draw[thick,red] (2,1.) circle (1pt);
\draw[thick,red] (2,0.5) circle (1pt);
\draw[thick,red] (2,0) circle (1pt);
\draw[blue,thick] (2.5,0.5) +(-1.5pt,-1.5pt) rectangle +(1.5pt,1.5pt);
\draw[dashed,thick] (1.5,-1) -- (1.5,2);
\draw[magenta] (1.5,0.5) node[right] { \color{magenta}$\iv$};
\draw (1,2) node[below] {Left};
\draw[step=0.5cm,gray,very thin] (3.5-0.001,-1+0.001) grid (7-0.001,2-0.001);
\filldraw[magenta] (4.5,0.5) circle (1pt);
\filldraw (4.5,1) circle (1pt);
\filldraw (5,1) circle (1pt);
\filldraw (5.5,1) circle (1pt);
\filldraw (6,1) circle (1pt);
\filldraw (6.5,1) circle (1pt);
\filldraw (4.5,0) circle (1pt);
\filldraw (5,0) circle (1pt);
\filldraw (5.5,0) circle (1pt);
\filldraw (6,0) circle (1pt);
\filldraw (6.5,0) circle (1pt);
\filldraw (5,0.5) circle (1pt);
\filldraw (5.5,0.5) circle (1pt);
\filldraw (4.5,-0.5) circle (1pt);
\filldraw (5,-0.5) circle (1pt);
\filldraw (5.5,-0.5) circle (1pt);
\filldraw (6,-0.5) circle (1pt);
\filldraw (6.5,-0.5) circle (1pt);
\filldraw (4.5,1.5) circle (1pt);
\filldraw (5,1.5) circle (1pt);
\filldraw (5.5,1.5) circle (1pt);
\filldraw (6,1.5) circle (1pt);
\filldraw (6.5,1.5) circle (1pt);
\draw[blue,thick] (4,0.5) +(-1.5pt,-1.5pt) rectangle +(1.5pt,1.5pt);
\draw[thick,red] (4,1) circle (1pt);
\draw[thick,red] (4,0.5) circle (1pt);
\draw[thick,red] (4,0) circle (1pt);
\draw[blue,thick] (3.5,0.5) +(-1.5pt,-1.5pt) rectangle +(1.5pt,1.5pt);
\draw[dashed,thick] (4.5,-1) -- (4.5,2);
\draw[magenta] (4.5,0.5) node[left] { \color{magenta}$\jv$};
\draw (5,2) node[below] {Right};
\end{tikzpicture}
\caption{Local stencil at point $\iv$ (left) or $\jv$ (right) on the interface for the fourth order accurate numerical interface approximations after \texttt{EP-decoupling} and \texttt{tangential decoupling} (dashed line -- interface, dots -- known interior values, square -- ghost points to be filled, circle -- second order prediction)}
\label{fig:ic4-stencil}
\end{figure}
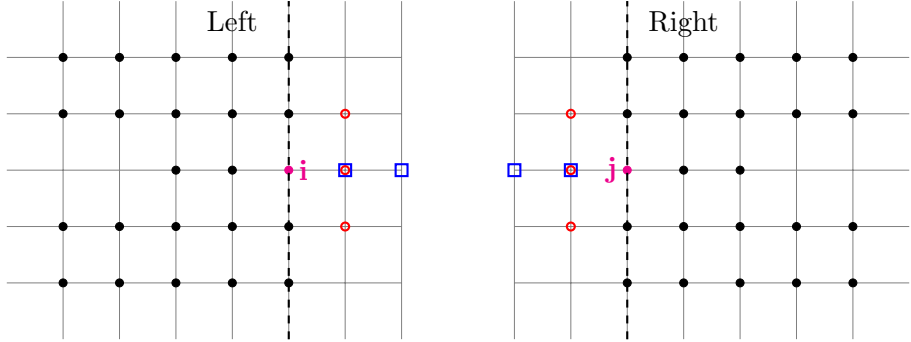

\mni
Guide to some variables:
\begin{flushleft}
   \texttt{u1(i1,i2,i3,0:2)}     : holds $E_x$, $E_y$ and $E_z$ at point $\iv=(i_1,i_2,i_3)$ on side 1, \\
   \texttt{u2(j1,j2,j3,0:2)}     : holds $E_x$, $E_y$ and $E_z$ at point $\jv=(j_1,j_2,j_3)$ on side 2, \\
   \texttt{rsxy1(i1,i2,i3,m,n)}  : $\displaystyle\f{\p r_m}{\p x_n}$ for a point $\iv=(i_1,i_2,i_3)$ on side 1, \\
   \texttt{rsxy2(j1,j2,j3,m,n)}  : $\displaystyle\f{\p r_m}{\p x_n}$ for a point $\jv=(j_1,j_2,j_3)$ on side 2, \\
   \texttt{axis1=0,1,2,3}        : interface side 1 is located on $r_{\rm axis1}=side1$ constant, side1=0,1\\
   \texttt{axis2=0,1,2,3}        : interface side 2 is located on $r_{\rm axis2}=side2$ constant, side2=0,1\\
   \texttt{is1 = 1 -2*side1}     : index shift vector on side1=0 (left) or side1=1 (right), \\
   \texttt{is2 = 1 -2*side2}     : index shift vector on side1=0 (left) or side1=1 (right), \\
   \texttt{dx1(0:2), dx2(0:2)}   : grid spacing for rectangular (Cartesian) grids on side 1 and side 2, \\
   \texttt{dr1(0:2), dr2(0:2)}   : unit-square grid spacing for curvilinear grids on side 1 and side 2, \\
\end{flushleft}

\subsection{Fourth-order accurate hierachical scheme, curvilinear grids}

Here is a listing of the fortran code that defines the Maxwell-MLA interface algorithm
for fourth-order accuracy on curvilinear grids. The code is written with the bpp macro preprocessor from Overture.
The following code implements the hierarchical interface approach by first using the second-order
accurate algorithm~\ref{sec:mlaInterface22c} to fill in the first line of ghost points before using the fourth-order accurate algorithm
given in~\ref{sec:mlaInterface24c}. Note that the residuals in the second-order accurate algorithm are being
computed with the fourth-order accurate approximations to the derivatives.
The stencil for the code is given in Figure~\ref{fig:ic4-stencil}

\begin{lstlisting}[numbers=left, stepnumber=1, firstline=1]
  else if( nd.eq.2 .and. orderOfAccuracy.eq.4 .and. gridType.eq.curvilinear )then
   #perl $DIM=2; $GRIDTYPE="curvilinear"; $ORDER=4;
    ! --------------- 4th Order Curvilinear ---------------
    ! ---- first satisfy the jump conditions on the boundary --------
    !    [ eps n.u ] = 0
    !    [ tau.u ] = 0
    !    [ w ] = 0 
    if( assignInterfaceValues.eq.1 )then
      boundaryJumpConditions(2,curvilinear)
    end if

   ! ----- assign ghost using jump conditions -----
   if( assignInterfaceGhostValues.eq.1 )then
    ! here are the real jump conditions for the ghost points
    ! 0  [ u.x + v.y ] = 0
    ! 1  [ n.(uv.xx + uv.yy) ] = 0
    ! 2  [ v.x - u.y ] =0 
    ! 3  [ tau.(v.xx+v.yy)/eps ] = 0
    ! 4  [ (u.xx+u.yy).x + (v.xx+v.yy).y ] = 0  OR [ (u.xx).x + (v.xx).y ] = 0 OR  [ (u.yy).x + (v.yy).y ] = 0 
    ! 5  [ {(Delta v).x - (Delta u).y}/eps ] =0  -> [ {(v.xxx+v.xyy)-(u.xxy+u.yyy)}/eps ] = 0
    ! 6  [ n.Delta^2 uv/eps ] = 0
    ! 7  [ tau.Delta^2 uv/eps^2 ] = 0 
    ! initialization step: assign first ghost line by extrapolation
    ! NOTE: assign ghost points outside the ends
    
    extrapolateGhost2dOrder4()

    ! ---- STAGE I: assign first ghost to 2nd-order accuracy -----

    ! *Note*: we are using fourth-order accurate derivatives in the residuals for the
    ! second-order scheme

    ! In parallel we add extra points in the tangential direction on parallel boundaries
    ! (otherwise we would use extrapolated values which is probably ok) 

    setIndexBoundsExtraGhost()
    if( dispersive.eq.0 )then
      assignInterfaceGhost22c()
    else if( useNonlinearModel.eq.0 )then
      assignDispersiveInterfaceGhost22c()
    else
      assignNonlinearInterfaceGhost22c()
    end if
    resetIndexBounds()         

    ! ---- STAGE II: assign two ghost to 4th-order accuracy -----
    ! Macro to assign ghost values:
    if( dispersive.eq.0 )then
      assignInterfaceGhost24c()
    else if (useNonlinearModel.eq.0) then
      ! dispersive case
      assignDispersiveInterfaceGhost24c()
    else
      assignNonlinearInterfaceGhost24c()
    end if    
    periodicUpdate2d(u1,boundaryCondition1,gridIndexRange1,side1,axis1)
    periodicUpdate2d(u2,boundaryCondition2,gridIndexRange2,side2,axis2)
    end if ! ----- end assign ghost using jump conditions -----
  end if

\end{lstlisting}

\subsection{assignNonlinearInterfaceGhost24c} \label{sec:mlaInterface24c}

Here is a listing of the fortran code that defines the Maxwell-MLA interface algorithm
for fourth-order accuracy on curvlinear grids. The code is written with the bpp macro preprocessor from Overture.

\begin{lstlisting}[numbers=left, stepnumber=1, firstline=1]
! --------------------------------------------------------------------------
! Macro: Assign interface ghost values, DIM=2, ORDER=4, GRID=Curvilinear
!         NONLINEAR DISPERSIVE CASE -- MLA 
! --------------------------------------------------------------------------
#beginMacro assignNonlinearInterfaceGhost24c()

 ! ****************************************************************
 ! ***********  DISPERSIVE, 2D, ORDER=4, CURVILINEAR **************
 ! ****************************************************************


! --- initialize some forcing functions ---
do n=0,nd-1
  fev1(n)=0.
  LfE1(n)=0.
  fEt1(n)=0.
  fEtt1(n)=0.

  fev2(n)=0.
  LfE2(n)=0.
  fEt2(n)=0.
  fEtt2(n)=0.

  fevx1(n)=0.
  fevy1(n)=0.
  fevx2(n)=0.
  fevy2(n)=0.
  if (dispersionModel1 .ne. 0) then
    do jv=0,numberOfPolarizationVectors1-1
      fpv1(n,jv)=0.
      LfP1(n,jv)=0.
      fPt1(n,jv)=0.
      fPtt1(n,jv)=0.

      fpvx1(n,jv)=0.
      fpvy1(n,jv)=0.
    end do
  endif
  if (dispersionModel2 .ne. 0) then
    do jv=0,numberOfPolarizationVectors2-1
      fpv2(n,jv)=0.
      LfP2(n,jv)=0.
      fPt2(n,jv)=0.
      fPtt2(n,jv)=0.

      fpvx2(n,jv)=0.
      fpvy2(n,jv)=0.
    end do
  endif
end do
! forcing functions for N
if (nonlinearModel1 .ne. 0) then
  do jv = 0,numberOfAtomicLevels1-1
      fnv1(jv) = 0.
      fntv1(jv) = 0.
  enddo
endif
if (nonlinearModel2 .ne. 0) then
  do jv = 0,numberOfAtomicLevels2-1
      fnv2(jv) = 0.
      fntv2(jv) = 0.
  enddo
endif


 ! =============== start loops ======================
 beginLoopsMask2d() 

   nn=nn+1

   ! here is the normal (assumed to be the same on both sides)
   an1=rsxy1(i1,i2,i3,axis1,0)   ! normal (an1,an2)
   an2=rsxy1(i1,i2,i3,axis1,1)
   aNorm=max(epsx,sqrt(an1**2+an2**2))
   an1=an1/aNorm
   an2=an2/aNorm
   tau1=-an2
   tau2= an1

   ! Evaluate the jump conditions using the wrong values at the ghost points 
   evaluateNonlinearInterfaceEquations2dOrder4()


     ! here is the matrix of coefficients for the unknowns u1(-1),v1(-1),u2(-1),v2(-1)
     ! Solve:
     !     
     !       A [ U ] = A [ U(old) ] - [ f ]


     ! write(debugFile,'(" interface:E: initialized,it=",2i4)') initialized,it
     if( .false. .or. (initialized.eq.0 .and. it.eq.1) )then
       ! form the matrix (and save factor for later use)
       if( nn.eq.0 )then
         write(*,'(" Interface42c: form matrix and factor, it=",i4)') it
       end if

       ! Equation 0: 
       ! 0  [ u.x + v.y ] = 0
       aa8(0,0,0,nn) = -is*8.*rsxy1(i1,i2,i3,axis1,0)*dr114(axis1)     ! coeff of u1(-1) from [u.x+v.y] 
       aa8(0,1,0,nn) = -is*8.*rsxy1(i1,i2,i3,axis1,1)*dr114(axis1)     ! coeff of v1(-1) from [u.x+v.y] 
       aa8(0,4,0,nn) =  is*   rsxy1(i1,i2,i3,axis1,0)*dr114(axis1)     ! u1(-2)
       aa8(0,5,0,nn) =  is*   rsxy1(i1,i2,i3,axis1,1)*dr114(axis1)     ! v1(-2) 

       aa8(0,2,0,nn) =  js*8.*rsxy2(j1,j2,j3,axis2,0)*dr214(axis2)     ! coeff of u2(-1) from [u.x+v.y] 
       aa8(0,3,0,nn) =  js*8.*rsxy2(j1,j2,j3,axis2,1)*dr214(axis2)  
       aa8(0,6,0,nn) = -js*   rsxy2(j1,j2,j3,axis2,0)*dr214(axis2) 
       aa8(0,7,0,nn) = -js*   rsxy2(j1,j2,j3,axis2,1)*dr214(axis2)  

      ! Equation 1:
      ! 1  [ u.xx + u.yy ] = 0
  
       setJacobian( aj1, axis1)

       dr0=dr1(axis1)
       ds0=dr1(axis1p1)
       aLap0 = lapCoeff4a(is,dr0,ds0)
       aLap1 = lapCoeff4b(is,dr0,ds0)


       ! dr1a(0:2) = dsBig in tangential directions if avoidInterfaceIterations=1
       ds0 =dr1a(axis1p1)
       aLapSq0 = lapSqCoeff4a(is,dr0,ds0)
       aLapSq1 = lapSqCoeff4b(is,dr0,ds0)

       setJacobian( aj2, axis2)
       dr0=dr2(axis2)
       ds0=dr2(axis2p1)
       bLap0 = lapCoeff4a(js,dr0,ds0)
       bLap1 = lapCoeff4b(js,dr0,ds0)

 
       ! dr2a(0:2) = dsBig in tangential directions if avoidInterfaceIterations=1
       ds0 = dr2a(axis2p1)
       bLapSq0 = lapSqCoeff4a(js,dr0,ds0)
       bLapSq1 = lapSqCoeff4b(js,dr0,ds0)


      ! Equation 1:
      aa8(1,0,0,nn) = an1*aLap0/mu1       ! coeff of u1(-1) from [n.(u.xx + u.yy)]
      aa8(1,1,0,nn) = an2*aLap0/mu1 
      aa8(1,4,0,nn) = an1*aLap1/mu1       ! coeff of u1(-2) from [n.(u.xx + u.yy)]
      aa8(1,5,0,nn) = an2*aLap1/mu1  
       
      aa8(1,2,0,nn) =-an1*bLap0/mu2       ! coeff of u2(-1) from [n.(u.xx + u.yy)]
      aa8(1,3,0,nn) =-an2*bLap0/mu2
      aa8(1,6,0,nn) =-an1*bLap1/mu2       ! coeff of u2(-2) from [n.(u.xx + u.yy)]
      aa8(1,7,0,nn) =-an2*bLap1/mu2

      ! Equation 2: 
      ! 2  [ v.x - u.y ] =0 
      !          a8(2,0) =  is*8.*ry1*dx114(axis1)
      !          a8(2,1) = -is*8.*rx1*dx114(axis1)    ! coeff of v1(-1) from [v.x - u.y] 
      !          a8(2,4) = -is*   ry1*dx114(axis1)
      !          a8(2,5) =  is*   rx1*dx114(axis1)
      !          a8(2,2) = -js*8.*ry2*dx214(axis2)
      !          a8(2,3) =  js*8.*rx2*dx214(axis2)
      !          a8(2,6) =  js*   ry2*dx214(axis2)
      !          a8(2,7) = -js*   rx2*dx214(axis2)

       curl1um1 =  is*8.*rsxy1(i1,i2,i3,axis1,1)*dr114(axis1)   ! coeff of u(-1) from v.x - u.y 
       curl1vm1 = -is*8.*rsxy1(i1,i2,i3,axis1,0)*dr114(axis1)   ! coeff of v(-1) from v.x - u.y 
       curl1um2 = -is*   rsxy1(i1,i2,i3,axis1,1)*dr114(axis1)   ! coeff of u(-2) from v.x - u.y 
       curl1vm2 =  is*   rsxy1(i1,i2,i3,axis1,0)*dr114(axis1)   ! coeff of v(-2) from v.x - u.y

       curl2um1 =  js*8.*rsxy2(j1,j2,j3,axis2,1)*dr214(axis2)   ! coeff of u(-1) from v.x - u.y 
       curl2vm1 = -js*8.*rsxy2(j1,j2,j3,axis2,0)*dr214(axis2)   ! coeff of v(-1) from v.x - u.y 
       curl2um2 = -js*   rsxy2(j1,j2,j3,axis2,1)*dr214(axis2)   ! coeff of u(-2) from v.x - u.y 
       curl2vm2 =  js*   rsxy2(j1,j2,j3,axis2,0)*dr214(axis2)   ! coeff of v(-2) from v.x - u.y

       aa8(2,0,0,nn) =  curl1um1/mu1
       aa8(2,1,0,nn) =  curl1vm1/mu1
       aa8(2,4,0,nn) =  curl1um2/mu1
       aa8(2,5,0,nn) =  curl1vm2/mu1

       aa8(2,2,0,nn) = -curl2um1/mu2  
       aa8(2,3,0,nn) = -curl2vm1/mu2    
       aa8(2,6,0,nn) = -curl2um2/mu2 
       aa8(2,7,0,nn) = -curl2vm2/mu2 


       ! -------------- Equation 3 -----------------------
       !   [ tau.{ (uv.xx+uv.yy)/eps -alphaP*P.tt } ] = 0
       !    P.tt = c4PttLEsum * L(E) + c4PttLLEsum* L^2(E) + ...
       ! coeff of P is not used, thus set to 0
       c4PttLEsum1 = 0.
       c4PttLLEsum1 = 0.
       c4PttLEsum2 = 0.
       c4PttLLEsum2 = 0.
       aa8(3,0,0,nn) = tau1*( aLap0*( 1./epsmu1 -alphaP1*c4PttLEsum1/epsmu1 ) - aLapSq0*alphaP1*c4PttLLEsum1/epsmu1**2 )
       aa8(3,1,0,nn) = tau2*( aLap0*( 1./epsmu1 -alphaP1*c4PttLEsum1/epsmu1 ) - aLapSq0*alphaP1*c4PttLLEsum1/epsmu1**2 )
       aa8(3,4,0,nn) = tau1*( aLap1*( 1./epsmu1 -alphaP1*c4PttLEsum1/epsmu1 ) - aLapSq1*alphaP1*c4PttLLEsum1/epsmu1**2 )
       aa8(3,5,0,nn) = tau2*( aLap1*( 1./epsmu1 -alphaP1*c4PttLEsum1/epsmu1 ) - aLapSq1*alphaP1*c4PttLLEsum1/epsmu1**2 )

       aa8(3,2,0,nn) =-tau1*( bLap0*( 1./epsmu2 -alphaP2*c4PttLEsum2/epsmu2 ) - bLapSq0*alphaP2*c4PttLLEsum2/epsmu2**2 )
       aa8(3,3,0,nn) =-tau2*( bLap0*( 1./epsmu2 -alphaP2*c4PttLEsum2/epsmu2 ) - bLapSq0*alphaP2*c4PttLLEsum2/epsmu2**2 )
       aa8(3,6,0,nn) =-tau1*( bLap1*( 1./epsmu2 -alphaP2*c4PttLEsum2/epsmu2 ) - bLapSq1*alphaP2*c4PttLLEsum2/epsmu2**2 )
       aa8(3,7,0,nn) =-tau2*( bLap1*( 1./epsmu2 -alphaP2*c4PttLEsum2/epsmu2 ) - bLapSq1*alphaP2*c4PttLLEsum2/epsmu2**2 )


      ! -------------- Equation 4 -----------------------
      !    [ (u.xx+u.yy).x + (v.xx+v.yy).y ] = 0

      setJacobian( aj1, axis1)

      ! dr1a(0:2) = dsBig in tangential directions if avoidInterfaceIterations=1
      dr0=dr1a(axis1)
      ds0=dr1a(axis1p1)

      aLapX0 = xLapCoeff4a(is,dr0,ds0)
      aLapX1 = xLapCoeff4b(is,dr0,ds0)

      bLapY0 = yLapCoeff4a(is,dr0,ds0)
      bLapY1 = yLapCoeff4b(is,dr0,ds0)

      setJacobian( aj2, axis2)

      ! dr2a(0:2) = dsBig in tangential directions if avoidInterfaceIterations=1
      dr0=dr2a(axis2)
      ds0=dr2a(axis2p1)

      cLapX0 = xLapCoeff4a(js,dr0,ds0)
      cLapX1 = xLapCoeff4b(js,dr0,ds0)

      dLapY0 = yLapCoeff4a(js,dr0,ds0)
      dLapY1 = yLapCoeff4b(js,dr0,ds0)


      ! 4  [ (u.xx+u.yy).x + (v.xx+v.yy).y ] = 0
  
      aa8(4,0,0,nn)= aLapX0*c1**2
      aa8(4,1,0,nn)= bLapY0*c1**2
      aa8(4,4,0,nn)= aLapX1*c1**2
      aa8(4,5,0,nn)= bLapY1*c1**2

      aa8(4,2,0,nn)=-cLapX0*c2**2
      aa8(4,3,0,nn)=-dLapY0*c2**2
      aa8(4,6,0,nn)=-cLapX1*c2**2
      aa8(4,7,0,nn)=-dLapY1*c2**2

      ! ---------------- Equation 5 (2nd-order) -----------------

      !   [ ( {(Delta v).x - (Delta u).y}/(epsmu) - alphaP*( Py.ttx - Px.tty) )/mu ] =0 
      !
      !     P.tt = c2PttLEsum * L(E)

      ! coeff of P is set to 0
      c2PttLEsum1 = 0.
      c2PttEsum1 = 0.
      c2PttLEsum2 = 0.
      c2PttEsum2 = 0.

      eqnCoeff = ( 1./epsmu1 - alphaP1*c2PttLEsum1/epsmu1 )/mu1 
      eqnCoeffb = -alphaP1*c2PttEsum1/mu1 ! added sept 16, 2018 
      aa8(5,0,0,nn)=-bLapY0*eqnCoeff + curl1um1*eqnCoeffb  
      aa8(5,1,0,nn)= aLapX0*eqnCoeff + curl1vm1*eqnCoeffb
      aa8(5,4,0,nn)=-bLapY1*eqnCoeff + curl1um2*eqnCoeffb 
      aa8(5,5,0,nn)= aLapX1*eqnCoeff + curl1vm2*eqnCoeffb

      eqnCoeff = ( 1./epsmu2 - alphaP2*c2PttLEsum2/epsmu2 )/mu2 
      eqnCoeffb = -alphaP2*c2PttEsum2/mu2 ! added sept 16, 2018 
      aa8(5,2,0,nn)=-(-dLapY0*eqnCoeff + curl2um1*eqnCoeffb)
      aa8(5,3,0,nn)=-( cLapX0*eqnCoeff + curl2vm1*eqnCoeffb)
      aa8(5,6,0,nn)=-(-dLapY1*eqnCoeff + curl2um2*eqnCoeffb)
      aa8(5,7,0,nn)=-( cLapX1*eqnCoeff + curl2vm2*eqnCoeffb)


       ! ------- Equation 6 -----
       !  [ nv.( c^2*Delta^2(E) - alphaP*Delta(Ptt) )/mu ] = 0 

       ! 6  [ n.Delta^2 u/eps ] = 0

       if( setDivergenceAtInterfaces.eq.0 )then
        ! use Eqn 6 
        ! NOTE: LE = c^2*Delta(E) and LLE = (c^4*Delta^2) E 
        ! Note: the coeff of L(E) in Delta(Ptt) is the coeff of E in Ptt
        ! Note: the coeff of LL(E) in Delta(Ptt) is the coeff of LE in Ptt

        c2PttEsum1 = 0.
        c2PttLEsum1 = 0.
        c2PttEsum2 = 0.
        c2PttLEsum2 = 0.

        aa8(6,0,0,nn) = an1*( aLapSq0/epsmu1 -alphaP1*( c2PttEsum1*aLap0 + c2PttLEsum1*aLapSq0/epsmu1 ) )/mu1
        aa8(6,1,0,nn) = an2*( aLapSq0/epsmu1 -alphaP1*( c2PttEsum1*aLap0 + c2PttLEsum1*aLapSq0/epsmu1 ) )/mu1
        aa8(6,4,0,nn) = an1*( aLapSq1/epsmu1 -alphaP1*( c2PttEsum1*aLap1 + c2PttLEsum1*aLapSq1/epsmu1 ) )/mu1
        aa8(6,5,0,nn) = an2*( aLapSq1/epsmu1 -alphaP1*( c2PttEsum1*aLap1 + c2PttLEsum1*aLapSq1/epsmu1 ) )/mu1

        aa8(6,2,0,nn) =-an1*( bLapSq0/epsmu2 -alphaP2*( c2PttEsum2*bLap0 + c2PttLEsum2*bLapSq0/epsmu2 ) )/mu2
        aa8(6,3,0,nn) =-an2*( bLapSq0/epsmu2 -alphaP2*( c2PttEsum2*bLap0 + c2PttLEsum2*bLapSq0/epsmu2 ) )/mu2
        aa8(6,6,0,nn) =-an1*( bLapSq1/epsmu2 -alphaP2*( c2PttEsum2*bLap1 + c2PttLEsum2*bLapSq1/epsmu2 ) )/mu2
        aa8(6,7,0,nn) =-an2*( bLapSq1/epsmu2 -alphaP2*( c2PttEsum2*bLap1 + c2PttLEsum2*bLapSq1/epsmu2 ) )/mu2
       end if

       ! ------- Equation 7 ------
       ! [ tv.( c^4*Delta^2(E) - alphaP*c^2*Delta(P.tt) - alphaP*P.tttt) ]=0 

       ! 7  [ tau.Delta^2 v/eps^2 ] = 0 
       ! Note: the coeff of L(E) in Delta(Ptt) is the coeff of E in Ptt
       ! Note: the coeff of LL(E) in Delta(Ptt) is the coeff of LE in Ptt

       c2PttEsum1 = 0.
       c2PttttLEsum1 = 0.
       c2PttLEsum1 = 0.
       c2PttttLLEsum1 = 0.
       c2PttEsum2 = 0.
       c2PttttLEsum2 = 0.
       c2PttLEsum2 = 0.
       c2PttttLLEsum2 = 0.

       coeffLap1   =              -alphaP1*(  c2PttEsum1 + c2PttttLEsum1  )/epsmu1
       coeffLapSq1 = 1./epsmu1**2 -alphaP1*( c2PttLEsum1 + c2PttttLLEsum1 )/epsmu1**2

       coeffLap2   =              -alphaP2*(  c2PttEsum2 + c2PttttLEsum2  )/epsmu2
       coeffLapSq2 = 1./epsmu2**2 -alphaP2*( c2PttLEsum2 + c2PttttLLEsum2 )/epsmu2**2

       aa8(7,0,0,nn) = tau1*( coeffLapSq1*aLapSq0 + coeffLap1*aLap0 )
       aa8(7,1,0,nn) = tau2*( coeffLapSq1*aLapSq0 + coeffLap1*aLap0 )
       aa8(7,4,0,nn) = tau1*( coeffLapSq1*aLapSq1 + coeffLap1*aLap1 )
       aa8(7,5,0,nn) = tau2*( coeffLapSq1*aLapSq1 + coeffLap1*aLap1 )

       aa8(7,2,0,nn) =-tau1*( coeffLapSq2*bLapSq0 + coeffLap2*bLap0 )
       aa8(7,3,0,nn) =-tau2*( coeffLapSq2*bLapSq0 + coeffLap2*bLap0 )
       aa8(7,6,0,nn) =-tau1*( coeffLapSq2*bLapSq1 + coeffLap2*bLap1 )
       aa8(7,7,0,nn) =-tau2*( coeffLapSq2*bLapSq1 + coeffLap2*bLap1 )



       ! save a copy of the matrix
       do n2=0,7
       do n1=0,7
         aa8(n1,n2,1,nn)=aa8(n1,n2,0,nn)
         ! a8(n1,n2)=aa8(n1,n2,0,nn)
       end do
       end do

     
       ! solve A Q = F
       ! factor the matrix
       numberOfEquations=8
       call dgeco( aa8(0,0,0,nn), numberOfEquations, numberOfEquations, ipvt8(0,nn),rcond,work(0))

     end if


     ! Save current solution to compare to new
     q(0) = u1(i1-is1,i2-is2,i3,ex)
     q(1) = u1(i1-is1,i2-is2,i3,ey)
     q(2) = u2(j1-js1,j2-js2,j3,ex)
     q(3) = u2(j1-js1,j2-js2,j3,ey)

     q(4) = u1(i1-2*is1,i2-2*is2,i3,ex)
     q(5) = u1(i1-2*is1,i2-2*is2,i3,ey)
     q(6) = u2(j1-2*js1,j2-2*js2,j3,ex)
     q(7) = u2(j1-2*js1,j2-2*js2,j3,ey)

     ! subtract off the contributions from the initial (wrong) values at the ghost points:
     do n=0,7
       f(n) = (aa8(n,0,1,nn)*q(0)+aa8(n,1,1,nn)*q(1)+aa8(n,2,1,nn)*q(2)+aa8(n,3,1,nn)*q(3)+\
               aa8(n,4,1,nn)*q(4)+aa8(n,5,1,nn)*q(5)+aa8(n,6,1,nn)*q(6)+aa8(n,7,1,nn)*q(7)) - f(n)
     end do

     ! solve A Q = F
     job=0
     numberOfEquations=8
     call dgesl( aa8(0,0,0,nn), numberOfEquations, numberOfEquations, ipvt8(0,nn), f(0), job)


     if( useJacobiUpdate.eq.0 )then
       u1(i1-is1,i2-is2,i3,ex)=(1.-omega)*u1(i1-is1,i2-is2,i3,ex) + omega*f(0)
       u1(i1-is1,i2-is2,i3,ey)=(1.-omega)*u1(i1-is1,i2-is2,i3,ey) + omega*f(1)
       u2(j1-js1,j2-js2,j3,ex)=(1.-omega)*u2(j1-js1,j2-js2,j3,ex) + omega*f(2)
       u2(j1-js1,j2-js2,j3,ey)=(1.-omega)*u2(j1-js1,j2-js2,j3,ey) + omega*f(3)

       u1(i1-2*is1,i2-2*is2,i3,ex)=(1.-omega)*u1(i1-2*is1,i2-2*is2,i3,ex) + omega*f(4)
       u1(i1-2*is1,i2-2*is2,i3,ey)=(1.-omega)*u1(i1-2*is1,i2-2*is2,i3,ey) + omega*f(5)
       u2(j1-2*js1,j2-2*js2,j3,ex)=(1.-omega)*u2(j1-2*js1,j2-2*js2,j3,ex) + omega*f(6)
       u2(j1-2*js1,j2-2*js2,j3,ey)=(1.-omega)*u2(j1-2*js1,j2-2*js2,j3,ey) + omega*f(7)
     else
       ! Jacobi-update
       wk1(i1-is1,i2-is2,i3,ex)=(1.-omega)*u1(i1-is1,i2-is2,i3,ex) + omega*f(0)
       wk1(i1-is1,i2-is2,i3,ey)=(1.-omega)*u1(i1-is1,i2-is2,i3,ey) + omega*f(1)
       wk2(j1-js1,j2-js2,j3,ex)=(1.-omega)*u2(j1-js1,j2-js2,j3,ex) + omega*f(2)
       wk2(j1-js1,j2-js2,j3,ey)=(1.-omega)*u2(j1-js1,j2-js2,j3,ey) + omega*f(3)

       wk1(i1-2*is1,i2-2*is2,i3,ex)=(1.-omega)*u1(i1-2*is1,i2-2*is2,i3,ex) + omega*f(4)
       wk1(i1-2*is1,i2-2*is2,i3,ey)=(1.-omega)*u1(i1-2*is1,i2-2*is2,i3,ey) + omega*f(5)
       wk2(j1-2*js1,j2-2*js2,j3,ex)=(1.-omega)*u2(j1-2*js1,j2-2*js2,j3,ex) + omega*f(6)
       wk2(j1-2*js1,j2-2*js2,j3,ey)=(1.-omega)*u2(j1-2*js1,j2-2*js2,j3,ey) + omega*f(7)
     end if
  

 endLoopsMask2d()
 ! =============== end loops =======================

 ! fill ghost lines using 4th order IC results

 if( useJacobiUpdate.ne.0 )then
   ! Jacobi-update: now fill in values 
   beginLoopsMask2d() 
     u1(i1-is1,i2-is2,i3,ex)=wk1(i1-is1,i2-is2,i3,ex)
     u1(i1-is1,i2-is2,i3,ey)=wk1(i1-is1,i2-is2,i3,ey)
     u2(j1-js1,j2-js2,j3,ex)=wk2(j1-js1,j2-js2,j3,ex)
     u2(j1-js1,j2-js2,j3,ey)=wk2(j1-js1,j2-js2,j3,ey)

     u1(i1-2*is1,i2-2*is2,i3,ex)=wk1(i1-2*is1,i2-2*is2,i3,ex)
     u1(i1-2*is1,i2-2*is2,i3,ey)=wk1(i1-2*is1,i2-2*is2,i3,ey)
     u2(j1-2*js1,j2-2*js2,j3,ex)=wk2(j1-2*js1,j2-2*js2,j3,ex)
     u2(j1-2*js1,j2-2*js2,j3,ey)=wk2(j1-2*js1,j2-2*js2,j3,ey)

     u1(i1-  is1,i2-  is2,i3,hz)=wk1(i1-  is1,i2-  is2,i3,hz)
     u2(j1-  js1,j2-  js2,j3,hz)=wk2(j1-  js1,j2-  js2,j3,hz)
     u1(i1-2*is1,i2-2*is2,i3,hz)=wk1(i1-2*is1,i2-2*is2,i3,hz)
     u2(j1-2*js1,j2-2*js2,j3,hz)=wk2(j1-2*js1,j2-2*js2,j3,hz)
   endLoopsMask2d()
 end if
      
 
#endMacro
\end{lstlisting}

\bibliographystyle{elsart-num}
\bibliography{bib/journal-ISI,bib/adegdm,bib/jwb,bib/jba,bib/henshaw,bib/henshawPapers,bib/DMX,bib/AVKrefs,bib/myrefs,bib/myrefs_qx}
 
\end{document}